\documentclass[12pt,a4paper]{article}
\usepackage[utf8x]{inputenc}
\usepackage[T1]{fontenc}
\usepackage[a4paper, margin=3cm]{geometry}

\usepackage[pdftex]{graphicx} 
\usepackage{calc} 
\usepackage{enumitem} 

\usepackage[all]{nowidow} 
\usepackage[protrusion=true,expansion=true]{microtype} 

\usepackage{lipsum} 
\usepackage{amsmath}
\usepackage{graphicx}
\usepackage{amssymb}
\usepackage{amsthm}
\usepackage{caption}
\usepackage{bbm}
\usepackage{dsfont}
\usepackage{accents}
\usepackage{mathrsfs}
\usepackage{enumitem}
\usepackage{tikz}
\usepackage{tikz-3dplot}
\usetikzlibrary{patterns}
\usetikzlibrary{arrows.meta}
\usepackage{subcaption}
\usepackage[dvipsnames]{xcolor}
\usetikzlibrary{shapes.misc}
\usetikzlibrary{calc}
\usetikzlibrary{positioning}
\usetikzlibrary{arrows.meta}
\usepackage{booktabs}
\usepackage{float}
\usepackage{todonotes}
\usepackage{mathtools}
\usepackage{enumitem}
\usepackage{mathtools}

\tdplotsetmaincoords{60}{120}

\usepackage[pdftex,linkcolor=black,pdfborder={0 0 0}]{hyperref} 
\numberwithin{equation}{section}

\newtheorem{Def}{Definition}[section]
\newtheorem{Theorem}[Def]{Theorem}
\newtheorem{Lemma}[Def]{Lemma}
\newtheorem{Prop}[Def]{Proposition}

\newtheorem{Ex}[Def]{Example}
\newtheorem{Remark}[Def]{Remark}

\newtheorem{Corollary}[Def]{Corollary}

\definecolor{cellblue}{HTML}{2C7FB8}
\definecolor{cellfill}{HTML}{A6CEE3}
\definecolor{boundarygreen}{HTML}{7FCDBB}
\definecolor{vertexgray}{HTML}{4D4D4D}
\definecolor{hiddenedge}{HTML}{969696}

\newcommand{\R}{\mathbb{R}}
\newcommand{\N}{\mathbb{N}}
\newcommand{\Z}{\mathbb{Z}}

\DeclareMathOperator{\E}{\mathbb{E}} 
 
\DeclareMathOperator{\supp}{\mathsf{supp}}
\DeclareMathOperator{\mylog}{\mathsf{log}}
\DeclareMathOperator{\myexp}{\mathsf{exp}}
\DeclareMathOperator{\Prob}{\mathbb{P}}
\DeclareMathOperator*{\argmin}{arg\,min}
\title{Slow and rapid mixing in $\Z_2$ hyperlattice gauge theory}
\author{Adrien Malacan \\
\small Mathematical Sciences, Chalmers University of Technology\\
\small \texttt{malacan@chalmers.se}}
\date{}

\begin{document} 
\maketitle 
\begin{abstract}
    The Ising model under Glauber dynamics on the torus is known to exhibit a phase transition in its mixing time. The mixing is exponentially slow when $\beta >\beta_c$ and logarithmically fast when $\beta <\beta_c$. In this paper, we generalize the slow-mixing result to $\Z_2$ hyperlattice gauge theory when $\beta$ is large enough. This is a higher-dimensional analogue of the Ising model, with spins on $k$-cells. Specifically, on the torus of any dimension $m \ge k+2$, we prove explicit exponential lower bounds on the mixing time of Glauber dynamics when $\beta$ is sufficiently large. More precisely, we prove that the result holds for a broader class of dynamics, allowing more than one $k$-cell to change its value at each step of the algorithm. Finally, we show that there are distinct slow and rapid mixing regimes, as for sufficiently small $\beta$ the Dobrushin condition is satisfied for Glauber dynamics under the heat-bath rule, which implies logarithmic mixing.
\end{abstract}
\section{Introduction}
In this section, we review the relevant literature, introduce the model and its associated dynamics, present our main results and outline the principal contributions of this paper.
\subsection{Literature}
Lattice gauge theory is a model from statistical mechanics, introduced in the 1970s independently by Wegner \cite{Wegner:1971app} and Wilson \cite{Wilson74}. Wilson introduced the model as a discretization of Yang-Mills theory; the rigorous construction of the corresponding continuum theory remains a major open problem in mathematical physics. Wegner constructed a model dual to the three-dimensional Ising model, now known as $\mathbb{Z}_2$ lattice gauge theory, which exhibits a phase transition without a local order parameter. He also introduced generalizations with spins on $k$-dimensional cells, known as hyperlattice gauge theory, which he called generalized Ising models. We refer to Figure \ref{fig:firstpicture} for a visualization of cells of dimensions $0,1,2$ and $3$.

Following the definition of lattice gauge theory as in \cite{Wilson74}, the topic was extensively studied in the physics literature in the following years, see, e.g., \cite{BalianDrouffeItzykson1974,FradkinSusskind1978,kogut1979introduction,Wilson74}. In recent years, an increasing amount of work from the probabilistic community has been done on these topics to obtain a better understanding of this model, see, e.g., \cite{Aizenman2025,Cao2020,Chatterjee2020WilsonLI,DuncanSchweinhart2025,DuncanSchweinhart2025Topological,LFVabelian,Forsstr_m_2023}. 

If we consider these models from a dynamical perspective, the Ising model (when $k=0$) provides a natural point of comparison. When considering continuous-time heat-bath Glauber dynamics on the $m$-dimensional torus, \cite{LubetzkySly2016} proved that when $\beta<\beta_c(m)$, where $\beta_c(m)$ is the critical value of the Ising model in $m$ dimensions, then the mixing time is logarithmic in the size of the torus, and that cutoff occurs. When $\beta >\beta_c(m)$,  phase coexistence leads to exponentially slow mixing, see \cite{CesiGuadagniMartinelliSchonmann1996} for the two-dimensional case with free boundary conditions and \cite{Bodineau2005,Pisztora1996} for estimates in higher dimensions. The mixing time at criticality was studied in the case $m=2$ in \cite{LubetzkySly2012}, where a polynomial upper bound was proved. Hence, these results show an interesting link between the dynamical properties and the model's equilibrium phase. We refer to \cite{RapidMixingGlauber,Martinelli1999,Weitz2005} for a discussion on the existence of a unique Gibbs state and fast mixing.

Concerning the mixing properties of three-dimensional lattice gauge theory, \cite{Stahl2026SlowMixing} outline a bottleneck argument showing exponentially slow mixing in dimension $m=3$, when considering single-update Glauber dynamics with the Metropolis rule, when $\beta>0$ is large enough. Their argument identifies regions of the state space that can only communicate through a set of configurations with exponentially small Gibbs probability (a so-called bottleneck). An important contribution of the paper is the extension of this argument when adding a weak Higgs field. We also mention the work of \cite{Thorngren2026GaugeAveraging} which proves stability under sufficiently small local perturbations of the Hamiltonian. 

These results also connect lattice gauge theory with the study of information storage. In quantum information theory, a canonical model is Kitaev's toric code \cite{Kitaevtoric}, a quantum error-correcting code closely related to lattice gauge theory. In this field, a central question is how long the encoded information can survive under thermal noise. For the toric code and its higher-dimensional generalizations, this question is related to the thermal properties of classical (non-quantum) $\Z_2$ lattice gauge models~\cite{MAZAC20122096}. 

In this paper, we will prove that the slow mixing argument at low temperature can be generalized to hyperlattice gauge theory, provided that the codimension of the model is at least two. We also explore classes of different dynamics. Moreover, if the updates are regular Glauber updates with the heat-bath rule, we prove logarithmic mixing at sufficiently high temperature.

\subsection{Model}
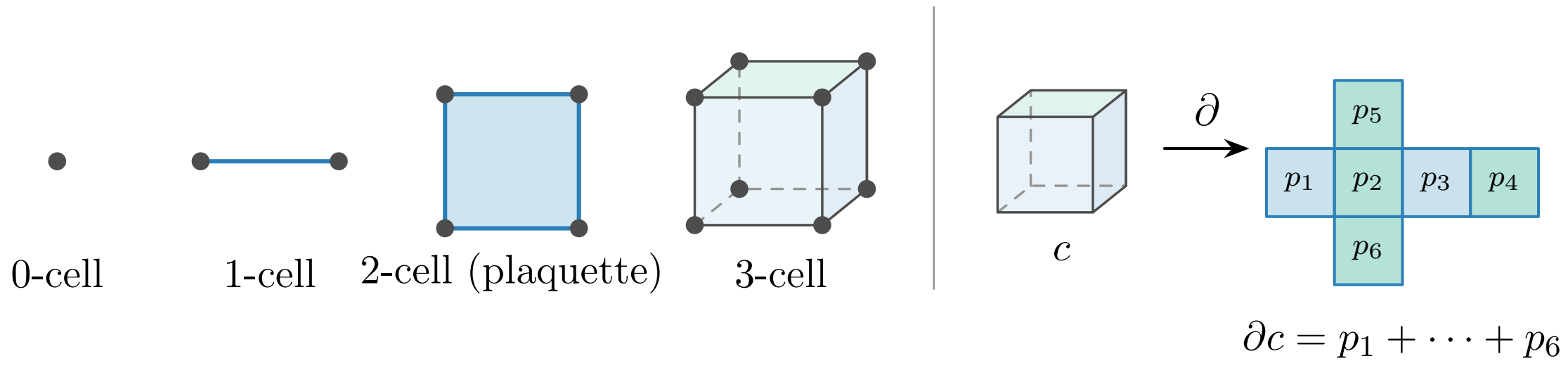
\begin{figure}
\begin{tikzpicture}[
  x=1cm,
  y=1cm,
  line cap=round,
  line join=round,
  cell edge/.style={draw=cellblue, line width=1.15pt},
  ordinary edge/.style={draw=vertexgray, line width=0.65pt},
  hidden edge/.style={draw=hiddenedge, line width=0.65pt, dashed},
  vertex/.style={circle, fill=vertexgray, inner sep=1.7pt},
  panel label/.style={font=\small, text=black},
  face label/.style={font=\scriptsize, text=black},
  >={Stealth[length=2.4mm]}
]

\node[vertex] at (0,0) {};
\node[panel label] at (0,-1.05) {$0$-cell};

\coordinate (e0) at (1.35,0);
\coordinate (e1) at (2.65,0);
\draw[cell edge] (e0)--(e1);
\node[vertex] at (e0) {};
\node[vertex] at (e1) {};
\node[panel label] at (2,-1.05) {$1$-cell};

\coordinate (p0) at (3.65,-0.63);
\coordinate (p1) at (4.91,-0.63);
\coordinate (p2) at (4.91, 0.63);
\coordinate (p3) at (3.65, 0.63);
\filldraw[fill=cellfill, fill opacity=0.55, cell edge]
  (p0)--(p1)--(p2)--(p3)--cycle;
\foreach \p in {p0,p1,p2,p3}{\node[vertex] at (\p) {};}
\node[panel label] at (4.28,-1.05) {$2$-cell (plaquette)};

\begin{scope}[shift={(6.00,-0.60)}]
  \coordinate (A)  at (0,0);
  \coordinate (B)  at (1.20,0);
  \coordinate (C)  at (1.20,1.20);
  \coordinate (D)  at (0,1.20);
  \coordinate (Ap) at (0.42,0.34);
  \coordinate (Bp) at (1.62,0.34);
  \coordinate (Cp) at (1.62,1.54);
  \coordinate (Dp) at (0.42,1.54);

  \fill[cellfill, opacity=0.24] (A)--(B)--(C)--(D)--cycle;
  \fill[cellblue, opacity=0.14] (B)--(Bp)--(Cp)--(C)--cycle;
  \fill[boundarygreen, opacity=0.22] (D)--(C)--(Cp)--(Dp)--cycle;

  \draw[hidden edge] (A)--(Ap)--(Bp);
  \draw[hidden edge] (Ap)--(Dp);
  \draw[ordinary edge] (A)--(B)--(C)--(D)--cycle;
  \draw[ordinary edge] (B)--(Bp)--(Cp)--(C);
  \draw[ordinary edge] (D)--(Dp)--(Cp);
  \draw[ordinary edge] (Bp)--(Cp);

  \foreach \p in {A,B,C,D,Ap,Bp,Cp,Dp}{\node[vertex] at (\p) {};}
  \node[panel label] at (0.81,-0.45) {$3$-cell};
\end{scope}

\draw[hiddenedge, line width=0.5pt] (8.25,-1.20)--(8.25,1.45);

\begin{scope}[shift={(8.85,-0.48)}, scale=0.78]
  \coordinate (A)  at (0,0);
  \coordinate (B)  at (1.15,0);
  \coordinate (C)  at (1.15,1.15);
  \coordinate (D)  at (0,1.15);
  \coordinate (Ap) at (0.40,0.32);
  \coordinate (Bp) at (1.55,0.32);
  \coordinate (Cp) at (1.55,1.47);
  \coordinate (Dp) at (0.40,1.47);

  \fill[cellfill, opacity=0.25] (A)--(B)--(C)--(D)--cycle;
  \fill[cellblue, opacity=0.16] (B)--(Bp)--(Cp)--(C)--cycle;
  \fill[boundarygreen, opacity=0.24] (D)--(C)--(Cp)--(Dp)--cycle;
  \draw[hidden edge] (A)--(Ap)--(Bp);
  \draw[hidden edge] (Ap)--(Dp);
  \draw[ordinary edge] (A)--(B)--(C)--(D)--cycle;
  \draw[ordinary edge] (B)--(Bp)--(Cp)--(C);
  \draw[ordinary edge] (D)--(Dp)--(Cp);
  \draw[ordinary edge] (Bp)--(Cp);
  \node[panel label] at (0.77,-0.47) {$c$};
\end{scope}

\draw[->, line width=0.75pt] (10.42,0.12)--(11.22,0.12)
  node[midway, above=2pt] {$\partial$};

\begin{scope}[shift={(11.38,-0.52)}]
  \def\s{0.64}

  \foreach \x/\name/\col in {
      0/p_1/cellfill,
      0.64/p_2/boundarygreen,
      1.28/p_3/cellfill,
      1.92/p_4/boundarygreen}{
    \filldraw[fill=\col, fill opacity=0.58, draw=cellblue, line width=0.75pt]
      (\x,0) rectangle ++(\s,\s);
    \node[face label] at ($(\x,0)+(0.5*\s,0.5*\s)$) {$\name$};
  }

  \filldraw[fill=boundarygreen, fill opacity=0.58,
            draw=cellblue, line width=0.75pt]
    (0.64,0.64) rectangle ++(\s,\s);
  \node[face label] at (0.96,0.96) {$p_5$};

  \filldraw[fill=boundarygreen, fill opacity=0.58,
            draw=cellblue, line width=0.75pt]
    (0.64,-0.64) rectangle ++(\s,\s);
  \node[face label] at (0.96,-0.32) {$p_6$};

  \node[panel label, anchor=north] at (1.28,-0.84)
    {$\partial c=p_1+\cdots+p_6$};
\end{scope}

\end{tikzpicture}
\caption{On the left-hand side we see $0,1,2,3$-cells. The right-hand side pictures the boundary of a $3$-cell denoted $c$, showing the 6 plaquettes in its boundary.}
\label{fig:firstpicture}

\end{figure}

We now introduce our model. Let $N \in \N_{\ge 3}$ and $m \in \N_{\ge 1}$. We let $\Z_N^m$ denote the discrete $m$-dimensional torus that contains $N$ vertices in each direction. In $\Z_N^m$, for $k \in \{0,\cdots,m-1\}$, we  call a $k$-cell the generalization of sites or edges to dimension~$k$. For example a $0$-cell is a site, a $1$-cell an edge, a $2$-cell a plaquette and a $3$-cell a cube. We refer to Figure \ref{fig:firstpicture} for a visualization of those and refer to Section \ref{sec:dec} for a more precise mathematical definition. 
We denote the set of $k$-cells in $\Z_N^m$ by $C_k(\Z_N^m)$. As we will work in the group $\Z_2$, the orientation of the cells doesn't play a role.
A $k$-form $\tau$ is a function assigning spins in $\Z_2$ to each $k$-cell $\tau : C_k(\Z_N^m) \to \Z_2$.
The set of $k$-forms is denoted by $C^k(\Z_N^m,\Z_2)$. Hyperlattice gauge theory is defined as a Gibbs probability measure on the set of $k$-forms. To introduce this measure, we first introduce the \textit{discrete exterior operator}, denoted $d_k$, which is a mapping from $C^k(\Z_N^m,\Z_2)$ to $C^{k+1}(\Z_N^m,\Z_2)$. For a $(k+1)$-cell $c$, we denote by $\partial_{k+1} c$ the formal sum of the $k$-cells in its boundary. In the following, we abuse notation and identify $\partial_k c$ with its support.
With this notation, we let
\[
d_k\tau(c) \coloneqq \sum_{a \in \partial_{k+1} c}\tau(a) \in \Z_2, \quad c\in C_{k+1}(\Z_N^m).
\]
In particular, we have $d_k\tau(c)=1$ iff there is an odd number of $k$-cells $a$ in the boundary of $c$ such that $\tau(a)=1$.  Let furthermore $\rho$ denote the unitary representation of $\Z_2$ such that $\rho(x)=(-1)^x, ~ x \in \Z_2$. We now define the Hamiltonian
of hyperlattice gauge theory as
\[
\mathcal{H}_{N,m,k}(\tau)\coloneqq - \sum_{c \in C_{k+1}(\Z_N^m)} \rho(d_k\tau(c)), ~\text{for} ~ \tau \in C^{k}(\Z_N^m,\Z_2).
\]
With this notation one defines hyperlattice gauge theory for $k$-cells, at inverse temperature $\beta>0$, as the probability measure
\begin{equation}
\label{eq:probmass}
\mu_{\beta,N,m,k}(\tau)=\frac{1}{Z_{\beta,N,m,k}}\myexp(-\beta\mathcal{H}_{N,m,k}(\tau)), \quad \tau \in C^{k}(\Z_N^m,\Z_2),
\end{equation}
where $Z_{\beta,N,m,k}$ is the partition function.  
Note that the case $k=0$ corresponds to the Ising model, and the case $k=1$ to $\Z_2$ lattice gauge theory also called Ising lattice gauge theory. Throughout this paper we might use the abbreviation $\mu=\mu_{\beta,N,m,k}$.\newline 
\subsection{Block dynamics}
\label{subsec:gd}
In this section, we introduce a family of Glauber dynamics, which generalize the \textit{heat-bath} dynamics, as well as the \textit{Metropolis-Hastings} algorithm. We first introduce the discrete-time version of this dynamics, before presenting the continuous-time version below.

To this end let $\mathcal{B}\coloneqq\{B_1,\dots,B_n\}$ be a subdivision of $C_k(\Z_N^m)$ into \textit{blocks}, such that for every $B \in \mathcal{B}$ one has $B \subseteq C_k(\Z_N^m)$ and $\bigcup_{B \in \mathcal{B}}B=C_k(\Z_N^m)$. Note that the blocks are not necessarily disjoint. At each discrete time step $n \in \N$, a block $B \in \mathcal{B}$ is selected according to the probabilities $(p_B)_{B \in \mathcal{B}}$.
Once a block $B$ is selected, given the state at time $n$ denoted $\tau_n$, we select the new proposal state $\tau'$ for $\tau_{n+1}$ according to a transition kernel $Q_B:C^k(\Z_N^m,\Z_2) \times C^k(\Z_N^m,\Z_2) \to [0,1]$, satisfying
\begin{equation}
\label{eq:condonkernel}
Q_B(\tau_n,\tau')=0 ~ \text{if} ~ \tau_n|_{B^c} \neq \tau'|_{B^c},
\end{equation}
i.e., we only accept changing the value of $\tau_n$ in the set of $k$-cells in $B \subset C_k(\Z_N^m)$. 
Hence the overall transition kernel has the following form
\begin{equation}
\label{eq:defQproposal}
Q(\tau_n,\tau')=\sum_{B \in \mathcal{B}} p_BQ_B(\tau_n,\tau').
\end{equation}
Note that the transformation $\tau_n \to \tau'$ might be quite general, and that $Q_B$ is itself not necessarily assumed to be reversible. 
Given the current state $\tau_n$ and the new state $\tau'$ proposed by $Q$, we denote by $\alpha(\tau_n,\tau') \in [0,1]$ the acceptance probability of this new state, meaning that with probability $\alpha(\tau_n,\tau')$ one has $\tau_{n+1}=\tau'$, and otherwise $\tau_{n+1}=\tau_n$. 
Now, from \cite{Generalmetropolis}, one has the following result, where $\mu$ is the measure from \eqref{eq:probmass}.
\begin{Lemma}[Theorem 2 in \cite{Generalmetropolis}]
\label{lem:firstlemmaall}
    If the target measure $\mu$, the transition kernel $Q$ and the acceptance probability $\alpha(\tau,\tau')$ satisfy the detailed balance condition for any $\tau,\tau' \in C^k(\Z_N^m,\Z_2)$
    \begin{equation}
    \label{eq:balancecond}
    \mu(\tau)Q(\tau,\tau')\alpha(\tau,\tau')=\mu(\tau')Q(\tau',\tau)\alpha(\tau',\tau),
    \end{equation}
    then the transition kernel $P$ defined by
    \begin{equation}
    \label{eq:defP}
    P(\tau,\tau')= Q(\tau,\tau')\alpha(\tau,\tau')+\mathbbm{1}\{\tau=\tau'\} \sum_{\tilde{\tau} \in C^k}(1-\alpha(\tau,\tilde{\tau}))Q(\tau,\tilde{\tau})
    \end{equation}
    is reversible with respect to the measure $\mu$. If this kernel $P$ is irreducible with respect to $\mu$ and aperiodic, then
   $(\tau_n)_{n \in \N}$ which evolves according to this transition kernel converges in law to $\mu$.
\end{Lemma}
We now construct continuous-time dynamics as follows. Equip $\mathcal{B}$ with an exponential clock of parameter $|\mathcal{B}|$, i.e., the number of blocks. Once the clock rings at time $t > 0$, denote $\tau_t^- \coloneqq \lim_{s \uparrow t}\tau_s$ and apply the same procedure as in the discrete case: one applies the kernel $P$ from \eqref{eq:defP} starting from $\tau_t^-$ to obtain a new state $\tau_t$, and one starts the clock again.
For any irreducible transition kernel $P$ as described in \eqref{eq:defP} such that condition \eqref{eq:balancecond} is satisfied, if $P$ is irreducible, then the law of the continuous dynamics $(\tau_t)_{t \ge 0}$ converges in distribution to $\mu$ by \cite[Theorem 20.1]{levin2017markov}.
Note that the generator of the Markov process $(\tau_t)_{t \ge 0}$ equals
\begin{equation}
\label{eq:defgenL1}
    \mathcal{L}=|\mathcal{B}|(P-I).
\end{equation}
This is due to the following. Let $f:C^k(\Z_N^m,\Z_2) \to \R$ be any function. Then the infinitesimal changes on the value of $f(\tau), \tau \in C^k$ are given by
    \begin{equation}
    \label{eq:generatorL}
    \mathcal{L}f(\tau)=|\mathcal{B}|\Big(\sum_{\tau' \in C^k} f(\tau')P(\tau,\tau')-f(\tau)\Big)=|\mathcal{B}|(P-I)f(\tau).
    \end{equation}
    Here the factor $|\mathcal{B}|$ comes from the rate at which the exponential clock rings.
For a collection of blocks $\mathcal{B}$, we define the size of its largest block as
\begin{equation}
\label{eq:deflambdab}
    \lambda(\mathcal{B}) \coloneqq \max_{B \in \mathcal{B}}|B|.
\end{equation}
The Markovian dynamics induced by $\mathcal{L}$ includes many different dynamics, as the following example shows. In the example, for a $k$-cell $c \in C_k(\Z_N^m)$, we denote by $\delta_c$ the $k$-form with support exactly $\{c\}$. For two $k$-forms $\tau, \tau'$ we denote by $\tau+\tau'$ their cell-by-cell sum. We now present three examples of irreducible dynamics, which are special cases of our construction above.
\begin{Ex}
\label{ex:firstexample}
\begin{enumerate}[label=(\alph*)]
    \item  Take each block to be a single $k$-cell, i.e., let the blocks be
    $\mathcal{B}=(B_c)_{c \in C_k(\Z_N^m)}$ where $B_c=\{c\}$. Furthermore, let 
    \[\alpha(\tau,\tau')\coloneqq\min(1, r(\tau',\tau) ),
\quad \text{where} ~~ r(\tau,\tau') \coloneqq \frac{\mu(\tau)Q(\tau,\tau')}{\mu(\tau')Q(\tau',\tau)},
\]
let $(p_c)_{c \in C_k(\Z_N^m)}$ be uniform and let $Q_c$ propose $\tau'=\tau+\delta_c$ meaning that one flips the value of $\tau$ on the $k$-cell $c$. Then, this is the well-studied Metropolis algorithm under Glauber updates.
\item Let $\mathcal{B}$ and $(p_c)_{c \in C_k(\Z_N^m)}$ be exactly as in (a). Define the transition proposals
\[
Q_c(\tau,\tau+\delta_c)=\frac{\mu(\tau+\delta_c)}{\mu(\tau)+\mu(\tau+\delta_c)}, ~ Q_c(\tau,\tau)=1-Q_c(\tau,\tau+\delta_c).
\]
Then, for $\alpha=1$, the condition \eqref{eq:balancecond} holds. Moreover, these are exactly the transition probabilities of the heat-bath dynamics under Glauber updates.
\item Choose any blocks $\mathcal{B}$. Take any probability distribution $(p_B)_{B \in \mathcal{B}}$ such that each $p_B >0$. Then, define
\[
Q_B(\tau,\cdot)=\mu(\cdot|\tau_{B^c}), ~ \tau \in C^k,
\]
i.e., one samples the new configuration according to the measure $\mu$ conditioned on being equal to the configuration $\tau$ outside of the cells in $B$. Then, exactly as in (b), \eqref{eq:balancecond} is satisfied for $\alpha=1$. In fact, this construction is a generalization of Glauber dynamics using the heat-bath rule on blocks. This construction appears in~\cite{BlancaEtAl2022,Martinelli1999}.
\end{enumerate}
\end{Ex}
Note that in this work, we allow the blocks to depend on the torus $\Z_N^m$. Thus the sequence will be denoted $(\mathcal{B}_N)_{N \in \N}$, such that $\mathcal{B}_N$ is a collection of blocks in $C_k(\Z_N^m)$. In our work, an important quantity will be $\lambda(\mathcal{B}_N)$ which upper bounds the largest Hamming distance on $k$-cells that can be achieved in one step by our dynamics. This value will be crucial for the proof of the main theorem. Concerning the blocks that are considered, we will assume the following upper bound: there exists $a>0$ such that
\begin{equation}
\label{eq:condonB_N}
|\mathcal{B}_N|<N^a.
\end{equation}
This assumption makes sense from a practical point of view, as the number of cells grows as a polynomial in $N$.

\begin{figure}
\centering
\begin{tikzpicture}[
  x=0.78cm, y=0.78cm,
  font=\small, line cap=round, line join=round,
  lgt/grid/.style={draw=black!35, line width=0.3pt},
  lgt/spin/.style={draw=red!75!black, line width=0.75pt},
  lgt/defect/.style={fill=orange!23, draw=none},
  lgt/selected edge/.style={draw=violet!85!black, line width=0.6pt,
    dash pattern=on 1.5pt off 1.1pt, rounded corners=1.2pt},
  lgt/selected vertex/.style={fill=violet!85!black, draw=violet!85!black,
    line width=0.4pt},
  lgt/arrow/.style={-{Stealth[length=2.5mm]}, line width=0.7pt},
  lgt/title/.style={anchor=west, font=\normalsize},
  lgt/celllabel/.style={font=\scriptsize}
]

\def\LGTGrid{%
  \draw[lgt/grid, step=1] (0,0) grid (4,4);
  \foreach \x in {0,...,4}{%
    \foreach \y in {0,...,4}{%
      \fill[black!65] (\x,\y) circle[radius=0.7pt];
    }
  }
}
\def\LGTSelectedEdge{%
  \draw[lgt/selected edge] (1.89,1.90) rectangle (2.11,3.10);
}
\def\LGTInitialDefects{%
  \fill[lgt/defect] (1,1) rectangle (2,2);
  \fill[lgt/defect] (1,2) rectangle (2,3);
}
\def\LGTSelectedVertex{%
  \filldraw[lgt/selected vertex]
    (2,2.105) -- (2.105,2) -- (2,1.895) -- (1.895,2) -- cycle;
  \node[lgt/celllabel, text=violet!85!black] at (2.24,2.23) {$v$};
}

\node[lgt/title] at (0,5.25)
  {\textbf{(a)} Update on a $k$-cell ($k=1$)};
\node at (2,4.43) {$\tau$};
\node at (12,4.43) {$\tau'$};

\begin{scope}
  \LGTInitialDefects
  \LGTGrid
  \draw[lgt/spin] (1,2) -- (2,2);
  \LGTSelectedEdge
  \node[lgt/celllabel, text=violet!85!black] at (2.23,2.86) {$c$};
  \node[lgt/celllabel] at (1.5,2.5) {$p_1$};
  \node[lgt/celllabel] at (2.5,2.5) {$p_2$};
\end{scope}

\begin{scope}[xshift=7.8cm]
  \fill[lgt/defect] (1,1) rectangle (2,2);
  \fill[lgt/defect] (2,2) rectangle (3,3);
  \LGTGrid
  \draw[lgt/spin] (1,2) -- (2,2) -- (2,3);
  \LGTSelectedEdge
  \node[lgt/celllabel, text=violet!85!black] at (2.23,2.86) {$c$};
  \node[lgt/celllabel] at (1.5,2.5) {$p_1$};
  \node[lgt/celllabel] at (2.5,2.5) {$p_2$};
\end{scope}

\node at (7,3.45) {$c\in B$};
\node at (7,2.70) {$\tau'=\tau+\delta_c$};
\draw[lgt/arrow] (4.55,2) -- (9.45,2);
\node at (7,1.16) {$d_1\tau'=d_1\tau+d_1\delta_c$};
\node[font=\footnotesize] at (7,0.48)
  {$p_1:\ 1\to0,\qquad p_2:\ 0\to1$};

\begin{scope}[yshift=-5.85cm]
  \node[lgt/title] at (0,5.25)
    {\textbf{(b)} Update generated by a $(k-1)$-cell};
  \node at (2,4.43) {$\tau$};
  \node at (12,4.43) {$\tau'$};

  \begin{scope}
    \LGTInitialDefects
    \LGTGrid
    \draw[lgt/spin] (1,2) -- (2,2);
    \LGTSelectedVertex
  \end{scope}

  \begin{scope}[xshift=7.8cm]
    \LGTInitialDefects
    \LGTGrid
    \draw[lgt/spin] (2,1) -- (2,3);
    \draw[lgt/spin] (2,2) -- (3,2);
    \LGTSelectedVertex
  \end{scope}

  \node at (7,3.45) {$S=\{v\}\subseteq \tilde{B}$};
  \node at (7,2.70) {$\tau'=\tau+d_0\delta_v$};
  \draw[lgt/arrow] (4.55,2) -- (9.45,2);
  \node at (7,1.16) {$d_1\tau'=d_1\tau$};
  \node[font=\footnotesize] at (7,0.48) {all plaquettes unchanged};
\end{scope}

\begin{scope}[yshift=-6.70cm, font=\footnotesize]
  \draw[lgt/spin] (0,0) -- (0.6,0);
  \node[anchor=west] at (0.80,0) {edge value $1$};
  \fill[lgt/defect] (4.50,-0.18) rectangle (4.94,0.18);
  \draw[lgt/grid] (4.50,-0.18) rectangle (4.94,0.18);
  \node[anchor=west] at (5.14,0) {plaquette value $1$};
  \draw[lgt/grid] (10.15,0) -- (10.70,0);
  \draw[lgt/selected edge] (10.10,-0.11) rectangle (10.75,0.11);
  \filldraw[lgt/selected vertex]
    (11.05,0.105) -- (11.155,0) -- (11.05,-0.105) -- (10.945,0) -- cycle;
  \node[anchor=west] at (11.35,0) {selected cell};
\end{scope}

\end{tikzpicture}
\caption{Two updates for $m=2$ and $k=1$, starting from the same configuration. \newline 
(a) Flipping an edge $c$ in $B$ flips the value of $\tau$ the two adjacent plaquettes $p_1$ and $p_2$.\newline 
  (b) The update generated by $S=\{v\}$ contained in $\tilde{B}$ flips the four edges incident to $v$ and leaves every plaquette unchanged.}
  \label{fig:dynamics}
\end{figure}
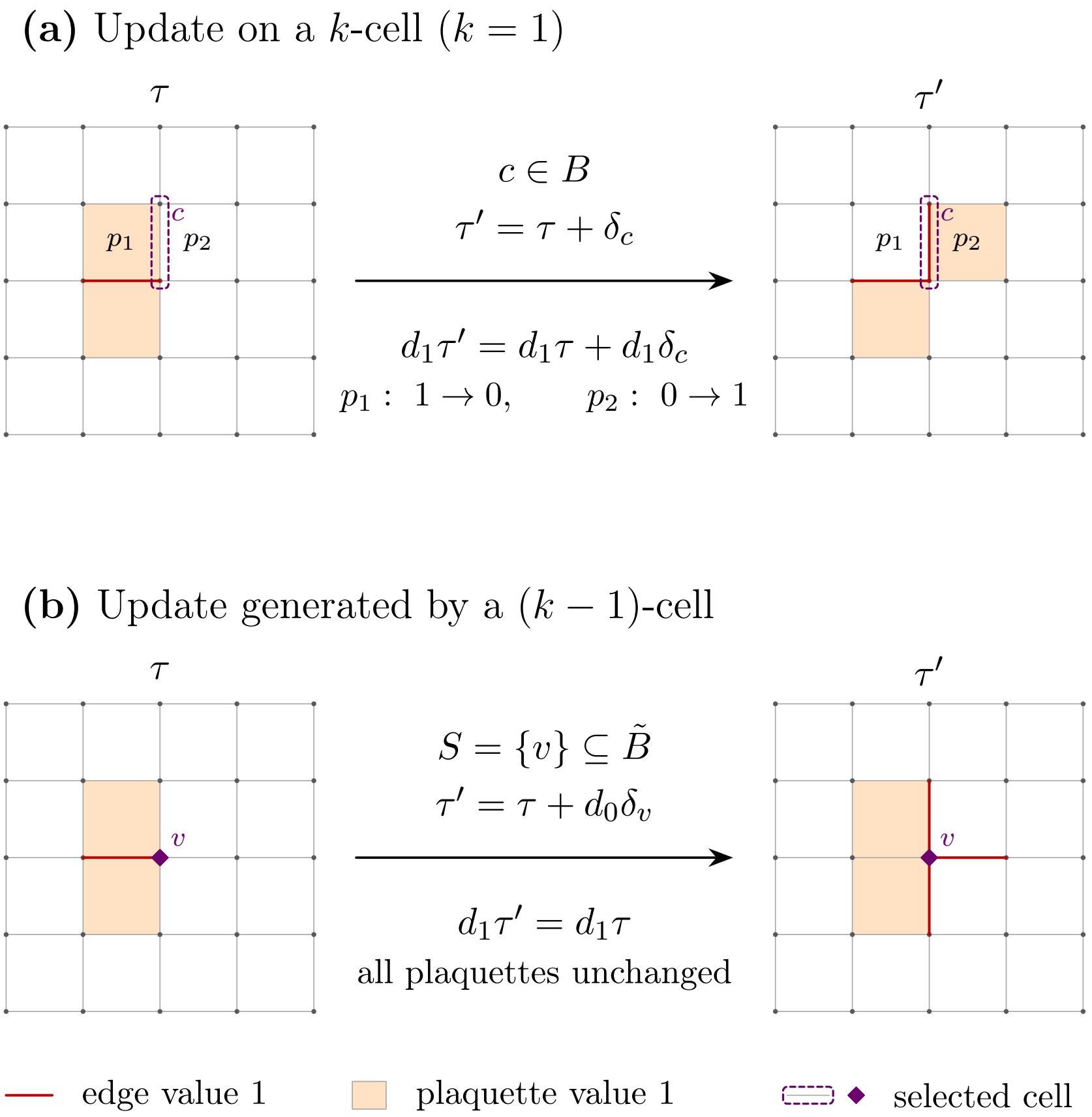

When thinking intuitively about adding dynamics to lattice gauge theory, one could consider two types of single-site Glauber updates. Either one flips the value of the $1$-form on an edge $e \in C_1(\Z_N^m)$, such that the updated configuration is $\tau'=\tau+\delta_e$, or one could consider making the model \textit{noisier}, in which case one could perform a so-called gauge transform and propose the new state as $\tau'=\tau+d\delta_v$, for $v \in C_0(\Z_N^m)$. This can be observed in Figure \ref{fig:dynamics}. Combining these two dynamics would intuitively increase the convergence speed. This can be generalized to any $k \ge 1$, where one can either act on the value on $k$-cells, or on the value of the discrete derivative of $(k-1)$-cells. Note that these dynamics are included in our example: one can choose some blocks to be the supports of discrete derivative of $(k-1)$-forms. Namely for any $D \subseteq C_{k-1}(\Z_N^m)$ pick 
\[
B_D \coloneqq \supp(d_{k-1}\sum_{c \in D} \delta_c)
\]
where $\delta_c$ is the $(k-1)$-form having support exactly at $c \in C_{k-1}(\Z_N^m)$. 
Note that these dynamics do not change the discrete derivative $d_j\tau$ and hence do not change energy. Indeed let $c \in C_{k-1}(\Z_N^m)$ and $\tau \in C^k$, then 
\[
d_k(\tau+d_{k-1}\delta_c)=d_k\tau+d_kd_{k-1}\delta_c=d_k\tau,
\]
as $d_kd_{k-1}=0$, we refer to Section \ref{sec:dec} for more details. One obtains $\mathcal{H}(\tau)=\mathcal{H}(\tau')$. We will see in the last section that these types of blocks do not play a role in the mixing time argument when $\beta$ is large. 

\begin{Def}
\label{def:markovprocess}
   If $\mathcal{L}$ is a generator as in \eqref{eq:defgenL1}, if the transition kernel $P$ has the form \eqref{eq:defP}, satisfies the condition of Lemma \ref{lem:firstlemmaall}, and is irreducible, then we say that the Markov process $(\tau_t)_{t \ge 0}$ with generator $\mathcal{L}$ satisfies condition $(A)$.
\end{Def}

\subsection{Mixing times and bottleneck ratio}
While the existence of a stationary measure towards which the law of a Markov chain converges is a classic result, the time needed to converge, called the \textit{mixing time}, is usually unknown, depends on the Markov process in question and is more complicated to analyze. This topic has attracted considerable attention from the probabilistic community in recent years, see e.g., \cite{levin2017markov,salez2025modernaspectsmarkovchains}. \newline 
We define the \textit{total variation distance} between two measures $\mu,\nu$ on a finite state space $S$, by
\[
\|\mu-\nu\|_{TV}\coloneqq\sup_{A \subseteq S}|\mu(A)-\nu(A)|.
\]
For a Markov process $(X_t)_{t \ge 0}$, we denote by $\Prob_{x}(X_t \in \cdot)$ the  law of the process at time $t \ge 0$ given that the process started in the state $x \in S$ at time $t=0$. We define the \textit{mixing time} with threshold $\epsilon \in (0,1)$ as the first time at which the total variation distance between the law of the chain and the stationary distribution $\mu$ is less than $\epsilon$ from every initial state, i.e.,
\[
t_{\text{mix}}(\epsilon)\coloneqq\inf \{t \ge 0 \colon \sup_{x \in S} \|\Prob_{x}(X_t \in \cdot)-\mu\|_{TV} < \epsilon\}.
\] 
The \textit{bottleneck ratio} is a tool to lower bound the mixing time. To describe this tool, let $P(x,\cdot)$ denote the transition probability in discrete time of the Markov process from the point $x$. Define the flow from stationarity for $x,y \in S$ as
\begin{equation}
\label{eq:edgemeasure}
\tilde{P}(x,y)\coloneqq \mu(x)P(x,y),
\end{equation}
and let the bottleneck ratio of a set $A \subseteq S$ such that $\mu(A) >0$ be given by
\[
\Phi(A)\coloneqq\frac{\tilde{P}(A,A^c)}{\mu(A)}.
\]
Finally, we let the bottleneck ratio of the entire Markov chain be defined as
\begin{equation}
\label{eq:bottleneckratio}
\Phi^\star \coloneqq \inf_{\substack{A \subset S\\\mu(A)\le \frac{1}{2}}} \Phi(A).
\end{equation}
Intuitively, a small bottleneck ratio means that relative to the probability of a set~$A$, the probability that the Markov chain leaves this set is low.
There is the following relation between the mixing time and the bottleneck ratio.
\begin{Theorem}[Theorem~7.4 and its proof in \cite{levin2017markov}]
\label{thm:bottleneckratiotomix}
Let $\epsilon \in (0,\frac{1}{2})$. Then the mixing time for the Markov chain $(X_n)_{n \in \N}$ with bottleneck ratio $\Phi^\star$ satisfies
\[
t_{\text{mix}}(\epsilon) \ge \frac{\frac{1}{2}-\epsilon}{\Phi^\star}.
\]
\end{Theorem}
\begin{Corollary}
    \label{cor:continuousdynamics}
    Let $\epsilon \in (0,\frac{1}{2})$, and consider the Markov process $(\tau_t)_{t \ge 0}$ with generator $\mathcal{L}$ defined in \eqref{eq:defgenL1} using the collection of blocks $\mathcal{B}$. Then, if the bottleneck ratio of the discrete dynamics $P$ defined in \eqref{eq:defP} is $\Phi^\star_P$, the mixing time of $(\tau_t)_{t \ge 0}$ satisfies
    \[
    t_{\text{mix}}(\epsilon) \ge \frac{\frac{1}{2}-\epsilon}{|\mathcal{B}|\Phi^\star_P}.
    \]
\end{Corollary}
\begin{proof}
    Consider the generator of the process $\mathcal{L}$ driven by the discrete dynamics $P$. 
    Computing the bottleneck ratio for the generator $\mathcal{L}$ for any $A \subseteq C^k$ we obtain
    \begin{align*}
    \Phi_{\mathcal{L}}(A)&=\frac{\sum_{\tau \in A, \tau' \in A^c}\mu(\tau)\mathcal{L}(\tau,\tau')}{\mu(A)} =|\mathcal{B}|\frac{1}{\mu(A)}\sum_{\tau \in A, \tau' \in A^c} \mu(\tau)(P(\tau,\tau')-I(\tau,\tau')) \\
    &= |\mathcal{B}|\frac{1}{\mu(A)}\sum_{\tau \in A, \tau' \in A^c} \mu(\tau)P(\tau,\tau')
    = |\mathcal{B}| \frac{\tilde{P}(A,A^c)}{\mu(A)}= |\mathcal{B}| \Phi_P(A),
    \end{align*}
    as the term $I(\tau,\tau')$ vanishes as $\tau \neq \tau'$ by definition. This implies that $\Phi^\star_{\mathcal{L}}=|\mathcal{B}|\Phi^\star_{P}$, and hence the result follows from Theorem \ref{thm:bottleneckratiotomix}.
\end{proof}
In the rest of the paper, we will consider a family of Markov processes which are indexed by the size of the torus $N \in \mathbb{N}$, and we write $t_{\text{mix}}^N(\epsilon)$ for the mixing time of the Markov process on $k$-forms in $C^k(\Z_N^m,\Z_2)$.
\subsection{Main Results}
We now present our main result, which proves exponentially slow mixing using a bottleneck argument at low temperature. We denote by $\ell \coloneqq m-k$ the codimension of the model.
\begin{Theorem}
\label{thm:lowerboundlgt}
    Consider $\mu_{\beta, N,m,k}$ the measure as defined in \eqref{eq:probmass} and let $k \in \{0, \dots, m-2\}$. Assume that the Markov process $(\tau_t)_{t \ge 0}$ satisfies Condition (A) stated in Definition \ref{def:markovprocess}.
    Then, if 
    \begin{equation}
    \label{eq:condbeta}
    \beta > \frac{\mylog\big((2(2k+3)(\ell-1))-1\big)+1}{2},
    \end{equation}
    and if the size of the largest block $\lambda(\mathcal{B}_N)$ in $(\mathcal{B}_N)_{N \in \N}$ has the following growth
    \begin{equation}
    \label{eq:conditionblock}
    \lambda(\mathcal{B}_N)= o\Big( \frac{N^{\ell}}{\log(N)^{\frac{\ell}{\ell-1}}} \Big),
    \end{equation}
    there exist $C,c>0$ depending only on $k,m,\beta$, and the entire sequence  $(\mathcal{B}_N)_{N \in \N}$, such that 
    \[
    t^N_{\text{mix}}(\epsilon) \ge \frac{C}{|\mathcal{B}_N|}\Big(\frac{1}{2}-\epsilon\Big)\myexp\Big(c\Big(\frac{N}{\lambda(\mathcal{B}_N)^\frac{1}{\ell}}\Big)^{\ell-1}\Big), ~~ \epsilon \in (0,\frac{1}{2}).
    \]
    In particular, if the size of the largest block satisfies $\lambda(\mathcal{B}_N)=\mathcal{O}(N^{\ell-\delta})$ for some $\delta >0$, then, the mixing time is exponential in the size of the torus, as $|\mathcal{B}_N|$ satisfies \eqref{eq:condonB_N}.
\end{Theorem}
In the case $k=1$, one recovers Ising lattice gauge theory, so that one has the following result as a special case. We give the result when considering the dynamics from Example \ref{ex:firstexample}(a) and (b).
\begin{Corollary}
    Let $\mu_{\beta,N,m,1}$ be the probability measure of Ising lattice gauge theory in any dimension $m \ge 3$ on the torus $\Z_N^m$. If the inverse temperature satisfies
    \[
    \beta > \frac{1}{2}\Big(\log\big(10(m-2)-1\big)+1\Big),
    \]
    there exist $C,c>0$, depending only on $\beta,m$ such that the mixing time under Glauber dynamics using the heat-bath or Metropolis rule on edges satisfies the following lower bound for $\epsilon \in (0,\frac{1}{2})$
    \[
    t_{\text{mix}}^N(\epsilon) \ge \frac{C}{N^m}\Big(\frac{1}{2}-\epsilon\Big)\exp\Big(cN^{m-2} \Big).
    \]
    Observe that we have the same lower bound if we allow gauge flips in our dynamics.
\end{Corollary}
\begin{Remark}
    In the case $k=0$ one recovers known results from the Ising model. For $m \ge 2$, we obtain that if 
    \[
    \beta > \frac{\mylog(6(m-1)-1)+1}{2}>\beta_c(m),
    \]
    then the mixing time on the torus under single-site Glauber updates using the Metropolis or heat-bath rule satisfies
    \[
    t_{\text{mix}}^N(\epsilon) \ge \frac{C}{N^m}\Big(\frac{1}{2}-\epsilon\Big)\exp\Big(cN^{m-1} \Big),
    \]
    for some $C,c>0$ depending only on $\beta$ and $m$. \newline 
\end{Remark}

\begin{figure}
\centering
\begin{tikzpicture}[x=1cm,y=1cm, font=\small,
    >=Stealth, every node/.style={align=center}]
  \definecolor{fast}{RGB}{27,112,120}
  \definecolor{slow}{RGB}{173,77,41}

  \fill[fast!10] (0,0) rectangle (4,1.85);
  \fill[black!3] (4,0) rectangle (8,1.85);
  \fill[slow!10] (8,0) rectangle (12,1.85);

  \node[fast,font=\bfseries] at (2,0.98)
    {Logarithmic mixing};
  \node at (2,0.43) {$t_{\mathrm{mix}}\leq C'\log(N)$};

  \node[text=black!55] at (6,0.98) {Intermediate regime};

  \node[slow,font=\bfseries] at (10,0.98)
    {Exponential mixing};
  \node at (10,0.43) {$t_{\mathrm{mix}}\geq CN^{-m}e^{cN^{\ell-1}}$};

  \draw[->,semithick] (0,0) -- (12.65,0)
    node[right] {$\beta$};
  \draw (0,-0.08) -- (0,0.08);
  \node[below=5pt] at (0,0) {$0$};
  \foreach \x/\lab in {4/1,8/2} {
    \draw[densely dashed,black!45] (\x,0.08) -- (\x,1.95);
    \draw (\x,-0.08) -- (\x,0.08);
    \node[below=5pt] at (\x,0) {$\beta_{\lab}$};
  }
\end{tikzpicture}
\caption{Proven phase diagram for the mixing time of hyperlattice gauge theory on the torus $\Z_N^m$, when the model is defined on $k$-cells, with $k \in \{0,\dots,m-2\}$ and the dynamics are Glauber dynamics with heat-bath updates.}
\label{fig:phasetransition}
\end{figure}
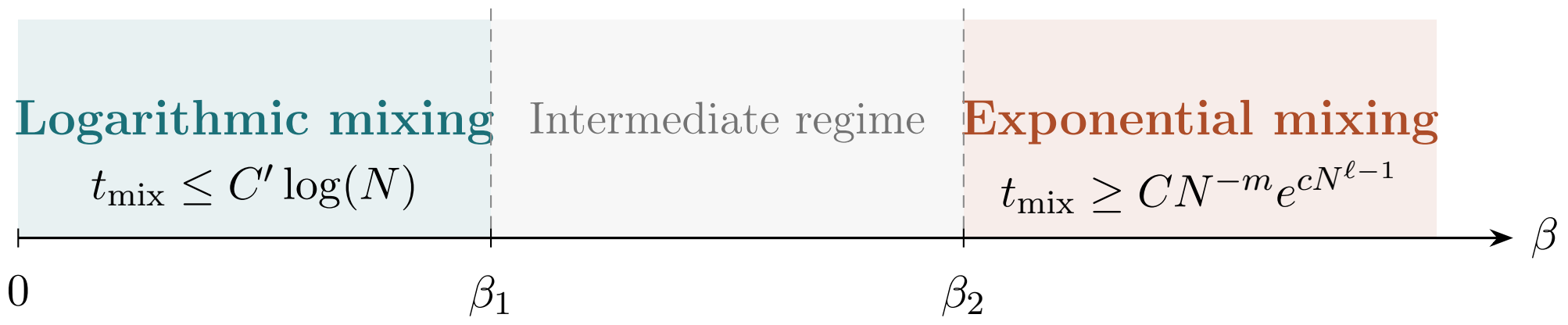

The result of Theorem \ref{thm:lowerboundlgt} does not hold true for small values of $\beta$, as the following Proposition shows.
\begin{Prop}
\label{prop:smallbeta}
    Consider the hyperlattice gauge theory probability measure $\mu_{\beta,N,m,k}$ on $\Z_N^m$ such that $k \in \{0 ,\dots, m-1\}$ and $m \ge 1$. Consider the Markov process $(\tau_t)_{t \ge 0}$ with generator $\mathcal{L}$ which corresponds to the transition kernel of Example \ref{ex:firstexample}(b), i.e., the heat-bath dynamics with Glauber updates. If 
    \begin{equation}
    \label{eq:smallbeta}
    \tanh(\beta) < \frac{1}{2(m-k)(2k+1)},
    \end{equation}
    then the mixing time is logarithmic in $N$, i.e., there exists $C$ that depends on $m,k,\beta,\epsilon$ only such that for any $\epsilon \in  (0,1)$
    \[
    t_{\text{mix}}^N(\epsilon) \le C\log(N).
    \]
\end{Prop}
Combining Theorem \ref{thm:lowerboundlgt} and Proposition \ref{prop:smallbeta} one deduces the existence of distinct rapid and slow mixing regimes of Glauber dynamics of hyperlattice gauge theory on the torus, as pictured in Figure~\ref{fig:phasetransition}. 
\subsection{Contribution and proof strategy}
The main contribution of this paper is to extend the bottleneck approach, outlined in \cite{Stahl2026SlowMixing}, to hyperlattice gauge theory in any dimension $k$ on the $m$-dimensional torus provided that the codimension satisfies $m-k \ge 2$. A first difficulty is that closed $(k+1)$-forms, which describe the defects, are dual to cycles of dimension $(m-k-1)$, and therefore don't need to be one-dimensional. This is the main reason for the work in Section~\ref{sec:dec} to obtain general isoperimetric inequalities to prove our result. The intuition for the proof of Theorem \ref{thm:lowerboundlgt} is the following. For each configuration $\tau \in C^k$ such that all connected components of $d\tau$ are small enough, we write $\tau=\tau_0+\tau_1$, where $d\tau_0=0$, in a way that $\tau_0$ is relatively close in Hamming distance to $\tau$. Their distance is quantified through the above-mentioned isoperimetric inequalities. On the other hand $\tau_1$ represents the \textit{perturbation} and induces the $(k+1)$-cells in the support of $d\tau$. We will see that one can categorize the defect-free representative $\tau_0$ into different classes (called \textit{cohomology classes}). In order to change the class to which the defect-free representative $\tau_0$ belongs, we prove that our dynamics need to pass through a configuration with high energy (meaning that  $d\tau$ contains connected components with large support), which is responsible for the exponentially slow mixing at low temperature. 

Moreover, we show that the slow mixing results are quite robust with respect to different dynamics, allowing the size of the blocks as defined in Section \ref{subsec:gd} to grow with the torus. A further difficulty is that a single block update may create several large vortex clusters simultaneously. This leads to a more general Peierls argument accounting for collections of large clusters.

At low values of $\beta$, we show that the Dobrushin conditions as presented in \cite{RapidMixingGlauber} are satisfied, when considering single-site Glauber updates using the heat-bath rule. This implies logarithmic mixing.

The rest of the paper is organized as follows. In Section~\ref{sec:dec}, we recall the basic definitions of discrete exterior calculus and topology and establish the isoperimetric results needed for the bottleneck argument. In Section~\ref{sec:rapidmixing}, we prove Proposition~\ref{prop:smallbeta} by verifying Dobrushin's condition for rapid mixing. In Section~\ref{sec:slowmixing}, we prove Theorem~\ref{thm:lowerboundlgt} using a bottleneck argument based on the tools developed in Section~\ref{sec:dec}. The Appendix based on an AI proof establishes a discrete isoperimetric inequality used in Section~\ref{sec:dec}.
\section{Discrete calculus and topology on the finite torus}
\label{sec:dec}
In this section we will introduce  some concepts of discrete exterior calculus and algebraic topology, following the notation of \cite{Hatcher2002}, among others. This section is crucial to define the model, and to properly write the bottleneck argument. We refer to \cite{eldridge2026} or \cite{LFVabelian} for an introduction to discrete exterior calculus in our $k$-cell models. 
\subsection{Discrete exterior calculus}
Consider the discrete torus $\Z_N^m$ in dimension $m \ge 2$ and $N \ge 3$. In this section, unless specified otherwise, we assume that $k \in \{0, \dots,m-1\}$.
Denote the standard coordinate vectors by $(e_1,\dots,e_m)$. For any $k \le m$, we call a $k$-cell an elementary $k$-dimensional cube that is defined by a base vertex $x \in \Z^m$ and a set $(e_{i_j})_{j=1}^k$ of coordinate directions spanning the cell. Hence, a $k$-cell will be written as
\[
[x;e_{i_1},\dots,e_{i_k}],
\]
and we denote the set of all $k$-cells by $C_k(\Z_N^m)$. Observe that $C_0(\Z_N^m)$ corresponds to the vertices in the torus, $C_1(\Z_N^m)$ to the edges and $C_2(\Z_N^m)$ to the plaquettes, recall Figure \ref{fig:firstpicture}. As we are working with the group $\Z_2$ throughout this paper, we do not need to consider orientated cells. 

A formal sum of $k$-cells with coefficients in $\Z_2$ is called a $k$-\textit{chain}. The set of $k$-chains over $\Z_N^m$ that have coefficients in $\Z_2$ is denoted by $C_k(\Z_N^m,\Z_2)$. A $k$-chain is written as
\begin{equation}
\label{eq:defkchain}
q=\sum_{c \in C_k(\Z_N^m)}a_cc, ~ c \in C_k(\Z_N^m),~ a_c \in \Z_2.
\end{equation}
\newline 
For $q \in C_k(\Z_N^m,\Z_2)$, we identify $q$ with its support for convenience. Namely, we write
\[
c\in q \iff c \in \supp q\iff q(c)=1.
\]
For a $k$-cell $c=[x;e_{i_1},...,e_{i_k}] \in C_k(\Z_N^m)$ we define the \textit{boundary} $\partial_k c$ as the following mapping 
\[
\partial_k c \coloneqq  \sum_{j=1}^k [x+e_{i_j};e_{i_1},...,e_{i_{j-1}},e_{i_{j+1}},...,e_{i_k}]+[x;e_{i_1},...,e_{i_{j-1}},e_{i_{j+1}},...,e_{i_k}].
\]
Then, $\partial_kc \in C_{k-1}(\Z_N^m,\Z_2)$. If $k=0$, we use the convention $\partial_0c=0$. We extend the definition of boundary to $k$-chains linearly. For $q=\sum_{c \in C_k(\Z_N^m)}a_cc \in C_k(\Z_N^m,\Z_2)$ one defines the \textit{boundary operator} as
\[
\partial_kq \coloneqq \sum_{c \in C_k(\Z_N^m)}a_c \partial_k c \in C_{k-1}(\Z_N^m,\Z_2).
\]
For $k$-cells $c \in C_k(\Z_N^m)$, we also define the \textit{upper boundary}  operator by
\[
\hat{\partial}_k c \coloneqq  \sum_{c' \colon c \in \partial_{k+1} c'}c' \in C_{k+1}(\Z_N^m,\Z_2).
\]
Similarly, one defines the \textit{upper boundary operator} for $k$-chains as follows
\begin{equation}
\label{eq:upperboundary}
\hat{\partial}_k \colon C_k(\Z_N^m,\Z_2)\to C_{k+1}(\Z_N^m,\Z_2), ~\hat{\partial}_k q\coloneqq  \sum_{c \in C_k(\Z_N^m)} a_c \hat{\partial}_kc.
\end{equation}
We define the kernel of $\partial_k$ as
\[
\ker \partial_k \coloneqq \{q \in C_k(\Z_N^m,\Z_2) \colon \partial_kq=0\},
\]
and the image $\text{Im}~\partial_{k+1} \subseteq C_k(\Z_N^m,\Z_2)$. We call the elements of $\ker \partial_k$ $k$-\textit{cycles} and the elements of $\text{Im}~\partial_{k+1}$ \textit{boundaries}.
With a similar notation as in \cite{Hatcher2002}, one defines the $k$-th \textit{homology group} as being the quotient group $H_k \coloneqq  \ker \partial_k / \text{Im}~\partial_{k+1}$. 
We say that two cycles $q,q' \in \ker \partial_k$ are \textit{homologous} if they are in the same equivalence class in $H_k$, i.e., if $q-q' \in \text{Im}~ \partial_{k+1}$.  
Finally a $k$-cycle $q$ is called \textit{null-homologous} if $[q]=0$, i.e., if there exists $Q \in C_{k+1}(\Z_N^m,\Z_2)$ such that $\partial_{k+1}Q =q$.  If a  $k$-cycle $q \in \ker \partial_k$ satisfies $[q] \neq 0$, then we call $q$ \textit{homologically non-trivial}.
\newline 
We denote the set of homomorphisms from $C_k(\Z_N^m,\Z_2)$ to $\Z_2$ by 
\[
C^k(\Z_N^m,\Z_2)\coloneqq \text{Hom}_{\Z_2}(C_k(\Z_N^m,\Z_2),\Z_2).
\] 
We call $\varphi$ a $k$-\textit{form}. For a $k$-form $\varphi$, one denotes the support of $\varphi$ as the following subset of $C_k(\Z_N^m)$,
\[
\supp \varphi \coloneqq  \{c \in C_k(\Z_N^m):\varphi(c)\neq 0\}.
\]
If $\varphi \in C^k(\Z_N^m,\Z_2)$ is a $k$-form and $q \in C_k(\Z_N^m,\Z_2)$ is a $k$-chain, we use the following abbreviation for the size of their supports
\[
|\varphi| \coloneqq |\supp \varphi|, ~ \text{and} ~ |q| \coloneqq |\supp q|.
\]
For a $k$-chain $q \in C_k(\Z_N^m,\Z_2)$ and $\varphi \in C^k(\Z_N^m,\Z_2)$ one defines
\begin{equation}
\label{eq:definitionphiq}
\varphi(q) \coloneqq  \sum_{c \in C_k(\Z_N^m)} a_c\varphi(c) \in \Z_2.
\end{equation}
We now define the \textit{discrete exterior operator} as the following map
\[
d_k:C^k(\Z_N^m,\Z_2) \to C^{k+1}(\Z_N^m,\Z_2), ~ d_k \varphi (\xi)\coloneqq  \sum_{c \in \partial_{k+1}\xi} \varphi(c), ~~ \xi \in C_{k+1}(\Z_N^m). 
\]
We call a $k$-form $\varphi$ \textit{closed} if $d_k\varphi=0$. \newline 
Observe that immediately by definition one has that $\partial_{k-1} \partial_k=0$ and $d_{k+1}d_k=0$. \newline
Note that for the rest of the paper, we adopt the following conventions, $C^{-1}=\{0\}, d_{-1} =0, d_m =0$. We might also use the abbreviation $C^k \coloneqq C^k(\Z_N^m,\Z_2)$.
\subsection{Results and duality}
\begin{Theorem}[Stokes' Theorem]
\label{thm:stokes}
    The operator $d_{k}$ is dual to the operator $\partial_{k+1}$, i.e., for any $\varphi \in C^k(\Z_N^m,\Z_2)$ and $q \in C_{k+1}(\Z_N^m,\Z_2)$, one has
    \[
    d_k\varphi(q)=\varphi(\partial_{k+1} q).
    \]
\end{Theorem}
One defines the \textit{cohomology group} $H^k(\Z_N^m,\Z_2)\coloneqq  \ker d_k / \text{Im}~d_{k-1}$. For the group $H^k(\Z_N^m,\Z_2)$, one has the following result.
\begin{Theorem}[{\cite[p.198]{Hatcher2002}}]
\label{thm:universalthmcohomology}
\[
H^k(\mathbb Z_N^m,\mathbb Z_2)
\cong
\operatorname{Hom}_{\mathbb Z_2}
\bigl(H_k(\mathbb Z_N^m,\mathbb Z_2),\mathbb Z_2\bigr).
\]
\end{Theorem}
The isomorphism in Theorem \ref{thm:universalthmcohomology} can be seen as follows. Let $[\varphi] \in H^k(\Z_N^m,\Z_2)$. Then $[\varphi]$ can be seen as a mapping from $[q] \in H_k(\Z_N^m,\Z_2)$ to $\Z_2$ through the following map 
\begin{equation}
\label{eq:homomorphicphiq}
([\varphi],[q])=\varphi(q) \in \Z_2,
\end{equation}
which by Theorem \ref{thm:universalthmcohomology} is independent of the choice of representative in the equivalence class.

Throughout the next results, an important tool will be duality. Hence we let $*$ denote a dual map, which maps a $k$-cell $c$ to a dual $(m-k)$-cell denoted $*c$. We refer to \cite[Section 2]{Forsstr_m_2023} for a detailed introduction to the construction of the dual lattice.  
\begin{Remark}
    In \cite[Section 2]{LFVabelian} by the construction of the dual cell, it holds that $*(*c)$ and $c$ are equal up to orientation. In our non-oriented setting as the gauge group is $\Z_2$, we have that $*(*c)=c$. It also means that we use the same notation for the primal and dual lattices, with each star map sending objects of one lattice to objects on the other.
\end{Remark}
For a set of $k$-cells $C \subseteq C_k(\Z_N^m)$ one defines
\[
*C\coloneqq \{*c : c \in C \}\subseteq C_{m-k}(\Z_N^m).
\]
For a $k$-chain $q \in C_k(\Z_N^m,\Z_2)$, we let  $*^cq \in C^{m-k}(\Z_N^m,\Z_2)$ be the unique $(m-k)$-form such that
\[
\supp *^cq = *\supp q.
\]
Furthermore, if $\varphi \in C^k(\Z_N^m,\Z_2)$ is a $k$-form, then one defines its dual as the unique $(m-k)$-chain $*^f\varphi \in C_{m-k}(\Z_N^m,\Z_2)$,
such that
\[
\supp *^f\varphi = * \supp \varphi.
\]
\begin{Lemma}
\label{lem:duality} 
    Let $k \in \{0,\dots,m\}$. For any $q \in C_k(\Z_N^m,\Z_2)$ and any $\varphi \in C^{k}(\Z_N^m,\Z_2)$, the following holds.
    \begin{enumerate}[label=(\alph*)]
\item $\displaystyle \qquad |\supp *^f\varphi|=|\supp \varphi|,\quad \text{and}\quad |\supp *^cq|=|\supp q|.$   
\item $\begin{aligned}[t]
\qquad \qquad *^c(*^f\varphi) = \varphi \in C^k(\Z_N^m,\Z_2) ~
\text{and}~ *^f(*^cq) = q \in C_k(\Z_N^m,\Z_2).
\end{aligned}$
\item For any $c \in C_k(\Z_N^m)$ one has $\supp *^c(\hat{\partial}_{k}c)=*\supp\hat{\partial}_{k}c=\supp \partial_{m-k}(*c)$.
\item \refstepcounter{equation}\label{eq:dualityrelation}%
\makebox[\linewidth]{%
    \hfil
    $\displaystyle *^c(\partial_kq)=d_{m-k}(*^cq) \quad \text{and} ~~\displaystyle *^f(d_k \varphi)=\partial_{m-k}(*^f\varphi).$
    \hfil
    \llap{\normalfont(\theequation)}%
}
    \end{enumerate}
\end{Lemma}
\begin{proof}
\begin{enumerate}[label=(\alph*)]
    \item  The result follows immediately from the definition of $*^c$ and $*^f$. 
    \item Since the forms and chains are $\Z_2$-valued, if we prove that the left- and right-hand sides have the same support, we prove the claim. To this end note that
    \[
    \supp *^c(*^f\varphi) = * \supp *^f\varphi = ** \supp \varphi= \supp \varphi.
    \]
    The second claim is proven in exactly the same way.
\item This follows immediately from the construction of the dual cells. Namely, if $c \in C_{k}(\Z_N^m), c' \in C_{k+1}(\Z_N^m)$, then  $c' \in \supp \hat{\partial}_kc$ iff $*c' \in \supp \partial_{m-k}*c$. This concludes the proof.
    \item We start proving the first equality. If we prove the claim for an arbitrary $k$-chain $q \in C_k(\Z_N^m,\Z_2)$ that has exactly one $k$-cell $c \in C_k(\Z_N^m)$ as support, then the result extends to any $k$-chain $q \in C_k(\Z_N^m,\Z_2)$ by linearity of $*^c,\partial_k, d_{m-k}$ and this finishes the proof. \newline 
    To this end, let $q \in C_k(\Z_N^m,\Z_2)$ be such that $\supp q=\{c\}$. Then, observe that the left-hand side and right-hand side of \eqref{eq:dualityrelation} are both elements of $C^{m-k+1}(\Z_N^m,\Z_2)$. Let $\xi \in C_{m-k+1}(\Z_N^m)$. Then,
\begin{align*}
&\xi \in \supp *^c(\partial_kq) \iff a \coloneqq *\xi \in \supp \partial_kq \iff c \in \supp \hat{\partial}_{k-1} a \\\iff &c \in *\supp \partial_{m-k+1}*a \iff *c \in \supp (\partial_{m-k+1} \xi).
\end{align*}
Here the first and last equivalences are immediate from the dual definition. For the second one, note that it follows directly from the upper boundary definition. The third equivalence is $(c)$. 
Now, observe that $\supp*^cq=\{*c\}$, hence the last equation is in turn equivalent to
\begin{align*}
&*c \in \supp (\partial_{m-k+1} \xi) \iff *^cq(\partial_{m-k+1} \xi)=1 \\ \iff &d_{m-k}(*^cq)(\xi)=1 \iff \xi \in \supp d_{m-k}(*^cq).
\end{align*}
Where the second equivalence follows from Stokes' Theorem. This concludes the proof of the first equality.

For the second equality, apply the first equality to the $(m-k)$-chain $*^f\varphi$. Together with (b) we deduce that
\[
*^c(\partial_{m-k}*^f\varphi)=d_k(\varphi).
\]
Now, one applies $*^f$ on both sides which gives the claim by (b). This concludes the proof.
\end{enumerate}
\end{proof}
Observe that the main result from the above Lemma is (d), meaning that the dual of a cycle is a closed form, and conversely.

Poincar\'e's duality theorem is one of the main tools in discrete exterior calculus, and takes the following form in our setting.
\begin{Theorem}{\cite[Theorem 3.30]{Hatcher2002}}
    \label{thm:poincareduality}
    The groups $H^k(\Z_N^m,\Z_2)$ and $H_{m-k}(\Z_N^m,\Z_2)$ are isomorphic
    \[
    H^k(\Z_N^m,\Z_2) \cong H_{m-k}(\Z_N^m,\Z_2).
    \]
\end{Theorem}
\begin{Remark}
    Observe that Theorem \ref{thm:universalthmcohomology} and Theorem \ref{thm:poincareduality} give two different perspectives on cohomology. The first Theorem allows us to view a cohomology class $H^k$ as a linear functional from $H_k$ to $\Z_2$, while the second identifies it with a homology class in $H_{m-k}$. That is consistent with our construction $*^f$ that maps closed $k$-forms to $(m-k)$-cycles on the dual lattice and acts on equivalence classes by $[\varphi] \mapsto [*^f\varphi]$.
\end{Remark} 
\subsection{(Co)homology on the finite torus and support condition}
In this section we establish some properties and results of homology and cohomology on the finite torus that will be used in the bottleneck argument for our lower bound on the mixing time. Although some of the results might be standard from the viewpoint of algebraic topology, we include proofs in our discrete setting to make it more accessible for an audience from probability or statistical mechanics. \newline 
The next results will help to connect null-homology (or null-cohomology) with a support condition. In order to prove the first support condition in Lemma \ref{lem:homologygroup}, we need to prove the existence of a $(k+1)$-chain having the sum of two specific $k$-chains as boundary, which is the next Lemma.
\begin{Lemma}
\label{lem:existenceofk+1chain}
Let $k \in \{0,\dots,m-1\}$, let $\mathcal{I} \coloneqq \{i_1,\dots,i_k\} \subsetneq [m]$ and let $\Gamma \in C_k(\Z_N^m,\Z_2)$ be a $k$-cycle such that
\begin{equation}
\label{eq:supportconditiongamma}
\supp \Gamma \subset \{ [x;e_{i_1},\dots,e_{i_k}]:x \in \Z_N^m \}. 
\end{equation}
For $l \in [m]\setminus \mathcal{I}$, let  $T_l:C_k(\Z_N^m) \to C_k(\Z_N^m)$ be the operator that shifts the cell by one in the direction $l$, i.e., $T_l[x;e_{i_1},...,e_{i_k}]=[x+e_l;e_{i_1},...,e_{i_k}]$ and let 
    \begin{equation}
    \label{eq:gammaandgamma'}
    \Gamma' \coloneqq  T_l\Gamma,
    \end{equation}
    where $T_l$ acts cell-by-cell in the support of $\Gamma$. \newline
    Then, there exists a $(k+1)$-chain $\Sigma \in C_{k+1}(\Z_N^m,\Z_2)$ such that $\partial_{k+1} \Sigma = \Gamma +\Gamma'$.
\end{Lemma}
\begin{proof}
Let $k \in \{0, \dots,m-1\}$, then for $l \in [m] \setminus \mathcal{I}$ and $c=[x;e_{i_1},...,e_{i_k}] \in C_k(\Z_N^m)$ such that $\{i_1, \dots,i_k\}=\mathcal{I}$ one defines the operator $P_l$ as follows
\[
     P_l(c) \coloneqq [x;e_{i_1},...,e_{i_k},e_l].
    \]
    We now claim that for $c \in \supp \Gamma$, it holds that $\partial_{k+1}P_l(c)+P_l(\partial_kc)=c+T_lc$. \newline 
    To see this note that
    \begin{align*}
    &\partial_{k+1}P_l(c)+P_l(\partial_kc)= \Big([x+e_l;e_{i_1},...,e_{i_k}]+[x;e_{i_1},...,e_{i_k}]\Big) +\\
        \sum_{r \in [k]} &\Big([x+e_{i_r};e_{i_1},\dots,e_{i_{r-1}},e_{i_{r+1}},\dots,e_{i_k},e_l]+[x;e_{i_1},\dots,e_{i_{r-1}},e_{i_{r+1}},\dots e_{i_k},e_l]\Big) \\
        +\sum_{r \in [k]} &\Big([x+e_{i_r};e_{i_1},\dots,e_{i_{r-1}},e_{i_{r+1}},\dots,e_{i_k},e_l]+[x;e_{i_1},\dots,e_{i_{r-1}},e_{i_{r+1}},\dots,e_{i_k},e_l]\Big) \\
        =&c+T_lc.
    \end{align*}
    As the second and third line cancel each other. This proves the claim. 

    In particular, for any $k$-chain $\Gamma \in C_k(\Z_N^m,\Z_2)$ described as in the Lemma, by linearity of both $T_l$ and $P_l$ one deduces that
    \[
    \partial_{k+1}P_l(\Gamma)+P_l(\partial_k \Gamma)=\Gamma+T_l\Gamma.
    \]
    If one defines 
    \[
    \Sigma \coloneqq P_l\Gamma=\sum_{c \in \supp \Gamma} P_lc \in C_{k+1}(\Z_N^m,\Z_2),
    \]
    then by linearity of $P_l$ and using the fact that $\Gamma$ is a $k$-cycle, i.e., $\partial_k \Gamma =0$, one obtains 
    \[
\partial_{k+1}\Sigma=\partial_{k+1}P_l\Gamma=\partial_{k+1}P_l\Gamma
+P_l(\partial_k \Gamma)=\Gamma+T_l\Gamma=\Gamma+\Gamma'. 
    \]
    This proves the claim.  
\end{proof}
In the next Lemma, we view cohomology classes as mappings as suggested by Theorem \ref{thm:universalthmcohomology} and Theorem \ref{thm:poincareduality}.
\begin{Lemma}
\label{lem:equivhomologyclassloops}
    Let $q_1,q_2 \in \ker \partial_k$ be two $k$-cycles. Then, $q_1$ and $q_2$ are in different homology classes $[q_1]\neq [q_2] \in H_k(\Z_N^m,\Z_2)$ if and only if there exists a homologically non-trivial $(m-k)$-cycle $[\Gamma] \in H_{m-k}(\Z_N^m,\Z_2)$ such that,
    \[
    *^cq_1(\Gamma) \neq *^cq_2(\Gamma).
    \]
\end{Lemma}
\begin{proof}
By Theorem \ref{thm:universalthmcohomology} and Theorem \ref{thm:poincareduality}, for $[q] \in H_k$ and $[\Gamma] \in H_{m-k}$, the mapping
\[
[q] \longmapsto \bigl([\Gamma]\longmapsto (*^c q)(\Gamma)\bigr)
\]
is a homomorphism from $H_k(\Z_N^m,\Z_2)$ to $\text{Hom}_{\Z_2}(H_{m-k}(\Z_N^m,\Z_2),\Z_2)$. In particular, $[q_1]\neq [q_2]$ if and only if $[q_1]$ and $[q_2]$ do not represent the same homomorphism in $\text{Hom}_{\Z_2}(H_{m-k}(\Z_N^m,\Z_2),\Z_2)$. Hence, it means that there is $[\Gamma] \in H_{m-k}(\Z_N^m,\Z_2)$ such that
\[
([*^cq_1])([\Gamma]) \neq ([*^cq_2])([\Gamma]).
\]
Moreover, $[\Gamma]$ is non-trivial, as if $[\Gamma]=0$, then one automatically has $([q])(\Gamma)=0$ for each $[q] \in H_k(\Z_N^m,\Z_2)$. This concludes the proof.
\end{proof}
Using the same proof as above but only considering $[\varphi_1],[\varphi_2] \in H^k(\Z_N^m,\Z_2)$, then by Theorem \ref{thm:universalthmcohomology} one has that
\[
[\varphi] \longmapsto \bigl([\Gamma]\longmapsto (\varphi)(\Gamma)\bigr)
\]
is a homomorphism from $H^k(\Z_N^m,\Z_2)$ to $\text{Hom}_{\Z_2}(H_{k}(\Z_N^m,\Z_2),\Z_2)$ from which we deduce the following result.
\begin{Corollary}
\label{cor:cohomologyclasses}
The classes $[\varphi_1],[\varphi_2] \in H^k(\Z_N^m,\Z_2)$ are distinct if and only if there exists a homologically non-trivial $[\Gamma] \in H_k(\Z_N^m,\Z_2)$ such that $\varphi_1(\Gamma) \neq \varphi_2(\Gamma)$.
\end{Corollary}
Picking $[q]=[q_1] \neq 0,~ \text{and}~[q_2]=0$ in Lemma \ref{lem:equivhomologyclassloops}, gives the following Corollary.
\begin{Corollary}
    \label{cor:cohomology111}
    Consider any non-zero class $[q] \in H_k(\Z_N^m,\Z_2)$. Then there exists a non-zero class $ [\Gamma] \in H_{m-k}$ such that 
    \[
    *^cq(\Gamma)=1.
    \]
    The converse holds true.
\end{Corollary}
Recall that for two $k$-forms $\varphi,\varphi' \in C^k(\Z_N^m,\Z_2)$, we denote by $\varphi+\varphi' \in C^k(\Z_N^m,\Z_2)$ their cell-by-cell sum
\begin{equation}
\label{eq:notaddition}
(\varphi + \varphi')(c)= \begin{cases}
    0, ~~ \varphi(c)=\varphi'(c), \\
    1, ~ \text{else}.
\end{cases}
\end{equation}
In the next result, we combine previous results to prove the following support bound, that separates homology classes by a support condition on the torus.
\begin{Lemma}
    \label{lem:homologygroup}
    Let $k \in \{0,\dots,m\}$ and let $q,q' \in C_k(\Z_N^m,\Z_2)$ be cycles such that $[q] \neq [q'] \in H_k(\Z_N^m,\Z_2)$. Then, 
    \[
    |\{c \in C_{m-k}(\Z_N^m): *^cq(c) \neq *^cq'(c)\}| \ge N^k. 
    \]
\end{Lemma}
\begin{proof}
    Let $\tau\coloneqq *^c(q +q')\in C^{m-k}(\Z_N^m,\Z_2)$. Observe that $\tau$ is closed by Lemma~\ref{lem:duality}(d): 
    \[
    d_{m-k}\tau=d_{m-k}(*^c(q+q'))=*^c(\partial_k(q+q'))=0.
    \]
    Moreover, note that $[q+q'] \neq 0 \in H_k(\Z_N^m,\Z_2)$. To see it, take $\tilde{\Gamma} \in C_{m-k}(\Z_N^m,\Z_2)$ such that $*^cq(\tilde{\Gamma}) \neq *^cq'(\tilde{\Gamma})$, which exists by Lemma~\ref{lem:equivhomologyclassloops}. Then 
    \[
    \tau(\tilde{\Gamma})=*^cq(\tilde{\Gamma})+*^cq'(\tilde{\Gamma})=1.
    \]
    This means that $[q+q'] \neq 0$ by Corollary \ref{cor:cohomology111}.
    
    Now, we claim that there exists $\Gamma \in C_{m-k}(\Z_N^m,\Z_2)$ such that $[\Gamma] \neq 0$, with $\tau(\Gamma)=1$ that has the form described in Lemma~\ref{lem:existenceofk+1chain}, namely that its support satisfies \eqref{eq:supportconditiongamma}. In words one searches for a coordinate torus $\Gamma$ and not only a generic $(m-k)$-cycle $\tilde{\Gamma}$.

    Let $I\subset[m]$ be any subset of size $m-k$. Then, define 
    \begin{equation}
    \label{eq:defGamma_IwithC_I}
    C_I \coloneqq \bigcup_{\substack{x \in \Z_N^m\\ x_j=0, j \notin I}}\{ [x;e_{i_1},\dots,e_{i_{m-k}}]:I=\{i_1,\dots , i_{m-k} \}\} \subset C_{m-k}(\Z_N^m).
    \end{equation}
    Let $\Gamma_I$ be the $(m-k)$-chain such that $\supp \Gamma_I=C_I$.
    
In particular, observe that $\partial_{m-k}\Gamma_I=0$ (as each $(m-k-1)$-cell in the boundary is counted twice and hence vanishes). Now, by \cite[Theorem 2.35]{Hatcher2002}, the homology on the torus is independent of the cell-complex studied. Hence one can consider the case where there is one vertex on the torus, and one edge in each direction connecting the vertex to itself. In this simplified cell-complex, there is for each $I$ as described above, a unique $(m-k)$-cell spanned by the directions in $I$. We call this cell $c_I$. In particular one has $\partial_{m-k}c_I=0$ as the $(m-k-1)$-faces of the cell cancel. One has a one-to-one correspondence between $\Gamma_I$ from above and $\gamma_I$ the $(m-k)$-cycle on this simplified cell-complex such that $\supp \gamma_I=\{c_I\}$. In particular, one deduces the following for the homology group on this simplified cell-complex:
\begin{equation}
\label{eq:basisHm-k}
H_{m-k} = \{\sum_{\substack{I \subset [m]\\ |I|=m-k }} a_I\gamma_I,~a_I \in \Z_2 \}.
\end{equation}
Subdividing each circle into $N$ edges gives the original lattice $\Z_N^m$. Under this subdivision, the cycle represented by $\gamma_I$ becomes $\Gamma_I$, the sum of the lattice cells in the corresponding coordinate torus. Since subdivision preserves homology by \cite[Theorem 2.35]{Hatcher2002}, the classes $([\Gamma_I])_I$ form a basis of $H_{m-k}(\Z_N^m,\Z_2)$.
\[
H_{m-k}(\Z_N^m,\Z_2) = \{\sum_{\substack{I \subset [m]\\ |I|=m-k }} a_I[\Gamma_I],~a_I \in \Z_2 \},
\]

Recall the $k$-cycle $\tilde{\Gamma}$ from above. By the basis decomposition, there exists coefficients $\tilde{a}_{I} \in \Z_2$ such that $[\tilde{\Gamma}] = \sum_{\substack{I \subset [m]\\ |I|=m-k }} \tilde{a}_I[\Gamma_I]$. In particular
\[
1=\tau(\tilde{\Gamma})=\sum_{\substack{I \subset [m]\\ |I|=m-k }} \tilde{a}_I\tau(\Gamma_I).
\]
meaning that there exists at least one $I$ such that $\tilde{a}_I=1$ and $\tau(\Gamma_I)=1$. Taking $\Gamma=\Gamma_I$ proves the claim: there is a coordinate $(m-k)$-cycle $\Gamma$ such that $\tau(\Gamma)=1$.

 Then, as $\tau(\Gamma)=1$, there exists at least one $c \in \supp \Gamma$ such that $c \in \supp \tau$. Now, define 
\[
X \coloneqq \{x \in \Z_N^m: x_i =0 ~\forall i \in I\},
\]
where $I$ is the set of coordinates spanning the cells in the support of $\Gamma$.
For any $x \in X$, we define
\[
    C^x \coloneqq \bigcup_{\substack{x' \in \Z_N^m\\ x'_j=x_j, j \notin I}}\{ [x';e_{i_1},\dots,e_{i_{m-k}}]:I=\{i_1,\dots , i_{m-k} \}\} \subset C_{m-k}(\Z_N^m).
\]
Accordingly, define $\Gamma^x$ as the $(m-k)$-chain such that $\supp \Gamma^x = C^x$. It holds that $(\Gamma^x)_{x \in X}$ are the parallel translates of $\Gamma$ in directions that do not span the cells.
Let $x \in X$ such that $\|x\|_1=1$. Then $\Gamma,\Gamma^x$ differ by a unit translation in a coordinate direction that does not span their cells and thus satisfy the condition of Lemma \ref{lem:existenceofk+1chain}. By that Lemma, there exists $\Sigma \in C_{m-k+1}(\Z_N^m,\Z_2)$ such that $\partial\Sigma=\Gamma+\Gamma^x$. As $d_{m-k}\tau=0$ one deduces that 
    \[
    0=d_{m-k}\tau(\Sigma)=\tau(\partial_{m-k+1} \Sigma)=\tau(\Gamma)+\tau(\Gamma^x) \implies \tau(\Gamma^x)=\tau(\Gamma)=1.
    \]
Using the same argument recursively one deduces that for each $x \in X$ : $\tau(\Gamma^x)=1$. Because of this, for each of the $(m-k)$-chains $\Gamma^x, ~{x \in X}$ there is at least one $(m-k)$-cell $c^x \in \supp \Gamma^x$ such that $c^x \in \supp \tau$. As for $x \neq x'$ one has $\supp \Gamma^x \cap \supp \Gamma^{x'} = \emptyset$, it follows that the cells $(c^x)_{x \in X}$ are pairwise distinct and hence
\[
|\{c \in C_{m-k}(\Z_N^m):*^cq(c) \neq *^cq'(c)\}|=|\supp \tau| \ge |\{c^x:x \in X\}|=|X|=N^k
\]
which concludes the proof of the Lemma.
\end{proof}
Following from the proof, one deduces the following corollary. 
\begin{Corollary}
    \label{cor:cohomology11}
    Consider any $k$-cycle $q \in \ker \partial_k$  on the torus $\Z_N^m$ such that $[q] \neq 0$. Then,
    \[
    |\supp q| \ge N^k.
    \]
\end{Corollary}
\begin{proof}
    Pick $[q]\neq 0$ and pick the chain with no support denoted $q'=0 \in C_{k}(\Z_N^m,\Z_2)$ as a representative for the null-homology class $[q']$ in Lemma \ref{lem:homologygroup}. This lemma, together with Lemma \ref{lem:duality}(a), proves the claim
    \[
    |\supp q|=|\supp *^cq|=|\supp (*^cq+*^c0)| \ge N^k.
    \]
\end{proof}
Lemma \ref{lem:homologygroup} has furthermore the following corollary on $k$-forms.
\begin{Corollary}
\label{cor:lowerbooundhomologyclasses}
    Let $\varphi_1,\varphi_2 \in C^{k}(\Z_N^m,\Z_2)$ be two closed $k$-forms such that $[\varphi_1] \neq [\varphi_2]$. Then
    \[
    |\{c \in C_k(\Z_N^m): \varphi_1(c) \neq \varphi_2(c) \}| \ge N^{m-k}.
    \]
\end{Corollary}
\begin{proof}
    By Theorem \ref{thm:poincareduality}, the respective duals of $\varphi_1$ and $\varphi_2$, the $(m-k)$-cycles $*^f\varphi_1$ and $*^f\varphi_2$, are in different homology classes in $H_{m-k}(\Z_N^m,\Z_2)$, i.e., $[*^f\varphi_1] \neq [*^f\varphi_2]$. 
    Hence, by Lemma \ref{lem:homologygroup}, we have 
    \[
    |\{c \in C_k(\Z_N^m): *^c*^f\varphi_1(c) \neq *^c*^f\varphi_2(c) \}| \ge N^{m-k}.
    \]
    By Lemma \ref{lem:duality}(b), it holds that $*^c*^f$ is the identity. This concludes the proof.
\end{proof}
\begin{Lemma}
    \label{lem:cohomology1}
    Let $\omega \in C^{k+1}(\Z^m_N,\Z_2)$ be a closed $(k+1)$-form. If
    \begin{equation}
    \label{eq:condomega}
    |\supp \omega| < N^{m-(k+1)},
    \end{equation}
     then there exists a $k$-form $\tau \in C^k(\Z^m_N,\Z_2)$ such that $d\tau = \omega$.
\end{Lemma}
\begin{proof}
    Consider the dual of $\omega$, denoted $q \coloneqq *^f\omega \in C_{m-k-1}(\Z_N^m,\Z_2)$. As $d_{k+1}\omega=0$, it follows by Lemma \ref{lem:duality}(d) that $\partial_{m-k-1} q=0$. Now by Lemma \ref{lem:duality}(a) one has
    \[
    |\supp q|=|\supp \omega| < N^{m-(k+1)},
    \]
    which implies by Corollary \ref{cor:cohomology11} that $q$ is null-homologous. Therefore there exists some $(m-k)$-chain $Q$ such that $\partial_{m-k}Q=q$. Then, for the dual of $Q$, $*^cQ \in C^k(\Z_N^m,\Z_2)$, it holds by Lemma~\ref{lem:duality}(d) for the first equality and (b) for the last equality that
    \[
    d_{k}(*^cQ)=*^c(\partial_{m-k}Q)=*^cq=*^c*^f\omega=\omega.
    \]
    Taking $\tau \coloneqq *^cQ$ completes the proof.
\end{proof}
In Lemma \ref{lem:cohomology1} we proved that if a closed $(k+1)$-form $\omega$ satisfies a support condition, then there exists a $k$-form $\tau$ such that $d\tau=\omega$. In the next lemma, we answer the question of whether the support of such a $\tau$ can be bounded in terms of the support of $\omega$. More generally, inequalities of this type are studied for $k$-cycles $q$, and the question is whether there exists a filling  $Q$ (meaning $\partial_{k+1} Q=q$) with controlled support. The dual translation of this is exactly our question on $k$-forms. We note that much work has been done on these types of questions, see for example \cite{Dotterrer2016FillingCube,Federer}. Such inequalities in suitable settings give that $|\supp Q|^{k}$ is upper bounded by a constant times $|\supp q|^{k+1}$.
We could not, however, find a directly applicable statement for the discrete torus $\Z_N^m$. With the help of AI, we establish in Appendix~\ref{sec:app} a 
proof in the setting of the discrete torus, showing bounds similar to those in other settings. This result leads to the following dual formulation.
\begin{Lemma}
\label{lem:dec2}
    Let $k \in \{0,\dots,m-2\}$ and let $\omega \in C^{k+1}(\Z^m_N,\Z_2)$ be a closed $(k+1)$-form that satisfies the support condition \eqref{eq:condomega}. Then, there exists a $k$-form $\tau \in C^k(\Z^m_N,\Z_2)$ such that $d\tau=\omega$ and 
    \begin{equation}
    \label{eq:boundonsupport}
|\supp \tau| \le \frac{3 \cdot 2^{m-k}}{2^{m-(k+1)}-1}m|\supp \omega|^{\frac{m-k}{m-k-1}}. 
    \end{equation}
\end{Lemma}
For the rest of this paper, we denote the constant in the isoperimetric inequality~\eqref{eq:boundonsupport} by
\[
\hat{C}_{m,k} \coloneqq \frac{3 \cdot 2^{m-k}}{2^{m-(k+1)}-1}m.
\]
\begin{Remark}
We note the following concerning the isoperimetric inequality. \begin{enumerate}
    \item Observe that for any $m \in \N_{\ge 2}$ and $k \in \{0,\dots,m-2\}$ one has the following bound, independent of $k$
\[
6m < \hat{C}_{m,k} \le 12m.
\]
\item The bound is not sharp. Take for example $k=1,m=3,N \ge 5$ and a closed $\omega \in C^{k+1}(\Z_N^m,\Z_2)$ with $|\omega|=4$. Then, one can prove that there exists a $1$-form that satisfies $d\tau=\omega$ having the property that $|\tau|=1$, see Figure \ref{fig:degofgraph}(a) for a graphical representation. Applying the above lemma only gives $|\tau| \le 576$.
\end{enumerate}
\end{Remark}
We now prove the Lemma using the result in Appendix \ref{sec:app}.
\begin{proof}
Using the same arguments as in the proof of Lemma \ref{lem:cohomology1} for $q \coloneqq *^f \omega$, one has that $[q]$ is null-homologous.
We apply Proposition \ref{prop:appB} which gives the existence of a filling $Q \in C_{m-k}(\Z_N^m,\Z_2)$ such that $\partial_{m-k}Q=q$ and 
\[ 
|\supp Q| \le \frac{3 \cdot 2^{m-k}}{2^{m-(k+1)}-1}m|\supp q|^{\frac{m-k}{m-(k+1)}}.
\]
Let $\tau\coloneqq *^cQ \in C^k(\Z_N^m,\Z_2)$. Then we have
\[
d_{k}\tau=d_{k}*^cQ=*^c\partial_{m-k}Q=*^c*^f\omega=\omega,
\]
by Lemma \ref{lem:duality} (b) and (d) and $|\supp Q|=|\supp \tau|$ by Lemma \ref{lem:duality}(a). This concludes the proof.
\end{proof}
We now present two results which will be relevant in later sections. The next lemma counts the number of cohomology classes on the discrete torus. 
\begin{Lemma}
    \label{lem:numbercohomologyclasses}
    There are
    \[
    |H^k(\Z_N^m,\Z_2)|=2^{\binom{m}{k}}
    \]
    cohomology classes.
\end{Lemma}
\begin{proof}
Recall that in the proof of Lemma \ref{lem:homologygroup} one has computed a basis for the homology group $H_{m-k}$ in \eqref{eq:basisHm-k}. The basis has $\binom{m}{m-k}=\binom{m}{k}$ elements, one for each $I \subset [m]$ such that $|I|=m-k$. 
From this we immediately deduce that $|H_{m-k}|=2^{\binom{m}{k}}$. The conclusion follows from Theorem \ref{thm:poincareduality}.
\end{proof}
We define the \textit{adjacency graph} $\mathcal{G}_{m}^k(N)$ as the graph having the set of $k$-cells in $\Z_N^m$ as vertices, i.e., $V(\mathcal{G}_{m}^k(N)) = C_k(\Z_N^m)$. Furthermore, for two distinct vertices $c_1,c_2 \in V(\mathcal{G}_{m}^k(N))$, we connect them with an edge if there exists a $(k+1)$-cell denoted $C$ such that $c_1,c_2 \in \partial C$. 
\begin{Lemma}
\label{lem:sizegraph}
    It holds that $|V(\mathcal{G}_m^k(N))|= \binom{m}{k}N^m$ and each vertex $v \in V(\mathcal{G}_m^k(N))$ has degree $2(2k+1)(m-k)$.
\end{Lemma}
\begin{proof}
    A $k$-cell $c=[x;e_{i_1},...,e_{i_k}]$ is uniquely defined by a vertex $x \in \Z_N^m$ and an unordered choice of $k$ directions among $m$. There are $N^m$ points in the torus $\Z_N^m$. Furthermore, one has $\binom{m}{k}$ possibilities of choosing an unordered subset of size $k$ from a set of size $m$. This establishes the first claim. 
    
    Now we prove the second claim. Pick any vertex $c \in V(\mathcal{G}_m^k(N))$. It is connected by an edge to all $k$-cells $c'$ such that $c' \in \partial \hat{\partial}c \setminus \{c\}$. Since $N \ge 3$, two distinct $k$-cells are faces of at most one common $(k+1)$-cell. Then, note that in $m$ dimensions, a $k$-cell is in the boundary of $2(m-k)$ distinct $(k+1)$-cells. Namely, a $k$-cell $c$ is defined by $k$ fixed directions spanning the cell. Hence, to get $C \in \hat{\partial}c$, one has to choose among $(m-k)$ coordinates to define a new $(k+1)$-cell $C$ such that $c \in C$. Then, once a direction is fixed there are two $(k+1)$-cells which have $c$ in their boundary. We deduce that $|\hat{\partial}c|=2(m-k)$.
    
    Moreover, for a $(k+1)$-cell $C$, one has $|\partial C|=2(k+1)$. This holds because when taking the boundary, one has $(k+1)$ directions to choose from, and for each fixed direction there are $2$ possible $k$-cells. Lastly, note that $c$ is already included in the boundary of each $(k+1)$-cell once, so that there are $2(k+1)-1=2k+1$ $k$-faces left.  
    One deduces that $|\partial \hat{\partial}c \setminus\{c\}| \le 2(2k+1)(m-k)$. Lastly, note that $c_1,c_2 \in \hat{\partial}c$ have only $c$ as a common face. In particular one has equality in the above equation
    \[
    |\partial \hat{\partial}c \setminus\{c\}| = 2(2k+1)(m-k),
    \]
    which concludes the proof.
\end{proof}

\section{Rapid mixing at high temperature}
\label{sec:rapidmixing}
This section is devoted to proving Proposition \ref{prop:smallbeta}.
This result will follow from the work on rapid mixing under Glauber updates in \cite{RapidMixingGlauber}. In particular, we will prove that the setup of the heat-bath dynamics as described in Example \ref{ex:firstexample}(b) satisfies the so-called Dobrushin condition, which implies rapid mixing. The proof will turn out to be similar to that for the Ising model, see \cite{Hayes2006} for example.
\newline 
Let us first introduce some notation. We denote the blocks by the cell they represent, namely $B_c=\{c\}$ for $c \in C_k(\Z_N^m)$, which is the cell selected by the dynamics. For $\tau \in C^k$, let $P^c(\tau,\cdot)$ denote the distribution of the value at $c$ of the new state $\tau'$. We furthermore denote by $\tau_c^{0}$ and $ \tau_c^1$ the configurations that equal $\tau$ except at $c \in C_{k}(\Z_N^m)$ and takes the values $0$ and $1$ at $c$, respectively. Using this notation, the heat-bath rule resamples the value of $\tau$ at $c$ according to the following probability
\begin{equation}
\label{eq:P^c}
P^c(\tau,x)\coloneqq\frac{\mu(\tau_c^x)}{\mu(\tau_c^0)+\mu(\tau_c^1)}, ~ x \in \Z_2.
\end{equation}
For $c\in C_k(\Z_N^m)$, we let $\tilde{C}_c^k$ be the set of all tuples $(\tau_1,\tau_2) \in C^k \times C^k$ such that $\tau_1(c')=\tau_2(c'),~ \text{for all} ~c' \neq c$, i.e., the set of pairs of configurations that agree except possibly on $c$. \newline 
Now for $a,b \in C_{k}(\Z_N^m)$, we let
\begin{equation}
\label{eq:defrab}
\hat{r}_{a,b} \coloneqq \max_{(\tau_1,\tau_2)\in \tilde{C}_a^k}d_{TV}(P^{b}(\tau_1,\cdot),P^{b}(\tau_2,\cdot)) \le 1.
\end{equation}
This quantity measures the maximal change in the heat-bath distribution at $b$, in total variation distance, when the configuration is modified only at $a$. Hence it measure the influences that the value at $a$ has on the updated value at $b$. Using this, one defines a \textit{dependency matrix} 
\begin{equation}
\label{eq:defdepnmatrix}
R = (r_{a,b})_{a,b \in C_k(\Z_N^m) \times C_k(\Z_N^m)}, ~ \text{where}~ r_{a,b} \ge \hat{r}_{a,b}.
\end{equation}
\newline 
The mathematical statement is the following. 
\begin{Lemma}[Corollary~18 in \cite{RapidMixingGlauber}]
\label{lem:dobrushincond}
    Let $P$ be the transition kernel from Example~\ref{ex:firstexample}(b), and denote by $(X_n)_{n \in \N}$ the discrete Markov chain whose dynamics are given by $P$. If a dependency matrix $R$ as defined in \eqref{eq:defdepnmatrix} satisfies the Dobrushin condition
    \begin{equation}
    \label{eq:dobrushincond}
    \|R\|_{1} \le \delta <1, ~ \text{for some} ~\delta \in (0,1),
    \end{equation}
    then it holds that
    \[
    t_\text{mix}^N(\epsilon) \le \frac{|C_k(\Z_N^m)|}{1-\delta} \log\Big(\frac{|C_k(\Z_N^m)|}{\epsilon}\Big), ~ \epsilon \in (0,1).
    \]
\end{Lemma}
In the following pages, we will prove the following bound for $\hat{r}_{a,b}$. 
\begin{Lemma}
    \label{lem:tvdistance}
Let $m \ge 1$, $k \in \{0,\dots,m-1\}$, and consider the transition kernel $P$ as described in Example \ref{ex:firstexample}(b). Consider $\hat{r}_{a,b}$ for $ a,b \in C_k(\Z_N^m)$ as defined in \eqref{eq:defrab}. Then,
    \[
0 \le \hat{r}_{a,b} \le  \tanh(\beta)\mathbbm{1}\{\exists q \in C_{k+1}(\Z_N^m) \colon a,b \in \partial q, a\neq b\}.    
\]
\end{Lemma}
Note the similarity with the result for the Ising model as in \cite[Observation 11]{Hayes2006} for example.

As a corollary, one deduces that the Dobrushin condition is satisfied for small enough $\beta$.
\begin{Corollary}
    \label{cor:depenmtrx}
Let $m \ge 1$, $k \in \{0,\dots,m-1\}$, and consider the transition kernel $P$ as described in Example \ref{ex:firstexample}(b). If $\beta$ satisfies \eqref{eq:smallbeta}, then there is a dependency matrix $R$ as in \eqref{eq:defdepnmatrix} that satisfies the Dobrushin condition \eqref{eq:dobrushincond}.
\end{Corollary}
\begin{proof}[Proof of Corollary \ref{cor:depenmtrx}]
    The proof follows immediately from Lemma \ref{lem:tvdistance} by defining $R=(r_{a,b})_{a,b \in C_{k}(\Z_N^m) \times C_{k}(\Z_N^m)}$ where
    \[
    r_{a,b} \coloneqq   \tanh(\beta)\mathbbm{1}\{\exists q \in C_{k+1}(\Z_N^m) : a,b \in \partial q, a\neq b\}.
    \]
Then, recall that
\[
\|R\|_{1} = \max_{1 \le b \le |C_k(\Z_N^m)|}\sum_{a=1}^{|C_k(\Z_N^m)|}|r_{a,b}|.
\]
For fixed $b \in C_k(\Z_N^m)$, note that there are $2(2k+1)(m-k)$ $k$-cells $a \in C_k(\Z_N^m)$ that share an upper boundary with $b$ by Lemma \ref{lem:sizegraph}. Hence, one deduces that
\[
\|R\|_{1}  \le 2(2k+1)(m-k)\tanh(\beta).
\]
If $\beta$ satisfies \eqref{eq:smallbeta}, then one can pick 
\[
\delta \coloneqq 2(2k+1)(m-k)\tanh(\beta)<1,
\]
so that the Dobrushin condition is satisfied for $\delta$, which concludes the proof.
\end{proof}
\begin{proof}[Proof of Lemma \ref{lem:tvdistance}]
    We are going to consider different cases for $a,b \in C_k(\Z_N^m)$. 
   
    Let us first assume that $a=b \in C_{k}(\Z_N^m)$. Then, $(\tau_1,\tau_2) \in \tilde{C}_a^k$ agree everywhere except possibly on $a$. As the probability distribution of the value at a $k$-cell given by $P^a$ in \eqref{eq:P^c} doesn't depend on the actual values of $\tau_1(a)$ and $\tau_2(a)$, we conclude that $P^a(\tau_1,\cdot)\overset{d}{=}P^a(\tau_2,\cdot)$, and hence $\hat{r}_{a,b}=0$. \newline 
    Now, let $a \neq b \in C_{k}(\Z_N^m)$. For $(\tau_1,\tau_2)\in \tilde{C}_a^k$, it holds that either $\tau_1=\tau_2$ or $\tau_1=\tau_2+\delta_a$. In the first case one has $d_{TV}(P^{b}(\tau_1,\cdot),P^{b}(\tau_2,\cdot))=0$. 
    Because of this, we now assume $\tau_2=\tau_1+\delta_a$. Note that for $\tau \in C^k$
    \begin{equation}
    \label{eq:ProbmeasureHB}
    d_{\text{TV}}(P^b(\tau,\cdot),P^b(\tau+\delta_a,\cdot))=|P^b(\tau,1)-P^b(\tau+\delta_a,1)|=|P^b(\tau,0)-P^b(\tau+\delta_a,0)|,
    \end{equation}
    where the equality comes from the definition of the total variation distance as $P^b(\tau,\cdot)$ only assigns probability to values in $\Z_2$.
    Now, note that the difference between $d\tau$ and $d(\tau_b^x)$ only affects the $(k+1)$-cells adjacent to $b$; hence, computing the Hamiltonian for $\tau_b^x$, one obtains
    \begin{align*}
    \mathcal{H}(\tau_b^x)&= -\sum_{q' \notin \hat{\partial b}}\rho(d\tau(q'))
    -\sum_{q \in \hat{\partial}b}\rho(d\tau_b^x(q))\\&= -\sum_{q' \notin \hat{\partial b}}\rho(d\tau(q'))-\rho(x)\sum_{q \in \hat{\partial}b}\prod_{c \in \partial q \setminus \{b\}}\rho(\tau(c)) =:  -A_b(\tau)-\rho(x)h_b(\tau).
    \end{align*}
From this we deduce that the probability of obtaining the configuration $\tau_b^x$ equals
\begin{align*}
P^b(\tau,x)&=\frac{\mu(\tau_b^x)}{\mu(\tau_b^0)+\mu(\tau_b^1)}=\frac{\myexp(\beta(A_b(\tau)+\rho(x)h_b(\tau)))}{\myexp(\beta(A_b(\tau)+h_b(\tau)))+\myexp(\beta(A_b(\tau)-h_b(\tau)))} \\ 
&= \frac{\myexp(\beta\rho(x)h_b(\tau))}{\myexp(-\beta h_b(\tau))+\myexp(\beta h_b(\tau))} =\frac{1+\rho(x)\tanh(\beta h_b(\tau)))}{2}.
\end{align*}
The last equality follows from the hyperbolic identities. 
This probability inserted in \eqref{eq:ProbmeasureHB} gives
\begin{align*}
    d_{TV}(P^{b}(\tau,\cdot),P^{b}(\tau+\delta_a,\cdot)) = \frac{1}{2}|\tanh(\beta h_b(\tau))-\tanh(\beta h_b(\tau+\delta_a))|.
\end{align*}
If $\hat{\partial}a\cap\hat{\partial}b= \emptyset$, then by definition $h_b(\tau_1)=h_b(\tau_2)$ and $\hat{r}_{a,b}=0$. Hence, let's assume that there exists $q_{a,b} \in C_{k+1}(\Z_N^m)$ such that $a,b \in \partial q_{a,b}$. 
Then, the difference between $h_b(\tau)$ and $h_b(\tau+\delta_a)$ comes solely from $q_{a,b}$.
There is a unique such $(k+1)$-cell as $N\ge3$, so that two $k$-cells share at most one upper boundary.
More precisely one has
\[
h_b(\tau)=\sum_{q \in \hat{\partial}b\setminus\{q_{a,b}\}}\prod_{c' \in \partial q \setminus \{b\}}\rho(\tau(c')) + \prod_{c \in \partial q_{a,b}\setminus\{b\}}\rho(\tau(c))\eqqcolon u+v.
\]
Then, one has that $d(\tau+\delta_a)(q)=d\tau(q)$ for $q' \in \hat{\partial}b \setminus\{q_{a,b}\}$ and $d(\tau+\delta_a)(q_{a,b}) \neq d\tau(q_{a,b})$. Hence, one deduces that
\[
h_b(\tau+\delta_a)=u-v.
\] 
Altogether, taking into account that $v=\pm1$ we deduce by the hyperbolic identities that
\begin{align*}
    &\frac{|\tanh(\beta h_b(\tau))-\tanh(\beta h_b(\tau+\delta_a))|}{2}=\frac{\tanh(\beta (u+1))-\tanh(\beta (u-1))}{2} \\
    &= \frac{\sinh(2\beta)}{\cosh(2\beta u)+\cosh(2\beta)} \le \frac{\sinh(2\beta)}{1+\cosh(2\beta)}=\frac{2\sinh(\beta)\cosh(\beta)}{2\cosh^2(\beta)}=\tanh(\beta).
\end{align*}
The first equality follows directly from the representation $h_b(\tau)=u\pm1$, and the fact that $\tanh$ is increasing. The second equality follows from the hyperbolic identity
\begin{align*}
\tanh(x)-\tanh(y)=\frac{\sinh(x-y)}{\cosh(x)\cosh(y)}=\frac{2\sinh(x-y)}{\cosh(x+y)+\cosh(x-y)}
\end{align*}
by taking $x=\beta (u+1)$ and $y=\beta (u-1)$. Then, the inequality follows from noting that $\cosh(x)\ge 1$ for all $x \in \R$. The third equality comes from hyperbolic identities, while the last equality is the definition of $\tanh$. 
Altogether, we have proved 
\[
0 \le \hat{r}_{a,b} \le  \tanh(\beta)\mathbbm{1}\{\exists q \in C_{k+1}(\Z_N^m) : a,b \in \partial q,~ a\neq b\},
\]
which proves the claim.
\end{proof}
We have now proven all ingredients needed for the following proof of Proposition~\ref{prop:smallbeta}.
\begin{proof}[Proof of Proposition \ref{prop:smallbeta}]
Note that by the definition of the generator in \eqref{eq:defgenL1}
\[
\mathcal{L}=|C_k(\Z_N^m)|(P-I),
\]
as each $k$-cell has a rate-one exponential clock. Hence, averaging the estimate of Lemma~\ref{lem:dobrushincond} over the Poisson number of updates, with total update rate equal to the number of $k$-cells $|C_k(\Z_N^m)|$, gives the following upper bound for the mixing time for $\epsilon \in (0,1)$
\begin{align*}
t_{\text{mix}}^N(\epsilon) &\le \frac{1}{1-\delta} \log\Big(\frac{|C_k(\Z_N^m)|}{\epsilon}\Big)\\
&=\frac{1}{1-\delta} \Big(\log\Big(\binom{m}{k}\Big)+m\log(N)-\log(\epsilon) \Big) \le C\log(N),
\end{align*}
for $C$ depending on $m,k,\beta$ and $\epsilon$. The value for $|C_k(\Z_N^m)|$ comes from Lemma \ref{lem:sizegraph}. This concludes the proof.
\end{proof}
\begin{Remark} \begin{enumerate}
    \item We observe that the space being a torus does not play a role in the above proof. With slightly different estimates, notably from Lemma \ref{lem:sizegraph}, the exact same proof works for other spaces such as the cubic box with no boundary condition.
    \item The proof of Proposition \ref{prop:smallbeta} allows the dimension of $k$-cells to be $k=m-1$, unlike in the bottleneck proof, detailed in the next section, where the highest dimension is $k=m-2$.
    \item Note that the proof of Proposition \ref{prop:smallbeta} explicitly takes into account the heat-bath rule. If one consider the Metropolis rule, then there would be a dependence on the value of a $k$-form at the cell that might be flipped, which would induce positive diagonal terms in $R$. This is explicitly not the case in the heat-bath rule as one samples the value at the updated cell independently of the value of the previous state at that cell.
\end{enumerate}
\end{Remark}
\section{Slow mixing at low temperature}
\label{sec:slowmixing}
In this section, we will prove that the state space $C^k(\Z_N^m,\Z_2)$ can be divided into \textit{bottles} that carry most of the mass under $\mu_{\beta,N,m,k}$ when $\beta$ is large enough, and \textit{bottlenecks} that connect these bottles, having a small mass under the probability measure. Our goal is to use a bottleneck argument on these bottles to prove the lower bound on the mixing time by applying Corollary \ref{cor:continuousdynamics}.
\subsection{Definitions: Vortex cluster and decoder}
To improve readability, recall that we defined by $\ell \coloneqq m-k \ge 2$ as the codimension of the model. Furthermore, recall that when it is clear from the context, we will write $C^k$ for $C^k(\Z_N^m,\Z_2)$ and $\mu$ for $\mu_{\beta,N,m,k}$. We also write $C^k_0$ for the set of closed $k$-forms and denote $d$ instead of $d_k$ when the dimension is clear. Now, we introduce some new definitions. 
For two distinct $k$-cells $c,c'$ we write $c \sim c'$ if there exists $q \in C_{k+1}(\Z_N^m)$ such that $c,c' \in \partial q$, which is equivalent to $c, c'$ being connected by an edge in the graph $\mathcal{G}_m^k(N)$.
\begin{Def}
\label{Def:vortexclusters}
    Let $\omega \in C_0^{k+1}$ be a closed $(k+1)$-form. For $c \in \supp \omega$ we define the vortex cluster of $c$ as the set
    \[
    \{c' \in \supp \omega: \exists c_1,\dots,c_n \in \supp \omega ~~\text{s.t.} ~ c \sim c_1, c_1 \sim c_2, \dots, c_n \sim c'\}.
    \]
Recalling the graph $\mathcal{G}_m^{k+1}(N)$, the vortex cluster of $c \in \supp \omega$ corresponds to the connected component of $\supp \omega$ in which $c$ lies. We also call $\omega \in C_0^{k+1}$ a vortex cluster if it induces a single connected component in $\mathcal{G}_m^{k+1}(N)$, and denote the set of all vortex clusters by $\mathcal{V}^{k+1}$. \newline
    Note that each $(k+1)$-cell $c$ is either in a specific vortex cluster or in none. We denote the number of distinct vortex clusters for $\omega \in C^{k+1}_0$ by $I(\omega)$. \newline 
    If we enumerate the vortex clusters of $\omega$ according to the lexicographic order, then for $i \in [I(\omega)]$, we denote by $\omega_i$ the $(k+1)$-form that has as support only the $i$-th vortex cluster, hence $I(\omega_i)=1$.
\end{Def}
Using the notation above and the cell-by-cell sum introduced in \eqref{eq:notaddition}, one has that
\[
\omega=\sum_{i=1}^{I(\omega)} \omega_i,
\]
where the $\omega_i$ are the vortex clusters in $\supp \omega$. 
Observe that each $\omega_i$ is closed. To see it, let $Q \in C_{k+2}(\Z_N^m)$, then
\[
0=d_{k+1}\omega(Q)=\sum_{i=1}^{I(\omega)}d_{k+1}\omega_i(Q)=\sum_{i=1}^{I(\omega)}\omega_i(\partial_{k+2}Q).
\]
The last equality follows from Theorem \ref{thm:stokes}. Now, recall that by construction, all faces of $Q$ that belong to $\supp \omega$ lie in the same connected component, i.e., there exists at most one $j \in [I(\omega)]$ such that $\supp \partial_{k+2}Q \cap \supp \omega_j \neq \emptyset$. In either case it implies $\omega_i(\partial_{k+2}Q)=0$ for all $i \neq j$ and hence $\omega_j(\partial_{k+2}Q)=0$.

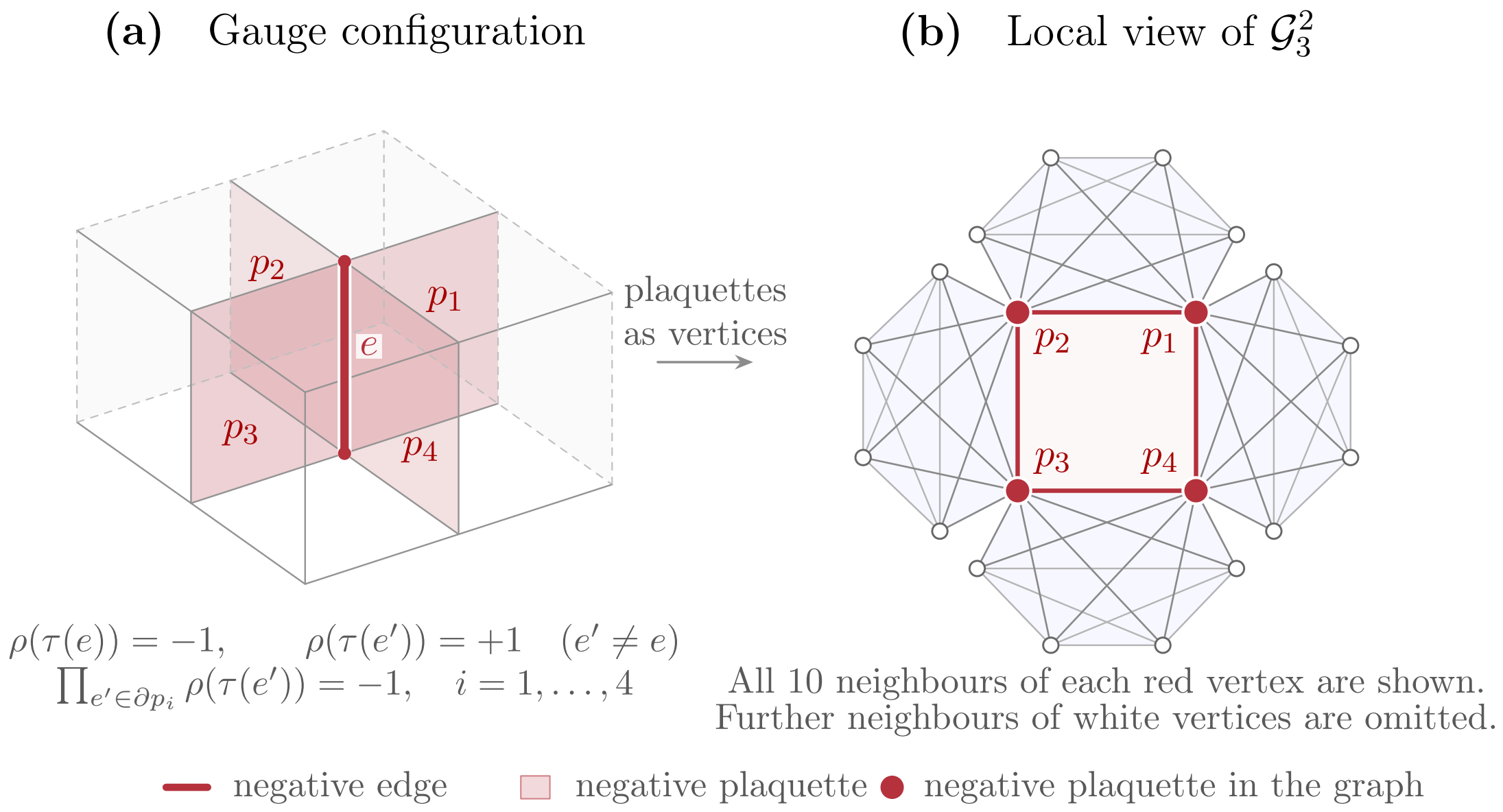
\begin{figure}
    \centering
\begin{tikzpicture}[
  font=\small,
  line cap=round,line join=round,
  lattice/.style={draw=black!42,line width=.45pt},
  rear/.style={draw=black!25,line width=.4pt,dash pattern=on 2pt off 2pt},
  g edge/.style={draw=black!29,line width=.45pt},
  g active/.style={draw=black!46,line width=.5pt},
  positive/.style={circle,draw=black!60,fill=white,line width=.55pt,
    minimum size=4.2pt,inner sep=0pt},
  negative/.style={circle,draw=white,fill=red!75!black,line width=.8pt,
    minimum size=7.2pt,inner sep=0pt},
  plabel/.style={text=red!65!black,font=\small,inner sep=1pt},
  panel/.style={font=\normalsize},
  note/.style={font=\footnotesize,text=black!65},
  cube patch/.style={fill=blue!3,draw=blue!12,line width=.4pt}
]
\definecolor{defect}{RGB}{181,49,60}

\node[panel] at (0,3.65) {\textbf{(a)}\quad Gauge configuration};

\begin{scope}[shift={(0,-.40)},
  x={(1.48cm,.48cm)},y={(-1.10cm,.78cm)},z={(0cm,1.85cm)}]

  \fill[black!2] (-1,1,0)--(1,1,0)--(1,1,1)--(-1,1,1)--cycle;
  \fill[black!2] (1,-1,0)--(1,1,0)--(1,1,1)--(1,-1,1)--cycle;

  \fill[defect,fill opacity=.19]
    (0,0,0)--(1,0,0)--(1,0,1)--(0,0,1)--cycle;
  \fill[defect,fill opacity=.12]
    (0,0,0)--(0,1,0)--(0,1,1)--(0,0,1)--cycle;
  \fill[defect,fill opacity=.22]
    (0,0,0)--(-1,0,0)--(-1,0,1)--(0,0,1)--cycle;
  \fill[defect,fill opacity=.15]
    (0,0,0)--(0,-1,0)--(0,-1,1)--(0,0,1)--cycle;

  \foreach \h in {0,1}{
    \draw[rear] (-1,1,\h)--(1,1,\h)--(1,-1,\h);
    \draw[lattice] (-1,1,\h)--(-1,-1,\h)--(1,-1,\h);
    \draw[lattice] (-1,0,\h)--(1,0,\h);
    \draw[lattice] (0,-1,\h)--(0,1,\h);
  }
  \foreach \a in {-1,0,1}{
    \foreach \b in {-1,0,1}{
      \ifnum\a=1
        \draw[rear] (\a,\b,0)--(\a,\b,1);
      \else
        \ifnum\b=1
          \draw[rear] (\a,\b,0)--(\a,\b,1);
        \else
          \draw[lattice] (\a,\b,0)--(\a,\b,1);
        \fi
      \fi
    }
  }

  \draw[white,line width=3.8pt] (0,0,0)--(0,0,1);
  \draw[defect,line width=2.3pt] (0,0,0)--(0,0,1);
  \fill[defect] (0,0,0) circle[radius=1.8pt];
  \fill[defect] (0,0,1) circle[radius=1.8pt];
  \node[text=defect,fill=white,fill opacity=.85,text opacity=1,
    inner sep=1pt,anchor=west] at (.07,0,.54) {$e$};

  \node[plabel] at (.66,0,.63) {$p_1$};
  \node[plabel] at (0,.67,.68) {$p_2$};
  \node[plabel] at (-.67,0,.28) {$p_3$};
  \node[plabel] at (0,-.67,.30) {$p_4$};
\end{scope}

\node[note] at (0,-2.22)
  {$\rho(\tau(e))=-1,\qquad \rho(\tau(e'))=+1\quad(e'\ne e)$};
\node[note] at (0,-2.67)
  {$\prod_{e'\in\partial p_i}\rho(\tau(e'))=-1,
    \quad i=1,\ldots,4$};

\draw[-{Stealth[length=5pt,width=4pt]},black!45,line width=.65pt]
  (3.02,.48)--(3.94,.48);
\node[note,align=center] at (3.48,.98) {plaquettes\\as vertices};

\begin{scope}[shift={(7.35,.10)}]
  \node[panel] at (0,3.55)
    {\textbf{(b)}\quad Local view of $\mathcal G_3^2$};

  \coordinate (p1) at (.86,.86);
  \coordinate (p2) at (-.86,.86);
  \coordinate (p3) at (-.86,-.86);
  \coordinate (p4) at (.86,-.86);

  \foreach \j/\angle in {1/0,2/90,3/180,4/270}{
    \begin{scope}[rotate=\angle]
      \coordinate (q\j a) at (-1.25,1.61);
      \coordinate (q\j b) at (-.54,2.35);
      \coordinate (q\j c) at (.54,2.35);
      \coordinate (q\j d) at (1.25,1.61);
      \path[cube patch] (-.86,.86)--(-1.25,1.61)--(-.54,2.35)
        --(.54,2.35)--(1.25,1.61)--(.86,.86)--cycle;
    \end{scope}
  }

  \foreach \j/\u/\v in {1/1/2,2/2/3,3/3/4,4/4/1}{
    \foreach \a/\b in {a/b,a/c,a/d,b/c,b/d,c/d}
      \draw[g edge] (q\j\a)--(q\j\b);
    \foreach \a in {a,b,c,d}{
      \draw[g active] (p\u)--(q\j\a);
      \draw[g active] (p\v)--(q\j\a);
    }
  }

  \fill[defect,fill opacity=.035] (p1)--(p2)--(p3)--(p4)--cycle;
  \draw[defect,line width=1.2pt] (p1)--(p2)--(p3)--(p4)--cycle;
  \foreach \j in {1,2,3,4}
    \foreach \a in {a,b,c,d} \node[positive] at (q\j\a) {};
  \foreach \i in {1,2,3,4} \node[negative,fill=defect] at (p\i) {};

  \node[plabel,anchor=north east] at (.74,.74) {$p_1$};
  \node[plabel,anchor=north west] at (-.74,.74) {$p_2$};
  \node[plabel,anchor=south west] at (-.74,-.74) {$p_3$};
  \node[plabel,anchor=south east] at (.74,-.74) {$p_4$};

  \node[note] at (0,-2.72) {All $10$ neighbours of each red vertex are shown.};
  \node[note] at (0,-3.07) {Further neighbours of white vertices are omitted.};
\end{scope}

\draw[defect,line width=2pt] (-1.72,-3.62)--(-1.32,-3.62);
\node[anchor=west,note] at (-1.22,-3.62) {negative edge};
\path[draw=defect!60,fill=defect!18,line width=.4pt]
  (1.70,-3.73) rectangle (1.98,-3.51);
\node[anchor=west,note] at (2.08,-3.62) {negative plaquette};
\node[negative,fill=defect] at (5.28,-3.62) {};
\node[anchor=west,note] at (5.44,-3.62) {negative plaquette in the graph};
\end{tikzpicture}
\caption{Example of a local view of the graph $\mathcal{G}_m^{k+1}$ for $m=3, k=1$. Here the gauge configuration $\tau$ is supported on a single edge $e$.}
\label{fig:degofgraph}
\end{figure}
Now, we define the value that divides the set of connected vortex clusters into \textit{small} and \textit{large} clusters. 
In order for the following value to be well-defined, recall that the block sequence $(\mathcal{B}_N)_{N \in \N}$ satisfies \eqref{eq:conditionblock}. It  means that there exists $N_0$ depending on $(\mathcal{B}_N)_{N \in \N}$, such that for every $N \ge N_0$ one has $\lambda(\mathcal{B}_N) \le \frac{N^\ell}{2}$. Now, for $N \ge N_0$, we let 
\begin{equation}
\label{eq:defg_n}
g_{\ell}(N)\coloneqq  \min\Big(N^{\ell-1},\Big(\frac{N^\ell-\lambda(\mathcal{B}_N)}{4\hat{C}_{m,k}\ell(\ell-1) \lambda(\mathcal{B}_N)}\Big)^{\frac{\ell-1}{\ell}}\Big),
\end{equation}
for a sequence of blocks $(\mathcal{B}_N)_{N \in \N}$ satisfying conditions \eqref{eq:condonB_N} and \eqref{eq:conditionblock}. Note that $g_\ell(N)$ depends on the dimension $m$, the dimension of the $k$-cells, the size of the torus $N$ and the block sequence used in the dynamics. However, we will just write $g_{\ell}(N)$ to simplify notation.

Throughout the remainder of this section,
we assume that $N \ge N_0((\mathcal{B}_N)_{N \in \N})$. The case for smaller values of $N$ will only be considered in the final proof. 

We now introduce the following state space, which contains all $k$-forms $\tau \in C^k$ such that each vortex cluster of its discrete derivative $d\tau$ has a support bounded by $g_\ell(N)$,
\begin{equation}
\label{eq:definitionOmega}
\Omega\coloneqq \{\tau \in C^k: d\tau = \sum_{i=1}^{I(d\tau)} \omega_i ~ \text{s.t.} ~ \forall i \colon |\supp \omega_i| < g_\ell(N)\}.
\end{equation}
We call $\Omega$ the set of \textit{bottles}. \newline 
Now, let
\begin{equation}
\label{eq:defbigW}
\mathcal{W}\coloneqq d(\Omega) =\{\omega \in C^{k+1}_0 \colon \omega=d\tau, \tau \in \Omega\}.
\end{equation}
Finally, we let
\[
\mathcal{V}_{\text{small}}^{k+1} \coloneqq \{\omega \in \mathcal{V}^{k+1}: |\supp \omega| < g_\ell(N) \}.
\]
By Lemma \ref{lem:dec2}, for each closed vortex cluster $\omega \in C^{k+1}_0$ such that $|\supp \omega| < g_\ell(N)$, there exists a $k$-form, which we denote $\tau$, that satisfies the support bound \eqref{eq:boundonsupport}, such that $d\tau=\omega$. For such a closed $(k+1)$-form, we define a mapping denoted $T : \mathcal{V}^{k+1}_{\text{small}} \to C^k$ which maps $\omega$ into a $k$-form such that 
\[
T(\omega)\in \argmin_{\tau:d\tau=\omega}|\supp \tau|.
\]
Whenever there is no unique choice for $T(\omega)$ we use the lexicographic order to uniquely fix the value of $T(\omega)$. Note that $T(\omega)$ satisfies the support bound \eqref{eq:boundonsupport}. 
\begin{Def}
    Let $\tau \in \Omega$. Then, one can rewrite $d\tau \in C^{k+1}$ as follows
    \[
    d\tau=\sum_{i=1}^{I(d\tau)} \omega_i, ~ \omega_i \in \mathcal{V}^{k+1}, i \in [I(d\tau)].
    \]
    Then, define the decoder as the mapping $\Delta : \Omega \to C^k$ 
\begin{equation}
\label{eq:defdecoder}
\Delta\tau \coloneqq \tau + \sum_{i=1}^{I(d\tau)} T(\omega_i).
\end{equation}
\end{Def}
\begin{Prop}
    Let $\tau \in \Omega$. Then $\Delta\tau$ is a $k$-form such that
    \[
    d(\Delta\tau)=0.
    \]
\end{Prop}
\begin{proof}
Note that by definition one has that $dT(\omega_i)=\omega_i$. In particular it implies
\[
d(\Delta\tau)=d\Big(\tau+\sum_{i=1}^{I(d\tau)} T(\omega_i)\Big) = d\tau +\sum_{i=1}^{I(d\tau)} dT(\omega_i)=d\tau+\sum_{i=1}^{I(d\tau)} \omega_i=d\tau+d\tau=0,
\]
which gives the claim.
\end{proof}
As the decoder maps $\tau \in \Omega$ to a closed $k$-form, we can intuitively see any configuration $\tau \in \Omega$ as a random noise $\tau_0$ with no defects (meaning $d\tau_0=0)$ perturbed by \textit{small} cluster of $(k+1)$-cells. The decoder is a useful way of mapping a configuration $\tau$ to a closed configuration $\tau_0$ to which one can apply the theory from Section \ref{sec:dec}. As we will see, the cohomology of $\Delta \tau=\tau_0$ will be relevant for our next arguments.
\subsection{Cohomology sectors}
Recall that by Lemma \ref{lem:numbercohomologyclasses}, there are $2^{\binom{m}{k}} =: \vartheta$ cohomology classes in $H^k$. As each $\Delta\tau$ is a closed $k$-form, each of them is in one of the $\vartheta$ cohomology classes in $H^{k}$. We can find some ordering, so that $(H_i^{k})_{i=1}^{\vartheta}$ is the set of the different cohomology classes.
Hence, we can define the following partition of $(\Omega_i)_{i \in [\vartheta]}$ of $\Omega$, where for each $i \in [\vartheta]$
\begin{equation}
\label{eq:omegaplusandomegaminus}
\Omega_i\coloneqq  \{\tau \in \Omega \colon \Delta\tau \in H_i^{k}\}.
\end{equation}
We define $H_1^k$ as being the null-cohomology class.
\begin{Ex}
 Let $\tau_1 \equiv 0 \in C^k$. Then $\Delta\tau_1=\tau_1$, and $\tau_1$ would be in the first cohomology class $H_1^k$, which is the null-homologous class, hence $\tau_1 \in \Omega_1$. Let's now assume that $\tau_2=\delta_c$, for some $c \in C_k(\Z_N^m)$ and that $N$ is large enough for $\tau_2$ to be in~$\Omega$. Then, $\Delta \tau_2 \equiv 0$ and hence $\Delta \tau_2 \in H_1^k \implies \tau_2 \in \Omega_1$. \newline
    If $\tau_3$ has all $k$-cells in its support, meaning $\tau_3 \equiv 1 \in C^k$, then $\Delta \tau_3 \equiv 1$ as well as $d\tau_3 \equiv0 $. Now let $(\Gamma_I)_{I\subseteq[m],\,|I|=k}$ denote the standard
coordinate $k$-cycles, as defined in \eqref{eq:defGamma_IwithC_I}, whose homology classes form a basis
of $H_k$, as in the proof of Lemma \ref{lem:homologygroup}.
Each $\Gamma_I$ contains $N^{k}$ cells, and hence
\[
    \tau_3(\Gamma_I)=N^{k} \pmod{2}.
\]
Hence if we assume that $k \ge 1$, then we have the following.
If $N$ is even, these evaluations all vanish. Since the classes $[\Gamma_I]$ form a basis of $H_k$, it follows from Corollary~\ref{cor:cohomologyclasses} that $\tau_3$ represents the null-cohomology class, and thus $\tau_3\in\Omega_1$. 
If $N$ is odd, these evaluations equal $1$, so $\tau_3$ lies in a non-trivial cohomology class by Corollary~\ref{cor:cohomologyclasses} and belongs to $\Omega_j$ for some $j\neq 1$.
In the case in which $k=0$, then observe that $\tau_3$ always lies in the non-trivial cohomology class (there are only two classes if $k=0$) by the above arguments, hence $\tau_3 \in \Omega_2$. This corresponds to the two ground states of the Ising model.
\end{Ex}

Now, we have the following result, concerning the probability mass of the sets $\Omega_i$.
\begin{Lemma}
\label{lem:samemass}
    Let $\Omega_i$ and $\Omega_j$, $i \neq j$, be defined as in \eqref{eq:omegaplusandomegaminus}. Then,
    \[
    \mu_{\beta,N,m,k}(\Omega_i)= \mu_{\beta,N,m,k}(\Omega_j).
    \]
\end{Lemma}
\begin{proof}
We will prove that there is a bijection $\phi$ from $\Omega_i$ to $\Omega_j$ such that $\mu(\phi(\tau_i))=\mu(\tau_i)$ for any $\tau_i \in \Omega_i$ where $\mu \coloneqq \mu_{\beta,N,m,k}$. \newline 
To this end, note that all cohomology classes contain by definition the same number of elements. Namely, pick any representative $\tau_i \in \ker d_k$ such that $\tau_i \in H_i^k$. Then, the set of all $k$-forms belonging to its cohomology class is
\[
\{\tau_i +d_{k-1}\eta \colon \eta \in C^{k-1}\}=H_i^k,
\]
and $H_i^k$ has cardinality $|\text{Im}~d_{k-1}|$. \newline 
Now, pick representatives $\tau_i \in H_i^k, \tau_j \in H_j^k$. Then, define
\[
\phi(\tau_i +d_{k-1}\eta) \coloneqq \tau_j+ d_{k-1}\eta.
\]
By the above, this is a bijection from $H_i^k$ to $H_j^k$. 
We now extend the bijection to $\Omega_i$ as follows. Let $\tau \in \Omega_i$ be such that 
\[
d\tau=\sum_{l=1}^{I(d\tau)}\omega_l.
\]
Then, using the definition of the decoder $\Delta$ recall that $
\tau=\Delta \tau+\sum_{l=1}^{I(d\tau)}T(\omega_l).
$
Now, define $\phi$ on $\Omega_i$ as follows
\[
\phi(\tau)\coloneqq\phi(\Delta \tau)+\sum_{l=1}^{I(d\tau)}T(\omega_l).
\]
In particular, $\tau$ and $\phi(\tau)$ have exactly the same $(k+1)$-defects, hence $\phi(\tau) \in \Omega$. \newline 
\textit{Claim.} The mapping $\phi: \Omega_i \to \Omega_j$ is well-defined and bijective. \newline 
\textit{Proof of the claim.} First we prove that $\phi(\tau)\in \Omega_j$ for $\tau \in \Omega_i$. To see this, one has to prove that $\Delta\phi(\tau) \in H_j^k$, hence one obtains
\begin{align*}
&\Delta \phi(\tau)=\Delta \Big(\phi(\Delta \tau)+\sum_{l=1}^{I(d\tau)}T(\omega_l) \Big) = \phi(\Delta \tau)+\sum_{l=1}^{I(d\tau)}T(\omega_l)+\sum_{l=1}^{I(d\tau)}T(\omega_l) =\phi(\Delta \tau) \in H_j^k,
\end{align*}
where the second equality comes from $\phi(\Delta\tau)$ being closed meaning $\Delta\phi(\Delta\tau)=\phi(\Delta\tau)$, and the last equality from $\Delta\tau \in H_i^k$. \newline 
To see that $\phi$ is injective, note first
\[
\phi(\tau)=\phi(\tau') \implies \phi(\Delta \tau)=\phi(\Delta \tau')\implies \Delta \tau=\Delta \tau'. 
\]
The first implication follows from the definition of $\phi$: $\phi(\tau)=\phi(\tau')$ implies that $d\tau=d\tau'$, and hence their decoded versions $\phi(\Delta \tau)$ and $\phi(\Delta \tau')$ must be equal. In particular $\Delta \tau=\Delta \tau'$ as $\phi$ is bijective on $H_i^k$. From this follows:
\[
\tau=\Delta \tau+\sum_{l=1}^{I(d\tau)}T(\omega_l)=\Delta \tau'+\sum_{l=1}^{I(d\tau)}T(\omega_l)=\tau'.
\]
This proves that $\phi$ is injective. 
To see that $\phi$ is surjective note that for any $\tau' = \Delta \tau' + \sum_{i=1}^{I(d\tau')} T(\omega_i)\in \Omega_j$, one defines 
\[
\tau\coloneqq\phi^{-1}(\Delta \tau')+\sum_{i=1}^{I(d\tau')} T(\omega_i) \in \Omega_i.
\]
By definition, one has $\phi(\tau)=\tau'$, which concludes the proof of the claim.
\newline 
Finally as $\phi(\tau)$ has exactly the same $(k+1)$-defects as $\tau$, it implies that $\mu(\phi(\tau))=\mu(\tau)$. Altogether one has
\[
\mu(\Omega_i)=\sum_{\tau \in \Omega_i}\mu(\tau)=\sum_{\tau \in \Omega_i}\mu(\phi(\tau))=\sum_{\tau' \in \Omega_j}\mu(\tau')=\mu(\Omega_j).
\]
This concludes the proof.
\end{proof}
Then, as each cohomology class is non-empty, each bottle is non-empty and furthermore from Lemma \ref{lem:samemass} one deduces \[
\mu(\Omega_i)=\frac{\mu(\Omega)}{\vartheta} \le \frac{1}{\vartheta} \le \frac{1}{2}.
\]
Hence, we deduce the following corollary.
\begin{Corollary}
\label{cor:probaofasetomegai}
    For any of the parameters $m,N,\beta,k$ it holds that for $i \in [\vartheta]$ the measure of the bottle $\Omega_i$ under $\mu$ satisfies
    \[
    \mu(\Omega_i) \in (0,\frac{1}{2}].
    \]
\end{Corollary}
\subsection{Bottleneck estimate}
Now, we define the following set, containing the configurations having at least one large vortex cluster.
\[
\Omega^L \coloneqq C^{k}\setminus\Omega=  \{\tau \in C^k(\Z_N^m,\Z_2) : d\tau = \sum_{i=1}^{I(d\tau)} \omega_i ~ \text{s.t.} ~ \exists i : |\supp \omega_i| \ge g_\ell(N)\}.
\]
As above, one defines 
\[
\mathcal{W}^L \coloneqq d(\Omega^L) = \{\omega \in C^{k+1}:\omega=d\tau, \tau \in \Omega^L\}.
\]
Lastly, we let
\[
\mathcal{Z}^L\coloneqq\{\omega \in C^{k+1}_0 \setminus \{0\} \colon \omega=\sum_{i=1}^{I(\omega)}\omega_i, ~ \forall i\in [I(\omega)] \colon |\supp \omega_i| \ge g_\ell(N)\}.
\]
We now present the following relations between these different sets. Recall the definition of $\Omega$, and $\mathcal{W}$ in \eqref{eq:definitionOmega} and \eqref{eq:defbigW}, respectively.
\begin{Lemma}
    \label{lem:thesets}
    Let $N \in \N_{\ge 3}$, $m \in \mathbb{N}_{\ge 2}$, and $k \in \{0,\dots,m-2\}$, and let $\mu$ be the probability measure of hyperlattice gauge theory with these parameters. Then, the following holds
    \begin{enumerate}[label=(\alph*)]
        \item One has
        \[
        |\ker d_k||\mathcal{W}|=|\Omega| ~ \text{and} ~ |\ker d_k||\mathcal{W}^L|=|\Omega^L|.
        \]
        \item It holds that
        \[
        |\ker d_k|\sum_{\omega \in \mathcal{W}}e^{-2\beta|\omega|}=\sum_{\tau \in \Omega}e^{-2\beta|d\tau|} \quad \text{and} \quad |\ker d_k|\sum_{\omega \in \mathcal{W}^L}e^{-2\beta|\omega|}=\sum_{\tau \in \Omega^L}e^{-2\beta|d\tau|}.
        \]
        \item There is an injection $\psi:\mathcal{W}^L \to \mathcal{W} \times \mathcal{Z}^L$.
        \item It holds that
        \[
        \sum_{\omega \in \mathcal{W}^L}\myexp(-2\beta|\omega|) \le \sum_{\omega_1 \in \mathcal{W}}\myexp(-2\beta|\omega_1|)\sum_{\omega_2 \in \mathcal{Z}^L}\myexp(-2\beta|\omega_2|).
        \]
    \end{enumerate}
\end{Lemma}
\begin{proof}
\begin{enumerate}[label=(\alph*)]
    \item Let $\omega  \in \mathcal{W}$. Then, by definition, there is $\tau$ such that $d\tau=\omega$. Now note that there are exactly $|\ker d_k|$ such $\tau$ with $\omega=d\tau$. Namely, for each $\varphi \in \ker d_k$ it holds that
    \[
    d_k(\tau+\varphi)=d_k\tau=\omega.
    \]
    Conversely, if $\varphi \notin \ker d_k$ one deduces that $d_k(\tau+\varphi) \neq d_k\tau = \omega$. 
    Hence, each element in $\mathcal{W}$ can be mapped to exactly $|\ker d_k|$ elements in $\Omega$. \newline 
    The proof for $\mathcal{W}^L$ is exactly the same. 
    \item From the definition of $\mathcal{W}$ for the first equality and $(a)$ for the second equality, one obtains that
\begin{align*}
    \sum_{\tau \in \Omega}e^{-2\beta|d\tau|}=\sum_{\omega \in \mathcal{W}}\sum_{\tau:\omega=d\tau}e^{-2\beta|\omega|}=|\ker d_k|\sum_{\omega \in \mathcal{W}}e^{-2\beta|\omega|}.
\end{align*}
The proof for $\mathcal{W}^L$ is the same and hence omitted here.
\item Let $\omega=\sum_{i=1}^{I(\omega)} \omega_i \in \mathcal{W}^L$. Let $I \coloneqq\{i:|\omega_i| \ge g_\ell(N)\} \subset [I(\omega)]$ and $J\coloneqq[I(\omega)] \setminus I$. Then, let $\psi:\mathcal{W}^L \to \mathcal{W} \times \mathcal{Z}^L$ be defined by
    \[
    \psi(\omega)\coloneqq\Big(\sum_{j \in J}\omega_j, \sum_{i \in I}\omega_i\Big).
    \]
    By construction, note that $\omega'\coloneqq\sum_{j \in J}\omega_j$ consists of vortex clusters whose support is bounded by $g_\ell(N)$. Hence by Lemma \ref{lem:cohomology1}, for $j \in J$, there exists $\tau_j \in C^k$ such that $d\tau_j=\omega_j$. Letting $\tau'\coloneqq\sum_{j \in J} \tau_j \in C^k$, we have  $d\tau'=\omega'$, and hence $\omega' \in \mathcal{W}$. 
Then, by definition, it is immediate that $\sum_{i \in I}\omega_i \in \mathcal{Z}^L$ as all vortex clusters have  support size greater than or equal to $g_\ell(N)$. Hence, $\psi(\omega) \in \mathcal{W} \times \mathcal{Z}^L$ is well-defined. 
Moreover, this function is injective. To see this let $\omega, \tilde{\omega} \in \mathcal{W}^L$ be such that $\psi(\omega)=\psi(\tilde{\omega} )$, and we write $\psi=(\psi_1,\psi_2)$. Then, since
\[
\omega=\psi_1(\omega)+\psi_2(\omega)
\]
one concludes that 
\[
\omega=\psi_1(\omega)+\psi_2(\omega)=\psi_1(\tilde{\omega})+\psi_2(\tilde{\omega})=\tilde{\omega} .
\]
This concludes the proof.
\item Using (c) and the function $\psi$ defined in its proof, one obtains that
\begin{align*}
    &\sum_{\omega \in \mathcal{W}^L}\myexp(-2\beta|\omega|) =  \sum_{\omega \in \mathcal{W}^L}\myexp(-2\beta|\psi_1(\omega)+\psi_2(\omega)|) \\ = &\sum_{\substack{(\omega_1,\omega_2)\in \psi(\mathcal{W}^L)}} \myexp(-2\beta|\omega_1|) \myexp(-2\beta|\omega_2|) \\ \le &\sum_{\substack{\omega_1 \in \mathcal{W} \\ \omega_2 \in \mathcal{Z}^L}} \myexp(-2\beta|\omega_1|) \myexp(-2\beta|\omega_2|),
\end{align*}
where the second equality follows from $\psi_1(\omega)$ and $\psi_2(\omega)$ being disjoint by definition and the inequality follows from injectivity of $\psi$.
\end{enumerate}
\end{proof}
\begin{figure}[htbp]
\centering
\begin{tikzpicture}[
    xscale=1.05, yscale=0.75,
    every node/.style={font=\large},
    line cap=round, line join=round
]
\begin{scope}[rotate=45]

\fill[red!22]
    (0.85,0.32)
    .. controls (0.40,0.32) and (0.32,0.40) .. (0.32,0.85)
    -- (-0.32,0.85)
    .. controls (-0.32,0.40) and (-0.40,0.32) .. (-0.85,0.32)
    -- (-0.85,-0.32)
    .. controls (-0.40,-0.32) and (-0.32,-0.40) .. (-0.32,-0.85)
    -- (0.32,-0.85)
    .. controls (0.32,-0.40) and (0.40,-0.32) .. (0.85,-0.32)
    -- cycle;

\def\BottlePath{
    (0.85,0.32)
    .. controls (1.25,0.30) and (1.55,0.42) .. (1.90,0.85)
    .. controls (2.40,1.65) and (4.20,1.65) .. (4.55,0.60)
    .. controls (4.70,0.22) and (4.70,-0.22) .. (4.55,-0.60)
    .. controls (4.20,-1.65) and (2.40,-1.65) .. (1.90,-0.85)
    .. controls (1.55,-0.42) and (1.25,-0.30) .. (0.85,-0.32)
}

\foreach \angle/\i in {90/1,0/2,180/3,270/4}{
    \begin{scope}[rotate=\angle]
        \fill[red!22] \BottlePath
            -- (0.70,-0.32) -- (0.70,0.32) -- cycle;

        \begin{scope}
            \clip \BottlePath -- cycle;
            \fill[blue!12] (1.55,-2) rectangle (5,2);
            \draw[dashed,thick] (1.55,-2) -- (1.55,2);
        \end{scope}

        \draw[thick] \BottlePath;
        \draw[thick] (0.85,0.32)
            .. controls (0.40,0.32) and (0.32,0.40) .. (0.32,0.85);

        \node at (3.15,0) {$\Omega_{\i}$};
    \end{scope}
}

\node at (0,0) {$\Omega^L$};
\end{scope}
\end{tikzpicture}
\caption{Schematic representation of four bottles
$\Omega_1,\ldots,\Omega_4$ and the central region $\Omega^L$ (red),
 the bottleneck. Dashed lines mark the interfaces
between bottles and $\Omega^L$.}
\label{fig:four-bottles}
\end{figure}
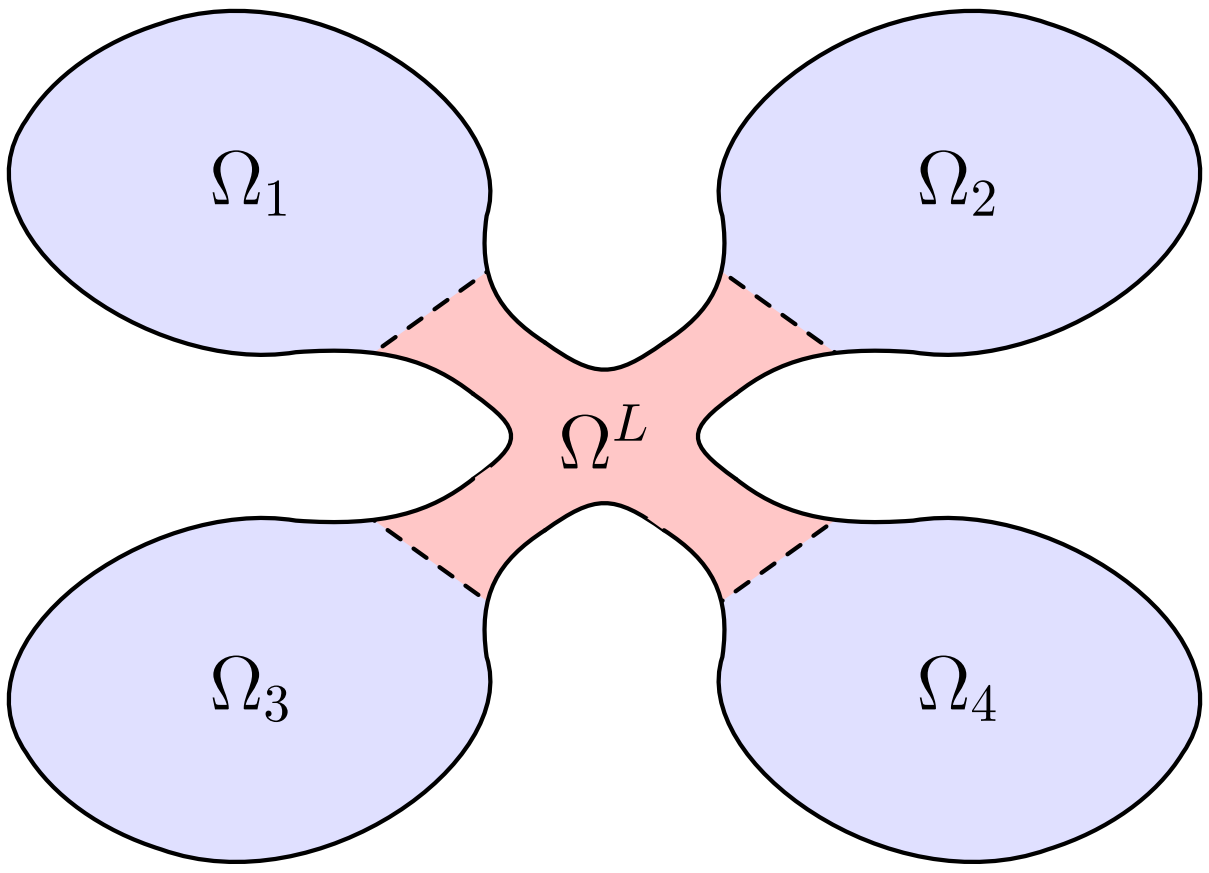
We introduce some new notation. For any two $k$-forms $\tau,\tau'$, we denote by $d_H$ the Hamming distance between them
\[
d_H(\tau,\tau') \coloneqq |\supp (\tau+\tau')|.
\]

Recall the dynamics of the Markov process introduced in Section \ref{subsec:gd}. We denote by $P$ the discrete transition kernel from \eqref{eq:defP}. Then, for such a transition kernel, we define for $A \subseteq C^k$
\[
\partial_PA \coloneqq \{\tau' \in C^k\setminus A \colon \exists \tau \in A:P(\tau,\tau')>0 \}.
\]
This corresponds to the set that can be attained from $A$ using the kernel $P$. 
\begin{Lemma}
\label{lem:partialomega}
    Consider the dynamics from Section \ref{subsec:gd}, using the transition kernel $P$ satisfying \eqref{eq:defP} on the sequence of blocks $(\mathcal{B}_N)_{N \in \N}$ that is assumed to satisfy \eqref{eq:conditionblock}. 
    Then, for any $i \in [\vartheta]$ the following holds
    \[
    \partial_P \Omega_i \subseteq \Omega^L.
    \]
    In other words, the only possibility of leaving a bottle with our dynamics is through the bottleneck $\Omega^L$. 
\end{Lemma}
In particular, Lemma \ref{lem:partialomega} shows that in order to leave $\Omega_i$, a large vortex cluster has to be created as depicted in Figure \ref{fig:changecohomologicalclass}.
\begin{figure}
\centering
\begin{tikzpicture}[
    scale=1.05,
    >=stealth,
    every node/.style={font=\small}
]

\draw[->, thick] (0,0) -- (9,0)
    node[right] {cohomology classes};

\draw[->, thick] (0,0) -- (0,5.8)
    node[above, rotate=90, anchor=south] {Energy of clusters};

\node[below] at (2,0) {$H_i^k$};
\node[below] at (7,0) {$H_j^k$};

\shade[
    bottom color=blue!18,
    top color=blue!6,
    draw=blue!30
]
    (2,1.2) ellipse (1.5 and 1.0);

\shade[
    bottom color=blue!18,
    top color=blue!6,
    draw=blue!30
]
    (7,1.2) ellipse (1.5 and 1.0);

\draw[dashed, very thick, orange!80!black]
    (0.2,4) -- (8.8,4)
    node[right] {$g_\ell(N)$};

\foreach \x/\y in {
    1.20/0.75,
    1.45/1.10,
    1.75/0.55,
    2.00/0.95,
    2.20/0.70,
    2.45/1.15,
    2.75/0.85
}{
    \fill[blue!55!black] (\x,\y) circle (1.8pt);
}

\foreach \x/\y in {
    6.20/0.90,
    6.45/0.65,
    6.75/1.15,
    7.00/0.55,
    7.25/0.85,
    7.50/1.20,
    7.75/0.75
}{
    \fill[blue!55!black] (\x,\y) circle (1.8pt);
}

\node[align=center] at (2,2.15)
    {$\Omega_i$\\ $H_i^k+$ small perturbations};

\node[align=center] at (7,2.15)
    {$\Omega_j$\\ $H_j^k+$ small perturbations};

\draw[->, very thick, red!70!black]
    (2.8,1.35)
    .. controls (3.4,2.5) and (3.7,4.5) ..
    (4.5,4.5)
    .. controls (5.3,4.5) and (5.6,2.5) ..
    (6.2,1.35);

\fill[red!80!black] (4.5,4.5) circle (2.2pt);

\node[above, align=center] at (4.5,4.62)
    {Trajectory from $\Omega_i$\\ to $\Omega_j$};

\end{tikzpicture}
\caption{Schematic trajectory between two bottles,
passing through the bottleneck $\Omega^L$.}
\label{fig:changecohomologicalclass}
\end{figure}
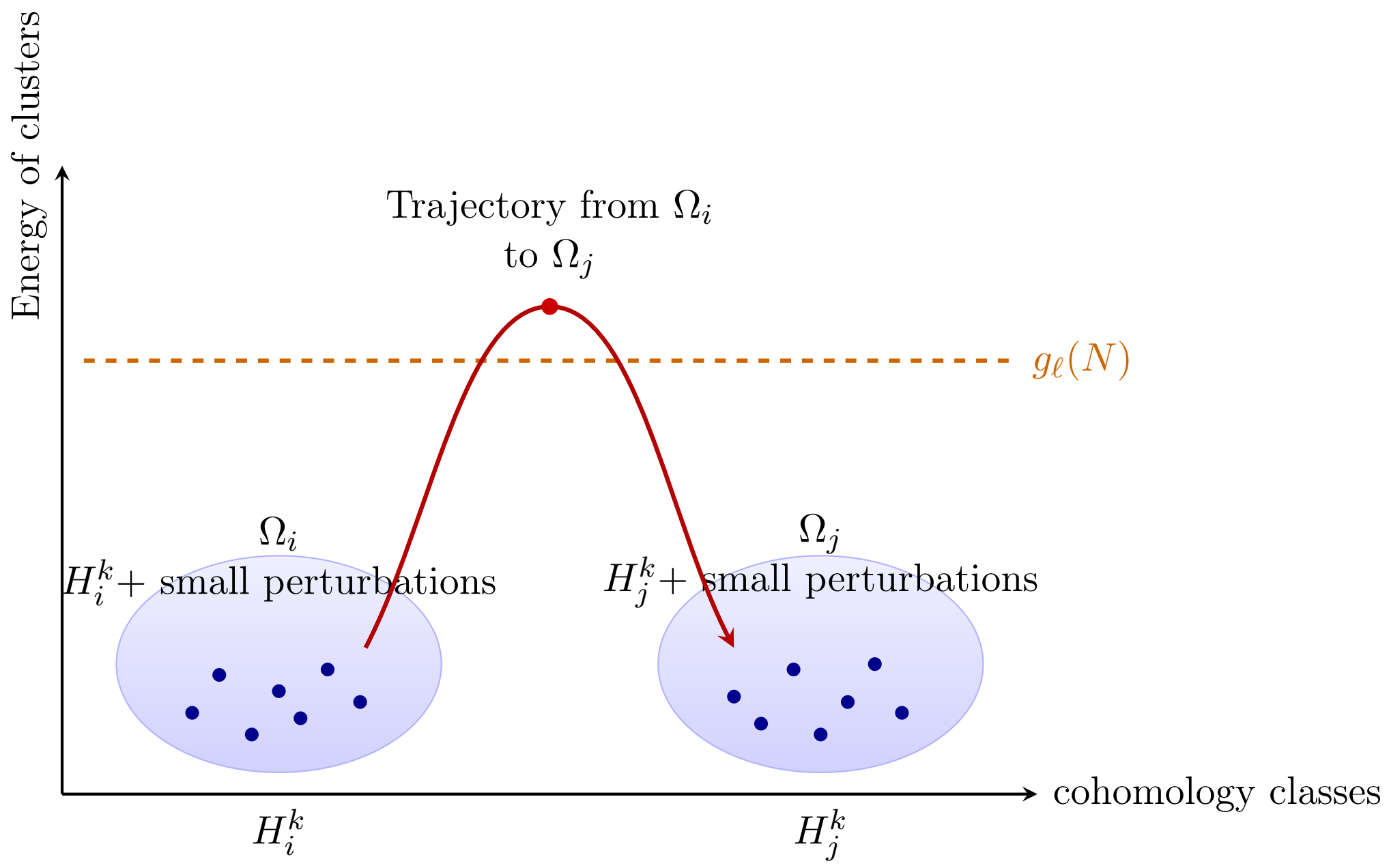
\begin{proof}
    Let $\tau \in \Omega_i$ and let $\tau' \in \partial_P\Omega_i$, such that $P(\tau,\tau')>0$.  Let's moreover denote the decomposition of $d\tau$ and  $d\tau'$ into vortex clusters by
    \[
        d\tau=\sum_{i=1}^{I(d\tau)}\omega_i, \quad
d\tau'=\sum_{j=1}^{I(d\tau')}\tilde{\omega}_j.
    \]
    Now, either there exists $j \in [I(d\tau')]$ such that $|\tilde{\omega}_j| \ge g_\ell(N)$ in which case $\tau' \in \Omega^L$. The other possibility is that all $\tilde{\omega}_j$ have their support bounded above by $g_\ell(N)$ in what case $\tau' \in \Omega$. Our goal is to prove that $\tau' \in \Omega \implies \tau' \in \Omega_i$, which is equivalent to $\Delta \tau$ and $\Delta \tau'$ being in the same cohomology class. By Corollary~\ref{cor:lowerbooundhomologyclasses}, we only have to prove the support condition $d_{H}(\Delta\tau,\Delta\tau')=|\Delta \tau + \Delta \tau'| < N^\ell$ to conclude the proof. 

    We let $\xi\coloneqq \tau+\tau' \in C^k$. Then 
    \begin{align}
    \label{eq:hammingsumdeltatautau'}
    |\Delta\tau+\Delta\tau'|=|\tau+\tau' + \sum_{i=1}^{I(d\tau)}T(\omega_i)+\sum_{j=1}^{I(d\tau')}T(\tilde{\omega}_j)|  \le |\xi|+ |\sum_{i=1}^{I(d\tau)}T(\omega_i) +\sum_{j=1}^{I(d\tau')}T(\tilde{\omega}_j)|.
    \end{align}
    Observe that $\xi$ is supported exactly on the $k$-cells which have been flipped by $P$. Hence the support of $\xi$ has size at most that of the support of the largest block. Hence
    \[
    |\xi| \le \lambda(\mathcal{B}_N). 
    \]
    Let's now focus on the second summand on the right-hand side of \eqref{eq:hammingsumdeltatautau'}. Our first goal is to count how many different vortex clusters can be affected by a flip at a single $k$-cell. Recall that a vortex cluster was defined by a configuration inducing a connected component on the graph $\mathcal{G}_m^{k+1}(N)$ in Definition \ref{Def:vortexclusters}.
    
    Observe that flipping the value of $\tau$ on a single $k$-cell $c \in C_k(\Z_N^m)$ only affects $d\tau$ in the adjacent $(k+1)$-cells $\hat{\partial}c \subset C_{k+1}(\Z_N^m)$. One has that $|\hat{\partial}_kc|=2\ell$: one can pick from $\ell=m-k$ directions to increase the dimension of the $k$-cell, and then one has two sides to choose from.
For $q \in \hat{\partial}_kc$, note that it is connected to at most $2(\ell-1)$ different vortex clusters: one for each upper boundary $Q \in \hat{\partial}_{k+1}q$. To see this, note that for $q$ to become connected (or disconnected) of a vortex cluster, there must be some $q' \sim q$ meaning that $q,q' \in \partial_{k+2}Q$ for some $(k+2)$-dimensional cell $Q$. However, we do not have to count all neighbors of $q$, since the other faces of $Q$ are already counted as part of the same vortex cluster. Hence the number of distinct vortex clusters with which $q$ interacts is at most the number of distinct cells in~$\hat{\partial}_{k+1}q$.
\newline 
Now, observe that each $Q \in \hat{\partial}_{k+1}q$ is counted exactly twice by this method, since it has two different $(k+1)$-faces containing $c$.
 Hence, one deduces that the number of vortex clusters which can be affected by the flip at $c \in C_k(\Z_N^m)$, in each of the configurations before and after the flip is bounded above by 
\[
\frac{2\ell2(\ell-1)}{2}=2\ell(\ell-1).
\]
 Now, at each step of our dynamics we flip at most $\lambda(\mathcal{B}_N)$ $k$-cells. This means that the total number of vortex clusters which are affected by the flips is at most $2\ell(\ell-1) \lambda(\mathcal{B}_N)$. 
    In particular, the sets of indices of the vortex clusters which are different in the two decompositions,
    \begin{align*}
   &A\coloneqq  \{i \in [I(d\tau)] : \nexists j \in [I(d\tau')]~\text{s.t.}~ \omega_i= \tilde{\omega}_j\} \quad \text{and}\\
   &B\coloneqq  \{j \in [I(d\tau')] : \nexists i \in [I(d\tau)]~\text{s.t.}~ \omega_i= \tilde{\omega}_j\},
    \end{align*}
    satisfy $\max\{|A|,|B|\} \le 2\ell(\ell-1) \lambda(\mathcal{B}_N)$. This implies that
    for all but at most $4\ell(\ell-1) \lambda(\mathcal{B}_N)$ vortex clusters it holds that $T(\omega_i)=T(\tilde{\omega_j})$ for some $i \in [I(d\tau)]$, and $j\in [I(d\tau')]$. 
    Hence the Hamming distance between $\Delta\tau$ and $\Delta\tau'$ comes from at most $4\ell(\ell-1) \lambda(\mathcal{B}_N)$ vortex clusters, each of size at most $g_\ell(N)$, as $\tau,\tau'$ are assumed to be in $\Omega$. Together these observations imply that
    \begin{align*}
   & |\sum_{i=1}^{I(d\tau)}T(\omega_i) +\sum_{i=1}^{I(d\tau')}T(\tilde{\omega}_i)| \\ \le &|\{c\in C_k(\Z_N^m): c \in T(\omega_i)~ \text{or}~ c\in T(\tilde{\omega}_j) ~ \text{for} ~ i\in A ~ \text{or} ~ j \in B\}| \\
   \le &4\ell(\ell-1) \lambda(\mathcal{B}_N) \Big(\max_{i \in A} |\supp T(\omega_i)| \vee \max_{j \in B} |\supp T(\tilde{\omega}_j)| \Big)\\  \le &4\hat{C}_{m,k}\ell(\ell-1) \lambda(\mathcal{B}_N)  \Big(\max_{i \in A} |\supp \omega_i| \vee \max_{j \in B} |\supp \tilde{\omega}_j| \Big)^{\frac{\ell}{\ell-1}} \\ < & 4\hat{C}_{m,k}\ell(\ell-1) \lambda(\mathcal{B}_N)g_\ell(N)^{\frac{\ell}{\ell-1}}. 
    \end{align*}
    The first inequality is a consequence of the considerations above as $T$ only acts differently on the vortex clusters $\omega_i$ when $i \in A$ or on $\tilde{\omega}_j$ when $j \in B$. The second inequality counts the number of these clusters and gives an upper bound on their size. The third inequality  comes from the isoperimetric inequality in Lemma~\ref{lem:dec2}, as $T$ maps $\omega_i$ to a $k$-form with minimal support. Finally, the fourth inequality comes from the fact that $\tau,\tau' \in \Omega$, and hence the support of all vortex clusters is bounded above by $g_\ell(N)$. Hence, one obtains overall that
    \[
    d_H(\Delta \tau,\Delta \tau') < \lambda(\mathcal{B}_N)+4\hat{C}_{m,k}\ell(\ell-1) \lambda(\mathcal{B}_N)g_\ell(N)^{\frac{\ell}{\ell-1}} \le N^\ell,
    \]
    by the definition of $g_\ell(N)$ in \eqref{eq:defg_n}. This concludes the proof.
\end{proof}
We now use this result to bound the bottleneck ratio, defined in \eqref{eq:bottleneckratio}.
\begin{Lemma}
\label{lem:bottleneckratiobound}
Consider the dynamics from Section \ref{subsec:gd}, using the transition kernel $P$ satisfying the conditions of Lemma \ref{lem:firstlemmaall}.
    Then, the bottleneck ratio for the Markov process with kernel $P$ satisfies
    \[
    \Phi_P^\star \le \frac{\mu(\partial_P\Omega_i)}{\mu(\Omega_i)} \le \frac{\mu(\Omega^L)}{\mu(\Omega_i)}.
    \]
\end{Lemma}
\begin{proof}
First note that $\mu(\Omega_i) \in (0,\frac{1}{2}]$ by Corollary \ref{cor:probaofasetomegai}.
    By definition of the bottleneck ratio, one has the following
    \begin{align*}
        \Phi^\star \le \Phi(\Omega_i)=\frac{\tilde{P}(\Omega_i,\partial_P\Omega_i)}{\mu(\Omega_i)}  \le \frac{\tilde{P}(\Omega_i,\Omega^L)}{\mu(\Omega_i)},
    \end{align*}
    where $P$ corresponds to the dynamics as described in Lemma \ref{lem:firstlemmaall}, and $\tilde{P}$ is defined in \eqref{eq:edgemeasure}. The last inequality follows from Lemma \ref{lem:partialomega}.
    Using the detailed balance condition \eqref{eq:balancecond}, and recalling the definition of $Q$ in \eqref{eq:defQproposal} one obtains that
    \begin{align*}
        &\frac{\tilde{P}(\Omega_i,\Omega^L)}{\mu(\Omega_i)}=\frac{1}{\mu(\Omega_i)}\sum_{\tau \in \Omega_i, \tau' \in \Omega^L}\mu(\tau)P(\tau,\tau') \\ =&\frac{1}{\mu(\Omega_i)} \sum_{\tau \in \Omega_i, \tau' \in \Omega^L}\mu(\tau)Q(\tau,\tau')\alpha(\tau,\tau') = \frac{1}{\mu(\Omega_i)}\sum_{\tau \in \Omega_i, \tau' \in \Omega^L}\mu(\tau) \frac{\mu(\tau') Q(\tau',\tau)\alpha(\tau',\tau)}{\mu(\tau)} \\
        =& \frac{1}{\mu(\Omega_i)} \sum_{\tau \in \Omega_i, \tau' \in \Omega^L} \mu(\tau') Q(\tau',\tau)\alpha(\tau',\tau) =  \frac{1}{\mu(\Omega_i)} \sum_{\tau' \in \Omega^L} \mu(\tau') \sum_{\tau \in \Omega_i}P(\tau',\tau) \\ \le& \frac{1}{\mu(\Omega_i)} \sum_{\tau' \in \Omega^L} \mu(\tau')  = \frac{\mu(\Omega^L)}{\mu(\Omega_i)}.
    \end{align*}
    Here the inequality comes from $\sum_{\tau \in \Omega_i}P(\tau',\tau) \le 1$ for any $\tau' \in C^k$. This concludes the proof.
\end{proof}
We now introduce the following two estimates, which give an upper bound on the bottleneck ratio. As a first step we give in the following lemma an upper bound for the measure of large vortex clusters.
\begin{Lemma}
    \label{lem:peierlsargument}
    Let $\omega \in \mathcal{V}^{k+1}$ be a vortex cluster as in Definition \ref{Def:vortexclusters}. Then, if 
$\beta$ satisfies the condition \eqref{eq:condbeta} and $(\mathcal{B}_N)_{N \in \N}$ satisfies \eqref{eq:conditionblock},
    there exist constants $C_1,c$, which depend on $m,k,\beta$ and on the entire sequence $(\mathcal{B}_N)_{N \in \N}$ but not on $N$, such that for all~$N \ge N_0$
    \[
    a_N\coloneqq\sum_{\substack{\omega \in \mathcal{V}^{k+1} \\ |\omega| \ge g_\ell(N)}} \myexp(-2\beta|\omega|) \le C_1 \myexp(-cg_\ell(N)).
    \]
\end{Lemma}
\begin{proof}
    Let $\omega' \in \mathcal{V}^{k+1}$ be such that $|\omega|\ge g_\ell(N)$. In order to use a counting argument, we want estimates on the number of such connected $\omega \in \mathcal{V}^{k+1}$ with this support size.
 As $\omega \in \mathcal{V}^{k+1}$,  $\supp \omega$ induces a connected subgraph of $\mathcal{G}_m^{k+1}(N)$ containing exactly $l=|\omega|$~vertices. Note that $\omega$ is uniquely determined by the subgraph of $\mathcal{G}_m^{k+1}(N)$ it induces. Hence, the problem reduces to counting the number of connected subgraphs of $\mathcal{G}_m^{k+1}(N)$ with $l$ vertices. \newline
    By \cite[Lemma 1]{KOMUSIEWICZ2021262}, if for a graph $G=(V,E)$ we denote by $d$ the maximum degree, then the number of connected induced subgraphs of $G$ containing a vertex $v \in V$ and of length $l$ is bounded above by $(e(d-1))^{l-1}$. 
   We apply this result to the graph $\mathcal{G}_m^{k+1}(N)$, which has $\binom{m}{k+1}N^m$ vertices while each vertex has degree $d=2(\ell-1)(2k+3)$ by Lemma~\ref{lem:sizegraph}. We obtain
    \begin{align*}
    a_N&=\sum_{\substack{\omega \in \mathcal{V}^{k+1} \\ |\omega| \ge  g_\ell(N) }} \myexp(-2\beta|\supp \omega|)= \sum_{l=\lceil g_\ell(N) \rceil}^\infty \sum_{\substack{\omega \in \mathcal{V}^{k+1} \\|\omega|=l}} \myexp(-2\beta l)    \\
    &= \sum_{l=\lceil g_\ell(N) \rceil}^\infty e^{-2\beta l} |\{\omega \in \mathcal{V}^{k+1}:|\omega|=l\}| \\
    &\le \sum_{l=\lceil g_\ell(N) \rceil}^\infty e^{-2\beta l} \binom{m}{k+1}N^m\Big(e\big(2(\ell-1)(2k+3)-1\big)\Big)^{l}.
    \end{align*}
    Let $\alpha(k,m,\beta)\coloneqq (2(\ell-1)(2k+3)-1)e^{-2\beta+1}$. Then if $\beta$ satisfies \eqref{eq:condbeta}, it holds that $|\alpha(k,m,\beta)| < 1$ which implies
    \begin{align*}
   & \sum_{l=\lceil g_\ell(N) \rceil}^\infty e^{-2\beta l} \binom{m}{k+1}N^m\Big(e\big(2(\ell-1)(2k+3)-1\big)\Big)^{l}   
   \\=&\sum_{l=\lceil g_\ell(N) \rceil}^\infty \binom{m}{k+1}N^m \alpha(k,m,\beta)^l  
     = \binom{m}{k+1}N^m\frac{\alpha(k,m,\beta)^{ \lceil g_\ell(N) \rceil}}{1-\alpha(k,m,\beta)} \\  \le &\binom{m}{k+1}\frac{N^m}{1-\alpha(k,m,\beta)}\alpha(k,m,\beta)^{g_\ell(N)}.
    \end{align*}
    \textit{Claim.} For $N$ large enough, there exist $C_1,c>0$ independent of $N$ such that if $(\mathcal{B}_N)_{N \in \N}$ satisfies condition \eqref{eq:conditionblock}, then
    \begin{equation}
    \label{eq:expdecayattheend}
    \binom{m}{k+1}\frac{N^m}{1-\alpha(k,m,\beta)} \alpha(k,m,\beta)^{g_\ell(N)} \le C\myexp(-cg_\ell(N)).
    \end{equation}
\textit{Proof of the claim.} Note that \eqref{eq:expdecayattheend} follows if one proves that there is $c >0$ such that for all sufficiently large $N$  the following is satisfied
\begin{equation}
\label{eq:theconditionfortheend}
m\log(N)+g_\ell(N)\log(\alpha(k,m,\beta)) \le -cg_\ell(N).
\end{equation}
Then $C_1>0$ can be chosen later such that this holds for all $N \ge N_0$. Observe that the choice of~$C$ might depend on the block sequence $(\mathcal{B}_N)_{N \in \N}$ as we only control the asymptotics of $(\mathcal{B}_N)_{N \in \N}$ by \eqref{eq:conditionblock}.
Now, note that if $\frac{\log(N)}{g_\ell(N)} \to 0$ then the inequality \eqref{eq:theconditionfortheend} holds for some $c>0$.  
This follows from the definition of $g_\ell(N)$ and from recalling the maximal size of the blocks $\lambda(\mathcal{B}_N)$ in \eqref{eq:conditionblock}. Namely by \eqref{eq:defg_n} one has
\begin{align*}
    &\frac{\log(N)}{g_\ell(N)}=\max \Big(\frac{\log(N)}{N^{\ell-1}}, \log(N)\Big(\frac{4\hat{C}_{m,k}\ell(\ell-1) \lambda(\mathcal{B}_N)}{N^\ell - \lambda(\mathcal{B}_N)} \Big)^{\frac{\ell-1}{\ell}}\Big) \\
    =& \max \Big(\frac{\log(N)}{N^{\ell-1}}, \Big(\frac{4\hat{C}_{m,k}\ell(\ell-1) \log(N)^{\frac{\ell}{\ell-1}} \lambda(\mathcal{B}_N)}{N^\ell - \lambda(\mathcal{B}_N)} \Big)^{\frac{\ell-1}{\ell}}\Big) \to 0.
\end{align*}
As $\ell \ge 2$, the first term goes to $0$ as $N \to \infty$. By \eqref{eq:conditionblock} it holds that $N^\ell-\lambda(\mathcal{B}_N) \sim N^\ell$ and that $\lambda(\mathcal{B}_N)\log(N)^{\frac{\ell}{\ell-1}} = o(N^\ell)$. This concludes the proof of the claim and the lemma.
\end{proof}
\begin{Lemma}
\label{lem:almostthefinalone}
Assume that $\beta$ satisfies \eqref{eq:condbeta} and let $(\mathcal{B}_N)_{N \in \N}$ satisfy \eqref{eq:conditionblock}. Then for $i \in [\vartheta]$ there exist constants $C_2,c>~0$, depending only on $m,k,\beta$ and the sequence of blocks $(\mathcal{B}_N)_{N \in \N}$ such that for all $N \ge N_0$ one has
\[
\frac{\mu_{N}(\Omega^L)}{\mu_{N}(\Omega_i)} \le C_1\myexp\big(-cg_\ell(N)\big).
\]
\end{Lemma}
For the next proof, we introduce the following notation,
\[
\mathcal{V}^{k+1}_{g_\ell(N)} \coloneqq \{\omega \in \mathcal{V}^{k+1} \colon |\omega| \ge g_\ell(N)\}.
\]
\begin{proof}
We have
    \begin{align*}
    \frac{\mu(\Omega^L)}{\mu(\Omega_i)} &= \vartheta\frac{\mu(\Omega^L)}{\mu(\Omega)}= \vartheta \frac{\sum_{\tau \in \Omega^L} \myexp(-2\beta|\supp d\tau|)}{\sum_{\tau \in \Omega} \myexp (-2\beta|\supp d\tau|)} \\ & = \vartheta \frac{|\ker d_k|\sum_{\omega \in \mathcal{W}^L}\myexp(-2\beta|\supp \omega|)}{|\ker d_k|\sum_{\omega \in \mathcal{W}}\myexp(-2\beta|\supp \omega|)} 
    \\ &\le \vartheta \frac{\sum_{\omega \in \mathcal{Z}^L} \myexp(-2\beta|\supp \omega|) \sum_{\omega \in \mathcal{W}} \myexp(-2\beta|\supp \omega|)}{\sum_{\omega \in \mathcal{W}} \myexp (-2\beta|\supp \omega|)} \\ &= \vartheta\sum_{\omega \in \mathcal{Z}^L} \myexp(-2\beta|\supp \omega|).
    \end{align*}
    The first equality comes from Lemma \ref{lem:samemass}; the second equality follows immediately from the definition. The third equality comes from Lemma \ref{lem:thesets}(b), while the inequality is Lemma \ref{lem:thesets}(d). \newline 
    Now, our goal is to prove exponential decay in $N$ of $\sum_{\omega \in \mathcal{Z}^L} \myexp(-2\beta|\supp \omega|)$ when~$\beta$ satisfies condition \eqref{eq:condbeta}. 
    To this end let $\omega =\sum_{i=1}^{I(\omega)} \omega_i \in \mathcal{Z}^L$, and let $A(\omega)=\{\omega_1, \dots,\omega_{I(\omega)}\}$. In particular, $\omega \mapsto A(\omega)$ is injective. Then
\begin{align*}
    &\sum_{\omega \in \mathcal{Z}^L} \myexp(-2\beta|\omega|) =\sum_{\omega \in \mathcal{Z}^L} \prod_{i=1}^{I(\omega)} \myexp(-2\beta|\omega_i|)=\sum_{\omega \in \mathcal{Z}^L} \prod_{\gamma \in A(\omega)}\myexp(-2\beta|\gamma|).
\end{align*}
Now, denote $\mathcal{A} \coloneqq \{A(\omega):\omega \in \mathcal{Z}^L\}$. Then the last sum can be rewritten as
\[
\sum_{\omega \in \mathcal{Z}^L} \prod_{\gamma \in A(\omega)}\myexp(-2\beta|\gamma|)=\sum_{A \in \mathcal{A}}\prod_{\gamma \in A} \myexp(-2\beta|\gamma|) \le \sum_{\emptyset \neq V \subset \mathcal{V}^{k+1}_{g\ell(N)}}\prod_{\gamma \in V} \myexp(-2\beta|\gamma|),
\]
here the inequality comes from the fact that each $A \in \mathcal{A}$ is a non-empty subset of~$\mathcal{V}^{k+1}_{g\ell(N)}$.
Next, one can rewrite the last expression as follows:
\[
\sum_{\emptyset \neq V \subset \mathcal{V}^{k+1}_{g\ell(N)}}\prod_{\gamma \in V} \myexp(-2\beta|\gamma|)=\prod_{\gamma \in \mathcal{V}^{k+1}_{g\ell(N)}}(1+\myexp(-2\beta|\gamma|))-1,
\]
where the equality comes from expanding the product and the $-1$ comes from the fact that we don't sum over the empty set.
Finally we apply the inequality $1+x \le e^x$ to each factor of the product and obtain
\begin{align*}
&\prod_{\gamma \in \mathcal{V}^{k+1}_{g\ell(N)}}(1+\myexp(-2\beta|\gamma|))-1 \le \prod_{\gamma \in \mathcal{V}^{k+1}_{g\ell(N)}} \myexp(e^{-2\beta|\gamma|})-1 \\ 
=&\myexp(\sum_{\gamma \in \mathcal{V}^{k+1}_{g\ell(N)}}e^{-2\beta|\gamma|})-1=\myexp(a_N)-1 \le a_Ne^{a_N}.
\end{align*}
By Lemma \ref{lem:peierlsargument} it holds that $a_N \to 0$, which implies that $a\coloneqq\|a_N\|_\infty$ is finite, and only depends on $\beta,m,k,$ and $(\mathcal{B}_N)_{N \in \N}$. Hence, we obtain
\[
\sum_{\omega \in \mathcal{Z}^L} \myexp(-2\beta|\supp \omega|) \le a_Ne^a.
\]
Using the constants $C_1,c$ from Lemma \ref{lem:peierlsargument}, one obtains the following bound for the bottleneck ratio
\[
\frac{\mu(\Omega^L)}{\mu(\Omega_i)} \le \vartheta e^aC_1\myexp(-cg_{\ell}(N))=:C_2\myexp(-cg_\ell(N)).
\]
This concludes the proof.
    \end{proof}
    \subsection{Proof of Theorem \ref{thm:lowerboundlgt}}
    We are now able to put the pieces together in order to prove Theorem \ref{thm:lowerboundlgt}.
\begin{proof}[Proof of Theorem \ref{thm:lowerboundlgt}]
Assume that $m \in \N_{\ge 2}, k \in \{0,\dots,m-2\}$, that $\beta$ satisfies \eqref{eq:condbeta}, and that the sequence of blocks $(\mathcal{B}_N)_{N \in \N}$ satisfies \eqref{eq:conditionblock}. 
Combining Lemma~\ref{lem:bottleneckratiobound} and Lemma~\ref{lem:almostthefinalone} one obtains the following bound for the bottleneck ratio of the Markov process with transition kernel $P$ for any $i \in [\vartheta]$, where $C_2,c$ are the two constants from Lemma~\ref{lem:almostthefinalone}
    \[
    \Phi^\star_P \le \frac{\mu(\Omega^L)}{\mu(\Omega_i)} \le C_2\myexp\big(-cg_\ell(N)\big).
    \]
From Corollary \ref{cor:continuousdynamics} one deduces that a lower bound for the mixing time under the generator $\mathcal{L}$ from \eqref{eq:defgenL1} is
\begin{equation}
\label{eq:endequationofinalproof}
t_{\text{mix}}^N(\epsilon) \ge \frac{1}{C_2|\mathcal{B}_N|}\Big(\frac{1}{2}-\epsilon\Big)\myexp\big(cg_\ell(N)\big), ~~ \epsilon \in (0,\frac{1}{2}).
\end{equation}
Where $C_1,c>0$ are independent on $N$. The claim now follows from the asymptotical behavior of $g_\ell(N)$ which grows at least as fast as $\Big(\frac{N}{\lambda(\mathcal{B}_N)^\frac{1}{\ell}}\Big)^{\ell-1}$, for large $N$, provided $(\mathcal{B}_N)_{N \in \N}$ satisfies \eqref{eq:conditionblock}. 

This proves the claim for all $N \ge N_0$. For all smaller values of $N \in \{3, \dots, N_0-1\}$, one can pick the constant $C>0$ such that for the constant $c$ from \eqref{eq:endequationofinalproof}, one has
\[
t_{\text{mix}}^N(\epsilon) \ge \frac{C}{|\mathcal{B}_N|}\Big(\frac{1}{2}-\epsilon\Big)\myexp\big(c\Big(\frac{N}{\lambda(\mathcal{B}_N)^\frac{1}{\ell}}\Big)^{\ell-1}\big), ~~ \epsilon \in (0,\frac{1}{2}).
\]
This concludes the proof.
\end{proof}
\begin{Remark}
    As mentioned in Section \ref{subsec:gd}, the dynamics we are considering also include a flip on $(k-1)$-cells. One could believe that including those dynamics would speed up convergence to stationarity (as we add more noise to the model), but it turns out to play no role in the bottleneck argument.
To see this assume that $\tau \in C^k$ is the state of the Markov process, and that $\tau'$ is the updated state, where values on $(k-1)$-cells have been flipped, i.e., 
\[
\tau'=\tau+d_{k-1}\delta_S,~ \delta_S \coloneqq \sum_{c \in S}\delta_c, \quad \text{for some}~  S \subseteq C_{k-1}(\Z_N^m). 
\]
Assume furthermore that $\tau \in \Omega_i$ for some $i \in [\vartheta]$. Then, since $d_k\tau=d_k\tau'$, it holds that $\tau' \in \Omega$, as no new $(k+1)$-defect is created. Now, let $\Gamma \in C_k(\Z_N^m,\Z_2)$ be any homologically non-trivial cycle. Then, note that
\[
\tau'(\Gamma)=\tau(\Gamma)+d_{k-1}\delta_S(\Gamma)\overset{\text{Thm.} \ref{thm:stokes}}{=}\tau(\Gamma)+\delta_S(\partial_k\Gamma) =\tau(\Gamma).
\]
Combined with the fact that $d\tau=d\tau'$ one has
\[
\Delta \tau'(\Gamma)=\tau'(\Gamma)+\sum_{i=1}^{I(d\tau)}T(\omega_i)(\Gamma)=\tau(\Gamma)+\sum_{i=1}^{I(d\tau)}T(\omega_i)(\Gamma)=\Delta \tau(\Gamma)
\]
so that $\Delta\tau'$ and $\Delta\tau$ are in the same cohomology class by Corollary \ref{cor:cohomologyclasses}. In particular, we proved that to go to another cohomology class, one needs to pass through the bottleneck, which is not possible with the dynamics on $(k-1)$-cells.
\end{Remark}
\appendix
\section{Proof of the Isoperimetric Inequality}
\label{sec:app}
In this section, we prove an isoperimetric inequality on $k$-cycles, which yields Lemma~\ref{lem:dec2} by duality. Despite an extensive search, we were not able to find an isoperimetric inequality in the form of this proposition in the literature. The following proof was initiated by AI and rewritten by the author. The principal idea is to use a random $(k+1)$-filling $Q$ of $q$ that is defined recursively. We will prove bounds on the expected support size of $Q$ and then deduce that there exists a realization whose support size is bounded above by the expectation, yielding the desired bound. In a first step we prove a linear bound with a constant depending on $N$, before adapting the construction to obtain a bound independent of $N$.
\begin{Prop}
    \label{prop:appB}
    Let $k \in \{1, \dots,m-1\}$ and let $q \in C_k(\Z_N^m,\Z_2)$ be any null-homologous $k$-cycle.
    Then there exists a $(k+1)$-chain $Q \in C_{k+1}(\Z_N^m,\Z_2)$ such that $\partial_{k+1}Q=q$ and 
    \[
    |\supp Q|\le \frac{3 \cdot 2^{k+1}}{2^k-1}m|\supp q|^{\frac{k+1}{k}}.
    \]
\end{Prop}
Note that $\frac{3 \cdot 2^{k+1}}{2^k-1}m \le 12m$, so the constant in Proposition~\ref{prop:appB} can be replaced by one independent of $k$. \newline 
We first prove the following Lemma, which gives a linear bound with a constant depending on~$N$.
\begin{Lemma}
    \label{lem:thelastlastone}
        Let $k \in \{1, \dots,m-1\}$ and let $q \in C_k(\Z_N^m,\Z_2)$ be any null-homologous $k$-cycle.
    Then there exists a $(k+1)$-chain $Q \in C_{k+1}(\Z_N^m,\Z_2)$ such that $\partial_{k+1}Q=q$ and 
    \[
    |\supp Q|\le mN|\supp q|.
    \]
\end{Lemma}
\begin{proof}
    If $q=0$, take $Q=0$. From now on, we assume that $q \neq 0$. Let $q\in C_k(\Z_N^m,\Z_2)$ be a null-homologous $k$-cycle. 
Let $c=[x;I]$ be any $k$-cell, for $x\in\Z_N^m$ and $I=\{e_{i_1},\ldots,e_{i_k}\}$ a set of coordinate vectors. Let
\[
a_1,\ldots,a_m \overset{\text{i.i.d.}}{\sim} \text{Unif}(\{0, \dots, N-1\}).
\]
For each $j\in\{1,\ldots,m\}$, we let $d_j(x)\in\{0,\ldots,N-1\}$ be the unique integer satisfying
\[
    x_j=a_j+d_j(x)\pmod N,
\]
where $x_j$ is the $j$-th coordinate of the vertex $x \in \Z_N^m$. Furthermore, let
\[
    x^j=x-d_j(x)e_j \in \Z_N^m,
\]
meaning that $x^j$ agrees with $x$ except that its $j$-th coordinate has been translated to $a_j$. Now, define the following operators on $k$-cells for $j \in [m]$
\[
P_j([x;I]) \coloneqq
\begin{cases}
    [x^j;I], & e_j\notin I,\\[2pt]
    0, & e_j\in I,\ d_j(x)<N-1,\\[2pt]
    \displaystyle\sum_{n=0}^{N-1}[x^j+ne_j;I],
        & e_j\in I,\ d_j(x)=N-1,
\end{cases}
\]
and
\[
H_j([x;I]) \coloneqq
\begin{cases}
    \displaystyle\sum_{n=0}^{d_j(x)-1}
        [x^j+ne_j;I\cup\{e_j\}], & e_j\notin I,\\[4pt]
    0, & e_j\in I.
\end{cases}
\]
Here, by convention, an empty sum is zero. Recall that these operators depend on the random sequence $(a_j)_{j=1}^m$, and are hence random as well.
The operators are defined on cells of every dimension $k \in \{0,\dots,m\}$ and get extended linearly to $k$-chains. Note that when $k=m$, we set $H_jc=0$.

Observe that $P_j$ either shifts a $k$-cell to the random coordinate $a_j$, maps it to $0$, or extends it into a $k$-chain around the torus in direction $j$. On the other hand, $H_j$ either suppresses a $k$-cell or maps it to a $(k+1)$-chain extending from $x^j$ to $x-e_j$, with spanning direction $I \cup \{e_j\}$. Now, we observe the following results about these operators. Recall that $\partial_{m+1}=0$ by convention.

\textit{Claim 1.}
For any $k \in \{0,\dots,m\}$ and any $j \in [m]$ one has $H_j\partial_k+\partial_{k+1}H_j=\mathrm{Id}+P_j$. \newline 
\textit{Proof of Claim 1.} We prove the claim for cells, and then extend it to $k$-chains by linearity. Let $c=[x;I]$, and assume first that $e_j\notin I$, so that by definition one obtains
\begin{align}
&\partial_{k+1}H_jc=\sum_{n=0}^{d_j(x)-1}[x^j+(n+1)e_{j};I]+[x^j+ne_{j};I] \notag \\ &+\sum_{e_i \in I}\sum_{n=0}^{d_j(x)-1} [x^j+ne_j+e_i;I\setminus \{e_i\}\cup\{e_j\}]+[x^j+ne_j;I\setminus \{e_i\}\cup\{e_j\}]. \label{eq:eqoperators}
\end{align}
Note that the first sum telescopes; hence one deduces that
\[
\sum_{n=0}^{d_j(x)-1}[x^j+(n+1)e_{j};I]+[x^j+ne_{j};I]=[x^j+d_j(x)e_j;I]+[x^j;I]=c+P_jc.
\]
On the other hand, since $e_j\notin I$, for every $e_i\in I$
we have
\[
d_j(x+e_i)=d_j(x),
\qquad
(x+e_i)^j=x^j+e_i.
\]
One deduces by definition that $H_j\partial_kc$ equals the second line of \eqref{eq:eqoperators}. These terms cancel in $\partial_{k+1}H_jc+H_j\partial_kc$, since the coefficients lie in $\mathbb{Z}_2$, leaving $c+P_jc$, which proves the claim in the case $e_j \notin I$.

Now, assume that $e_j \in I$. Then $H_j$, applied to each cell such that $e_j \in I$, vanishes:
\begin{equation}
\label{eq:j'enaimarredesequationsjeveuxfinir}
\partial_{k+1}H_jc+H_j\partial_kc=H_j([x+e_j;I\setminus\{e_j\}]+[x;I\setminus\{e_j\}]).
\end{equation}
Assume that $d_j(x) <N-1$. Then $P_jc=0$. This also means $d_j(x+e_j)=d_j(x)+1<N$ and $(x+e_j)^j=x^j$; hence the sum telescopes and
\begin{align*}
&H_j([x+e_j;I\setminus\{e_j\}]+[x;I\setminus\{e_j\}])\\ 
=&\sum_{n=0}^{d_j(x)} [x^j+ne_j;I]+\sum_{n=0}^{d_j(x)-1} [x^j+ne_j;I]=[x^j+d_j(x)e_j;I]=[x;I]=c,
\end{align*}
which concludes this case. \newline 
Finally, assume that $d_j(x)=N-1$, meaning that $x^j=x-(N-1)e_j=x+e_j$ as addition is taken modulo $N$.  
Recall that as $e_j \in I$, \eqref{eq:j'enaimarredesequationsjeveuxfinir} holds. 
As $x+e_j=x^j$ it holds that $d_j(x+e_j)=0$ and hence $H_j([x+e_j;I\setminus\{e_j\}])=0$. This leads to
\[
H_j([x+e_j;I\setminus\{e_j\}]+[x;I\setminus\{e_j\}])= \sum_{n=0}^{N-2}[x^j+ne_j;I]=\sum_{n=0}^{N-1}[x^j+ne_j;I]+[x;I]=P_jc+c.
\]
By linearity one extends the claim to $k$-chains. 

\textit{Claim 2.} For any $j \in [m]$ and $k \in [m]$ it holds that $\partial_kP_j=P_j\partial_k$. \newline 
\textit{Proof of Claim 2.}
Observe that the boundary operator applied twice vanishes (see Section \ref{sec:dec}). Hence, as Claim 1 holds for any $j \in [m]$ and any dimension $k \in \{0, \dots,m\}$ one deduces that
\[
\partial_kH_j\partial_k=\partial_kH_j\partial_k+\partial_k\partial_{k+1}H_j=\partial_k+\partial_kP_j,
\]
and that for $k \in \{1,\dots,m\}$
\[
\partial_{k}H_j\partial_k=H_j\partial_{k-1}\partial_k+\partial_{k}H_j\partial_k=\partial_k+P_j\partial_k.
\]
This concludes the proof of Claim 2. 

\textit{Claim 3.} For any chain $q \in C_k(\Z_N^m,\Z_2)$ one has $|H_jq| \le N|q|$ and $\E[|P_jq|] \le |q|$. \newline 
\textit{Proof of Claim 3.}
The first inequality follows immediately from the definition. For the second inequality, we first consider a cell $c=[x;I] \in C_k(\Z_N^m)$. If $e_j \notin I$, then by definition $|P_jc|=|c|=1$. If $e_j \in I$, recall that we picked $a_j \sim \text{Unif}\{0,\dots,N-1\}$. The image $P_jc$ is a $k$-chain with support size $N$ if and only if $a_j=x_j+1 \pmod{N}$, an event of probability $\frac{1}{N}$; otherwise, $P_jc=0$. Hence $
\E[|P_jc|]=N\frac{1}{N}=1=|c|$. 
Finally, we extend these estimates to $k$-chains using the triangle inequality, which proves the claim.

Now, we recursively define the following sequence of $k$-chains for $j \in [m]$
\[
q_0\coloneqq q ~~ \text{and} ~~ q_j \coloneqq P_jq_{j-1}. 
\]
\textit{Claim 4.} Let $k \in [m-1]$. Then $q_m=0$. \newline 
\textit{Proof of Claim 4.} First, note that using Claim 2, each $q_j$ is a $k$-cycle. This can be seen recursively, e.g., $\partial_kq_1=\partial_kP_1q_0=P_1\partial_kq_0=0$.
Now, note that for a cell $c=[x;I]$ of any dimension, applying $P_m \cdots P_1$ to that cell can either suppress the cell, or the cell doesn't get suppressed, meaning that for each $e_j \in I$, one has $d_j(x)=N-1$. In this case the cell gets mapped into a product of full coordinate cycles around the torus, which has zero boundary. In either case, one has $\partial_kP_m \cdots P_1 c =0$.

As $q$ is null-homologous, there exists $Q \in C_{k+1}(\Z_N^m,\Z_2)$ such that $\partial Q=q$. Using Claim 2 and the observation above, one obtains
\[
q_m=P_m \dots P_1q=P_m \cdots P_1\partial_{k+1}Q=\partial_{k+1}P_m \cdots P_1Q=\sum_{c \in Q}\partial_{k+1}P_m \cdots P_1c=0.
\]
This proves the claim.
Now, one defines the (random) $(k+1)$-chain
\[
B\coloneqq\sum_{j=1}^mH_jq_{j-1}.
\]
As each $q_j$ is a cycle, one gets from Claims 1 and 4 that every realization of $B$ is a filling of $q$:
\[
\partial_{k+1}B=\sum_{j=1}^m\partial_{k+1}H_jq_{j-1}=\sum_{j=1}^m(\text{Id}+P_j)q_{j-1}=\sum_{j=1}^m q_{j-1}+q_j=q_0+q_m=q.
\]
Now, observe that for $j \in [m]$, $P_j$ acts only in the $j$-th direction and leaves the other directions unchanged. Using Claim 3 and the independence of $(a_i)_{i=1}^m$, which implies in particular that $q_{j-1}$ is independent of $a_j$, we get
\begin{equation}
\label{eq:supportofq_j}
\E[|q_j|~|~a_1,...,a_{j-1}]=\E[|P_jq_{j-1}|~ |~a_1,...,a_{j-1}]\le |q_{j-1}|.
\end{equation}
This implies $\E[|q_{j}|] \le \E[|q_{j-1}|]$ and hence $\E[|q_j|] \le |q|$ for each $j \in [m]$. From that, one deduces using Claim 3 and the preceding expectation bound, we obtain
\begin{align}
\label{eq:appattheend}
    \E[|B|] \le \sum_{j=1}^m\E[|H_jq_{j-1}|] \le N\sum_{j=1}^m\E[|q_{j-1}|] \le Nm|q|.
\end{align}
This concludes the proof of the preliminary bound, as there is a realization of $B$ denoted by $Q$ that satisfies $|Q| \le Nm|q|$. 
\end{proof}
\begin{proof}[Proof of Proposition \ref{prop:appB}]
If $q=0$ take $Q=0$. Let $0 \neq q \in C_k(\Z_N^m,\Z_2)$ satisfy the assumptions of Proposition \ref{prop:appB}. First, observe that whenever $|q|$ is large enough with respect to the size of the torus, then the bound in Lemma~\ref{lem:thelastlastone} is sharper than the one in Proposition~\ref{prop:appB}. To see this, suppose that $|q|>\Big(\frac{N}{6}\Big)^k$. Then $N < 6|q|^\frac{1}{k}$ and hence
\[
|Q| \le Nm|q|<6m|q|^\frac{k+1}{k} < \frac{2^{k}}{2^k-1}6m|q|^\frac{k+1}{k},
\]
and we have the claim in this case. Now, we have to prove the claim when $|q| \le \Big(\frac{N}{6}\Big)^k$. 
The proof will be quite similar to the above proof, with the main difference that we will partition each coordinate direction into multiple blocks, their number depending on the support of $q$. First we introduce the quantities
    \[ 
    r\coloneqq\lceil 2|q|^{\frac{1}{k}}\rceil \le 3|q|^{\frac{1}{k}} \le \frac{N}{2}, \quad \text{and}~~ w \coloneqq \lfloor \frac{N}{r} \rfloor \ge 2.
    \]
    Then, define 
    \[
    l_u \coloneqq \begin{cases}
        r, ~~ \text{if} ~ u<w, \\
        N-(w-1)r, ~ \text{if} ~ u=w.
    \end{cases}
    \]
As above, let $a_1, \dots, a_m \sim \text{Unif}(\{0,\dots,N-1\})$ be independent coordinates. Fix a direction $j \in [m]$ and partition the coordinates into blocks $B_j^u$, defined as 
\[
B_j^u \coloneqq \{a_j^u+t: t \in \{0,\dots,l_u-1\}\}, \quad a_j^u \coloneqq a_j+(u-1)r \pmod{N}. 
\]
Let $x \in \Z_N^m$. We denote by $x_j$ its $j$-th coordinate and by $u(x) \in \{1,\dots,w\}^m$ the vector of indices of the blocks containing the coordinates of $x$, so that for every $j \in [m]$, $x_j \in B_j^{u_j(x)}$. 
Let $d_j(x)\in\{0,\ldots,l_{u_j(x)}-1\}$ be the unique integer satisfying
\[
x_j=a_j^{u_j(x)}+d_j(x)\pmod N,
\]
and set $x^j \coloneqq x-d_j(x)e_j$. 
Hence, the structure is very similar to the one defined in the proof above, with the main difference that we have cut each direction at $w$ (random) locations $(a_j^u)_{u=1}^w$, instead of one, with block lengths $l_u \in \{r, \dots,2r-1\}$ for $u=1,\dots,w$. 

Now, consider a $k$-cell, $c\coloneqq[x;I]$. One defines the (new) operators $P_j$ and $H_j$ as follows.
\[
P_j([x;I]) \coloneqq
\begin{cases}
    [x^j;I],
        & e_j\notin I,\\[4pt]
    0,
        & e_j\in I,\ d_j(x)<l_{u_j(x)}-1,\\[6pt]
    \displaystyle\sum_{n=0}^{l_{u_j(x)}-1}
        [x^j+ne_j;I],
        & e_j\in I,\ d_j(x)=l_{u_j(x)}-1,
\end{cases}
\]
and
\[
H_j([x;I]) \coloneqq
\begin{cases}
    \displaystyle\sum_{n=0}^{d_j(x)-1}
        [x^j+ne_j;I\cup\{e_j\}],
        & e_j\notin I,\\[6pt]
    0,
        & e_j\in I.
\end{cases}
\]
We see the similarity with the previously defined operators $P_j,H_j$. Here the difference is that the operators act only on the block containing $x$ and do not produce full coordinate cycles. The base vertices of the cells in $P_jc$ have their $j$-th coordinate in $B_j^{u_j}$.

Following the same argument as in the proofs of Claims 1 and 2 in the above proof, for $j\in[m]$ and $k \in \{0,\dots,m\}$ for the first identity and $k \in [m]$ for the second, one has
\begin{equation}
    \partial_{k+1}H_j+H_j\partial_{k}=Id+P_j ~~ \text{and} ~ \partial_{k}P_j=P_j\partial_{k}.
\end{equation}
Now, following the same proof as in Claim 3, one deduces for a $k$-chain $q \in C_k(\Z_N^m,\Z_2)$ that
\begin{equation}
\label{eq:AppsuppHj}
|H_jq| \le 2r|q|,~ \text{and}~~ \E[|P_jq|] \le |q|.
\end{equation}
The first inequality comes from $l_u \le 2r$ and for the second inequality, observe that when $e_j \in I$ for a cell $c=[x;I]$, then
\begin{align*}
&|P_jc|=\sum_{u=1}^w l_u\mathbbm{1}\{x_j=a_j^u+l_u-1\} \implies \\
&\E[|P_jc|]=\sum_{u=1}^w l_u\Prob(x_j=a_j^u+l_u-1)=\frac{1}{N}\sum_{u=1}^wl_u=1.
\end{align*}
The triangle inequality extends this estimate to $k$-chains.
\newline
Now, we define the following sequence of $k$-chains. Let
\[
q_0 \coloneqq q, ~~q_j \coloneqq P_jq_{j-1}.
\]
Observe that unlike above, $q_m$ need not vanish as $\partial_kP_m \cdots P_1c$ is not necessarily $0$. In particular, if $q_m$ does not vanish it must satisfy the following support condition.
\textit{Claim 5.} If $q_m \neq 0$, then $|q_m| \ge r^k$. \newline 
\textit{Proof of Claim 5.} Let $u \in \{1,\dots,w\}^m$ denote the indices of the blocks over the torus and $I$ be a set of $k$ coordinate vectors. Let $\Gamma_{u,I}$ be a $k$-chain such that 
\[
\supp \Gamma_{u,I} \coloneqq \{[y;I]: y_j=a_j^{u_j}, e_j \notin I, ~ y_j \in B_j^{u_j}, e_j \in I\}.
\]
This means that whenever $e_j \notin I$, the coordinate $y_j$ is fixed at the start of the block $B_j^{u_j}$, and whenever $e_j \in I$, $y_j$ ranges over the entire block $B_j^{u_j}$. In particular, as each block has size $l_{u_j} \ge r$, one has $|\Gamma_{u,I}| \ge r^k$. \newline 
Now, let $c=[x;I]$. If $P_m \dots P_1c$ doesn't vanish, then $d_j(x)=l_{u_j(x)}-1$ for all $e_j \in I$. In particular, following the definition of $(P_j)_{j=1}^m$, $c$ gets mapped to the chain $\Gamma_{u(x),I}$, where $u(x)$ contains the block indices in which $x$ lies. This means that $|P_m\cdots P_1c| = |\Gamma_{u(x),I}| \ge r^k$. \newline 
Now, observe that $\supp \Gamma_{u,I} \cap \supp \Gamma_{v,J} \neq \emptyset$ if and only if $I=J$ and $u=v$. The first point is clear by definition: a cell belonging to both supports must have the same set of coordinate directions in each, so $I=J$. Let $[y;I] \in \supp \Gamma_{u,I} \cap \supp \Gamma_{v,I}$. Then for all $e_j \notin I$ one has that the $j$-th coordinate $y_j=a_j^{u_j}$ is fixed. Moreover, for $e_j \in I$, one has $y_j \in B_j^{u_j} \cap B_j^{v_j} \implies u_j=v_j$ as the blocks are disjoint. One obtains that $\Gamma_{u,I},\Gamma_{v,J}$  have disjoint support unless $u=v$ and $I=J$. Now, we can rewrite
\[
q_m=P_m\cdots P_1q=\sum_{c\in \supp q}P_m\cdots P_1c=\sum_{u,I}\alpha_{u,I}\Gamma_{u,I}, \quad \alpha_{u,I} \in \Z_2.
\]
Then $q_m \neq 0$ implies that there is at least one $\alpha_{u,I} =1$. As $(\Gamma_{v,J})_{v,J}$ have disjoint support one deduces that $|q_m|\ge |\Gamma_{u,I}| \ge  r^k$, which concludes the proof of the claim. \newline
Now, for any $j \in [m]$, using the same argument as in \eqref{eq:supportofq_j} one has
\begin{equation}
\label{eq:appsuppofqi}
\E[|q_j|] \le |q|.
\end{equation}
From  Claim 5 and \eqref{eq:appsuppofqi} one deduces that
\[
\Prob(q_m \neq 0) \le \frac{|q|}{r^k} \le \frac{1}{2^k},
\]
where the second inequality follows from $2^k|q| \le r^k$. In particular, 
\begin{equation}
\label{eq:boundonq_m}
\Prob(q_m=0) \ge 1-\frac{1}{2^k}=\frac{2^k-1}{2^k}.
\end{equation}
Now, define the (random) $(k+1)$-chain
\[
B\coloneqq\sum_{j=1}^mH_jq_{j-1}.
\]
Observe that since $\partial_k$ and $P_j$ commute, each $q_j$ is a $k$-cycle. It implies that $H_j\partial_kq_{j-1}=0$ and hence
\[
\partial_{k+1}B=\sum_{j=1}^m\partial_{k+1}H_jq_{j-1}=\sum_{j=1}^m(\text{Id}+P_j)q_{j-1}=\sum_{j=1}^m (q_{j-1}+q_j)=q_0+q_m.
\]
Thus, whenever $q_m$ vanishes, $B$ is a filling of $q$. We compute the following upper bound on the expected support size of $B$.
\begin{equation}
\label{eq:A7lookslikeclosetotheend}
\E[|B|] \le \sum_{j=1}^m\E[|H_jq_{j-1}|] \overset{\eqref{eq:AppsuppHj}}{\le} \sum_{j=1}^m 2r\E[|q_{j-1}|] \overset{\eqref{eq:appsuppofqi}}{\le} 2rm|q|.
\end{equation}
 Among the realizations of $B$ for which $q_m=0$, choose one with the smallest support size and denote it by $Q$. Then $Q$ is a filling of $q$ and we can upper bound its support size as follows:
\begin{align*}
&|Q| \le \E[|B||q_m=0] =  \frac{\E[|B|\mathbbm{1}\{q_m=0\}]}{\Prob(q_m=0)} \\
&\le \frac{\E[|B|]}{\Prob(q_m=0)} \le \frac{2^k}{2^k-1}2mr|q| \le \frac{2^k}{2^k-1}6m|q|^{\frac{1}{k}}|q|=\frac{3 \cdot 2^{k+1}}{2^k-1}m|q|^{\frac{k+1}{k}}.
\end{align*}
The first inequality follows from the  choice of $Q$ with minimal support size. The first equality is the definition, while the second inequality follows from $|B|\ge |B|\mathbbm{1}\{q_m=0\}$. The third inequality comes from \eqref{eq:boundonq_m} and \eqref{eq:A7lookslikeclosetotheend}. The fourth inequality follows from $r \le 3|q|^\frac{1}{k}$. \newline 
This concludes the proof.
\end{proof} 

\newpage
\bibliography{References}
\end{document}